\input epsf

\documentclass[12pt]{article}

\usepackage{latexsym,amsmath,amssymb,amsfonts,amscd,multirow}
\usepackage{epsfig,verbatim,epstopdf,graphics}
\usepackage{color,leftidx}
\usepackage{slashbox}
\usepackage{tabu}
\usepackage{subfigure}
\newtheorem{thm}{Theorem}[section]

\newcommand{\bx}{\mathbf{x}}

\newcommand{\bE}{\mathbf{E}}

\newfont{\iams}{msbm9}

\newcommand{\commentbis}[1]{}
\newcommand{\be}{\begin{eqnarray}}
\newcommand{\ee}{\end{eqnarray}}
\newcommand{\beno}{\begin{eqnarray*}}
\newcommand{\eeno}{\end{eqnarray*}}
\newcommand{\barr}[1]{\begin{array}{#1}}
\newcommand{\earr}{\end{array}}

\newcommand{\df}{\partial}

\newcommand{\mE}{{\mathcal E}}
\def\Real{\mathbb R}
\newcommand{\beq}{\begin{equation}}
\newcommand{\eeq}{\end{equation}}
\newcommand{\beqa}{\begin{eqnarray}}
\newcommand{\eeqa}{\end{eqnarray}}

\newcommand{\bnva}{{\bn_{v_\alpha}}}
\newcommand{\ah}{{\alpha, h}}
\newcommand{\Kx}{{K_x}}

\newcommand{\Ka}{{K_\alpha}}
\newcommand{\Kva}{{K_{v_\alpha}}}
\newcommand{\bv}{{\bf v}}
\newcommand{\bB}{{\bf B}}
\newcommand{\bJ}{{\bf J}}

\newcommand{\bU}{{\bf U}}

\newcommand{\bn}{{\bf n}}
\newcommand{\Ox}{{\Omega_x}}
\newcommand{\Ov}{{\Omega_v}}
\newcommand{\Ove}{{\Omega_{v_e}}}
\newcommand{\Ovi}{{\Omega_{v_i}}}
\newcommand{\Ova}{{\Omega_{v_\alpha}}}
\newcommand{\Rn}{{\mathbb{R}^n}}

\newcommand{\na}{{n+1}}
\newcommand{\mT}{{\mathcal T}}
\newcommand{\mG}{{\mathcal G}}
\newcommand{\mS}{{\mathcal S}}
\newcommand{\mU}{{\mathcal U}}
\newcommand{\mW}{{\mathcal W}}
\newcommand{\mZ}{{\mathcal Z}}
\newcounter{nsez}

\title
{  Numerical Study of  the Two-Species Vlasov-Amp\`{e}re System: Energy-Conserving Schemes  and the Current-Driven Ion-Acoustic Instability}

\author{  Yingda Cheng
\thanks{Department of Mathematics, Michigan State University,
East Lansing, MI 48824 U.S.A.
 {\tt ycheng@math.msu.edu}}
  \and
Andrew J. Christlieb
\thanks{
Department of Mathematics and Department of Electrical and Computer Engineering,
Michigan State University,
East Lansing, MI 48824 U.S.A.
{\tt christli@msu.edu}}
   \and
Xinghui Zhong
\thanks{Corresponding author. Department of Mathematics, Michigan State University,
East Lansing, MI 48824 U.S.A.
 {\tt zhongxh@math.msu.edu}}
}

\date{\today}

\begin{document}

%======================================

\maketitle

\begin{abstract}

In this paper, we propose energy-conserving Eulerian solvers for the   two-species Vlasov-Amp\`{e}re (VA) system and apply the methods to simulate  current-driven ion-acoustic instability. The algorithm is generalized from our previous work for the single-species VA system \cite{cheng_va} and Vlasov-Maxwell (VM) system  \cite{zhong_vm}. The main feature of the schemes is their ability to preserve the total particle number and total energy on the fully discrete level regardless of mesh size. Those are desired properties of numerical schemes especially for long time simulations with under-resolved mesh. The conservation is realized by explicit and implicit energy-conserving temporal discretizations, and the discontinuous Galerkin (DG) spatial discretizations.
We benchmarked our algorithms on a test example to check the one-species limit, and the current-driven ion-acoustic instability. To simulate the current-driven ion-acoustic instability, a slight modification for the implicit method is necessary to fully decouple the split equations. This is achieved by a Gauss-Seidel type iteration technique. Numerical results verified the conservation and performance of our methods.
\end{abstract}


\section{Introduction}

In this paper, we propose energy-conserving Eulerian solvers for the   two-species Vlasov-Amp\`{e}re (VA) system and apply the methods to simulate  current-driven ion-acoustic instability. The two-species VA model describes the evolution of the distribution functions for  a single species of electrons and ions under the influence of the self-consistent electric field. Accurate numerical simulation for this system is crucial for the understanding of ion-acoustic waves, ion-acoustic turbulence in fusion plasmas and magnetic reconnection in space plasmas.

In the literature, a class of well-established methods for the Vlasov equation is the particle-in-cell (PIC) methods  \cite{Birdsall_book1991, Hockney_book1981}. In PIC methods, the macro-particles are advanced in a Lagrangian framework, while the field equations are solved on a mesh.
The main advantage of the PIC method is its relatively low cost for high dimensional problems, but it suffers from  statistical noise built intrinsically in those methods. 
Our approach in this paper is to use a grid-based Vlasov solver, which does not have statistical noise and can resolve the low-density regions more accurately. While there are abundant literature on grid-based Vlasov solver for single-species VA or Vlasov-Poisson (VP) system, e.g. \cite{chengknorr_76, Sonnendrucker_1999, Boris_1976, Filbet_PFC_2001, Horne2001182, crouseilles_va, elkina2006new}, there are relatively fewer published works  for the two-species system. In \cite{eliasson2004dynamics, eliasson2004production, eliasson2010numerical}, Fourier transformed methods are used to compute two-species VP system for electron and ion holes. In \cite{labrunie2004numerical}, the hydrodynamic and quasi-neutral limits of the   two-species VP system are studied by the  finite difference WENO method.
In \cite{Petkaki2003, ionacoustic}, the MacCormack method is employed for calculation of anomalous resistivity and the nonlinear evolution of ion-acoustic instability. A detailed study of the comparison of the MacCormack method and a PIC method for anomalous resisitivity in current-driven ion-acoustic waves can be found in \cite{thesision}.

One of the main focus of this paper is to develop  fully discrete energy-conserving methods. The total energy is a nonlinear quantity that depends on the distribution functions of both species as well as the electric field. To achieve energy conservation, special care must be taken to design both the temporal and spatial discretizations. In this work, we generalize
our previous methods for the single-species VA system \cite{cheng_va} and Vlasov-Maxwell (VM) system  \cite{zhong_vm} to the two-species system. The main feature of our method is that it can preserve the total energy on the fully discrete level regardless of mesh size. This is advantageous for long time simulation, guaranteeing no generation of annulation of spurious energy, and avoiding artifacts such as plasma self heating or cooling \cite{cohen1989performance}. Previously, several PIC methods have been proposed to conserve the total energy for the single-species VA system \cite{chen2011energy} and the VM system \cite{Markidis20117037}. Finite difference and DG methods \cite{Filbet_compare_2003, ayuso2012high} were proposed to conserve the total energy of VP systems on the semi-discrete level. Our method is the first Eulerian solver to achieve conservation of total energy and particle number simultaneously for the two-species system. This is done by using the newly developed energy-conserving temporal discretizations \cite{cheng_va, zhong_vm}, and the discontinuous Galerkin (DG) spatial discretizations \cite{Cockburn_2001_RK_DG}.
In fact, the approach in this paper can be easily adapted to multi-species VA or VM systems as well.

For the two-species system, one of the other computational challenges besides conservation is the multiscale nature of the problem. Because ions are much heavier than electrons,  electrons move faster and the temporal scale for electrons is smaller than that of the ions. For efficient calculations, hybrid and multiscale particle codes have been developed \cite{brackbill1985multiple,lipatov2002hybrid}. We mention in particular   the 
 implicit particle methods \cite{denavit1981time, cohen1982implicit} and electron sub-cycling techniques \cite{adam1982electron}. In this paper, we aim at resolving the physical phenomena that happen at the electron time scale. Therefore the typical time step $\Delta t$ satisfies $w_{pe} \Delta t \propto O(1)$, where $w_{pe}$ is the electron plasma frequency. In the velocity space, we 
 take the common approach of  choosing different computational domain for the velocity space of the electrons and ions, and taking larger grids in electrons than ions. This is allowed because the two species are only coupled together through the electric field. We want to remark that in some applications, it would  be natural to follow the slower ion time scales. In those scenarios, the electron equation becomes stiff and  a multiscale temporal solver would be necessary. However, we do not attempt to address this issue in the current paper and leave it to our future work.

The rest of this paper is  organized as follows: in Section \ref{sec:equation}, we describe the equations under consideration.
In Section \ref{sec:method}, we develop our energy-conserving  schemes and discuss their  properties.  The additional term involving the spatial average of the current density will cause the split equation to be globally coupled.  To resolve this issue, a Gauss-Seidel iteration is employed. Section \ref{sec:numerical} is devoted to numerical results, including the test of one-species limit and the simulations of   the current-driven ion-acoustic waves (CDIAW), in which we perform numerical tests  on an ensemble of 100 VA simulations with random phase perturbations to investigate the anomalous resistivity with a reduced mass ratio.
Finally, we conclude with a few remarks in Section \ref{sec:conclusion}.

%\textcolor{red}{reference, unstructured mesh, nodes,  vector notations}

\section{The Two-Species VA System}
\label{sec:equation}
In this section, we describe the two-species VA system and its dimensionless version.
The two-species VA system for electrons and ions is given by
 \begin{subequations}
 \label{tsva}
 \begin{align}
&\displaystyle\partial_t f_\alpha+  \bv \cdot \nabla_{\bx} f_\alpha  +\frac{q_\alpha}{m_\alpha}\bE  \cdot \nabla_\bv f_\alpha= 0~, \quad (\bx, \bv) \in (\Lambda_x, \Real^n)\label{tsva:1}\\
&\partial_t \bE = - \frac{1}{\epsilon_0} (\bJ-\bJ_{ext}) , \quad \bx \in \Lambda_x\label{tsva:2}
\end{align}
\end{subequations}
%In this section, we will introduce the basic models. 
%We consider the evolution of electron probability distribution function in the presence of a uniform background of ions.
%
where $\alpha=e, i$;  $e$ for electrons and  $i$ for ions. $\Lambda_x \subset\Real^n$ is the physical domain. $f_\alpha(\bx,\bv,t)$ is the probability distribution function  of the particle species $\alpha$.  $m_\alpha$ denotes the particle mass of species $\alpha$. $q_i=-q_e=e$  is the magnitude of the electron charge. In \eqref{tsva:2},
$$\bJ=\sum_\alpha q_\alpha \int_{\Real^n} f_\alpha(\bx, \bv, t) \bv d\bv$$
is the total current density of the two species. $\bJ_{ext}$ is the external current that may be generated by the gradients of an external magnetic field $\bJ_{ext}=\nabla_\bx \times \bB_{ext}$.

%In this paper, we are interested in the physical phenomena that happen at the scale of the electrons. Therefore,  the rescaling of \eqref{tsva} is performed according to the typical time and spatial scales of  electrons.
Given that density,  time and space variables
are in units of the background electron number density $n_0$,  the electron plasma period $\omega_{pe}^{-1}=\left(\displaystyle\frac{n_0e^2}{\varepsilon_0m_e}\right)^{-1/2}$ 
 and the electron Debye radius $\lambda_{De}=\left(\displaystyle\frac{\varepsilon_0k_BT_e}{n_0e^2}\right)^{1/2}$, respectively,
the distribution function $f_\alpha$  is scaled by $n_0/{V_{T_e}}$, 
where 
$V_{T_e}=(k_B T_e/m_e)^{1/2}$ is the electron thermal speed, $T_e$ is the electron
temperature; the electric field $E$ and the current density  are scaled by
$k_B T_e/e\lambda_{De}$  and $n_0 e V_{T_e}$, respectively. 
Keeping the same notations $t, \bx, \bv$ and $f_\alpha, \bE$ for the rescaled unknowns and variables, \eqref{tsva} becomes
the dimensionless two-species VA system   
 \begin{subequations}
 \label{tsvar}
 \begin{align}
&\displaystyle\partial_t f_\alpha+  \bv \cdot \nabla_{\bx} f_\alpha  + \mu_\alpha\bE  \cdot \nabla_\bv f_\alpha= 0~,  \quad (\bx, \bv) \in (\Ox, \Real^n), \quad \alpha=e,i \label{tsvar:1}\\
&\partial_t \bE = -(\bJ-\bJ_{ext}), \quad \bx \in \Ox \label{tsvar:3}
\end{align}
\end{subequations}
where $\mu_\alpha=\displaystyle\frac{q_\alpha m_e }{em_\alpha}$, i.e. $\mu_e=-1$, $\mu_i=\frac{m_e}{m_i}$. $\bJ=\bJ_i-\bJ_e,$ with $\bJ_\alpha=\int_{\Real^n} f_\alpha(\bx, \bv, t) \bv d\bv$.

 The two-species VA system \eqref{tsvar} conserves many physical quantities, such as the  total particle number for each species
$\int_\Ox \int_\Rn  f_\alpha \, d\bv d\bx ,\; \alpha=e,i $, the entropy $\int_\Ox \int_\Rn  f_\alpha \ln(f_\alpha) \, d\bv d\bx $ and any integral of functions of $f_\alpha$, as well as   the total energy 
$$TE=\frac{1}{2} \int _\Ox \int_{\Real^n} f_e |\bv|^2 d\bv d\bx +\frac{1}{2\mu_i} \int _\Ox \int_{\Real^n} f_i |\bv|^2 d\bv d\bx +\frac{1}{2}\int_\Ox |\bE|^2 d\bx, $$
if $\int_\Ox \bE \cdot \bJ_{ext} \,d \bx=0$. This is  true when no external current is present, i.e.  $\bJ_{ext}=0$, as well as for the CDIAW discussed  in Section \ref{sec:cdiaw}. In CDIAW, the external current $\bJ_{ext}$ is a constant chosen to balance the internal current such that $\partial \bE_0/\partial t=0$ \cite{ionacoustic}, where $\bE_0$ denotes the spatially averaged electric field. In \eqref{tsvar}, this is equivalent to letting $\bJ_{ext}=\bJ_0=\displaystyle\frac{1}{|\Ox|}\int_\Ox \bJ \,d \bx$, where $\bJ_0$ denotes the spatially averaged current for $\bJ$. With $\bE_0(t=0)=0$, we will get $\bE_0\equiv0$, therefore $\int_\Ox \bE \cdot \bJ_{ext} \,d \bx=0$, and the energy conservation is implied. For simplicity, in the rest of the paper, we only consider the situation of $\bJ_{ext}=0$ or $\bJ_{ext}=\bJ_0$, with $\bE_0(t=0)=0$.
 In particular, we adopt the following notation
 \begin{subequations}
 \label{tsvartwo}
 \begin{align}
&\displaystyle\partial_t f_\alpha+  \bv \cdot \nabla_{\bx} f_\alpha  + \mu_\alpha\bE  \cdot \nabla_\bv f_\alpha= 0~,  \label{tsvartwo:1}\\
&\partial_t \bE = -\bJ \{+\bJ_0\}, \label{tsvartwo:3}
\end{align}
\end{subequations}
to incorporate the discussion of both cases, where the inclusion of $\{+\bJ_0\}$ is for the CDIAW simulations with $\bJ_{ext}=\bJ_0$.

\section{Numerical Methods}
\label{sec:method}

In this section, we  develop energy-conserving  numerical methods for the two-species VA system \eqref{tsvartwo}. Our methods are generalized from the energy-conserving  methods introduced in \cite{cheng_va, zhong_vm} for one-species VA and VM systems. By proper design, the methods  in this paper can achieve similar conservation properties as those in \cite{cheng_va, zhong_vm}, and can be readily adapted to multi-species VA systems.

\subsection{Temporal discretizations}
\label{sec: time}

In \cite{cheng_va, cheng_vp}, second and higher order temporal discretizations are introduced for the one-species systems. The unique features are that those methods are designed to preserve the discrete total energy. For simplicity, in this paper, we will only consider two types of second-order time stepping methods for the two-species system: one being the fully explicit method, and the other one being the fully implicit method with operator splitting.

The explicit method 
is given  as follows
\begin{subequations}
\label{scheme1}
\begin{align}
&\frac{f_\alpha^{n+1/2}-f_\alpha^n}{\Delta t/2}   +\bv \cdot \nabla_{\bx} f_\alpha^n  +\mu_\alpha  \bE^n \cdot \nabla_\bv f_\alpha^n = 0,\quad \alpha=e,i\label{scheme1:1}\\[2mm]
%&\frac{f^{n+1/2}-f^n}{\Delta t/2}   +\bv \cdot \nabla_{\bx} f^n  +  \gamma\ \bE^n \cdot \nabla_\bv f^n = 0~, \label{scheme1:12}\\
&\frac{\bE^{n+1}-\bE^n}{\Delta t}=-\bJ^{n+1/2} \,\{+\bJ^{n+1/2}_0\},  \quad \textrm{where}\, \,\bJ^{n+1/2}= \int_{\Real^n} (f_i^{n+1/2} -f_e^{n+1/2}) \bv d\bv \label{scheme1:2}\\
&\qquad  \qquad \qquad \qquad \qquad \qquad \qquad   \qquad  \qquad \bJ^{n+1/2}_0= \frac{1}{|\Ox| }\int_{\Ox} \bJ^{n+1/2} d\bx\notag\\
&\frac{f_\alpha^{n+1}-f_\alpha^n}{\Delta t}  + \bv \cdot \nabla_{\bx} f_\alpha^{n+1/2}   + \frac{1}{2} \mu_\alpha (\bE^n+\bE^{n+1})  \cdot \nabla_\bv f_\alpha^{n+1/2} = 0,~ \label{scheme1:3}\
%&\frac{f^{n+1}-f^n}{\Delta t}  + \bv \cdot \nabla_{\bx} f^{n+1/2}   +\gamma\   \frac{1}{2}(\bE^n+\bE^{n+1})  \cdot \nabla_\bv f^{n+1/2} = 0~. \label{scheme1:32}
\end{align}
\end{subequations}
and the term $\{+\bJ^{n+1/2}_0\}$ in \eqref{scheme1:2} is for the case of $\bJ_{ext}=\bJ_0$.
Similar to \cite{cheng_va}, we denote the scheme above to be {\bf Scheme-1}$(\Delta t)$, namely, this means $(f_e^\na,f_i^\na, \bE^\na)= \textnormal{\bf Scheme-1}(\Delta t) (f_e^n,\,f_i^n,\, \bE^n)$.

The fully implicit method is based on the energy-conserving operator splitting for \eqref{tsvartwo} as follows:

\[
  \textrm{ (a)} \left\{
  \begin{array}{l}
\partial_t f_\alpha+  \bv \cdot \nabla_{\bx} f_\alpha   = 0 ,\;\alpha=e,i\\
\partial_t \bE = 0 \,,
  \end{array} \right.\quad
  \textrm{(b)} \left\{
  \begin{array}{l}
\partial_t f_\alpha  +  \mu_\alpha \bE  \cdot \nabla_\bv f_\alpha  = 0,\quad \alpha=e,i \\
\partial_t \bE = -\bJ\{+\bJ_0\},
  \end{array} \right.
\]
Both the split equations maintain the same energy conservation as the original system,
$$\frac{d}{dt} (\int_\Ox \int_{\Real^n} f_e|\bv|^2 d\bv d\bx+\frac{1}{\mu_i}\int_\Ox \int_{\Real^n} f_i|\bv|^2 d\bv d\bx +\int_\Ox |\bE|^2 d\bx ) = 0.$$ In particular,
\[
  \textrm{ (a)} \left\{
  \begin{array}{l}
\displaystyle\frac{d}{dt} \int_\Ox \int_{\Real^n} f_\alpha |\bv|^2 d\bv d\bx = 0~, \\[5mm]
\displaystyle\frac{d}{dt}  \int_\Ox |\bE|^2 d\bx = 0 ,
  \end{array} \right.\]
  \[
  \textrm{(b)} \frac{d}{dt} (\int_\Ox \int_{\Real^n} f_e|\bv|^2 d\bv d\bx+\frac{1}{\mu_i}\int_\Ox \int_{\Real^n} f_i|\bv|^2 d\bv d\bx+ \int_\Ox |\bE|^2 d\bx ) = 0 \{+2\int_\Ox \bJ_0\cdot \bE dx\}=0,
\]
where in the last equality we have used the assumption $\bE_0^0=\bE_0(t=0)=0$, therefore $\bE_0(t)\equiv0$ in the case of $\partial_t \bE=-\bJ+\bJ_0$.

Now we will solve each of the subequations by the second-order implicit midpoint method. In particular, we denote the implicit midpoint method for system (a),
\begin{subequations}
\label{schemea}
\begin{align}
&\frac{f_\alpha^{n+1}-f_\alpha^n}{\Delta t} +  \bv \cdot \nabla_{\bx} \frac{f_\alpha^n+f_\alpha^\na}{2}=0,\\[2mm]
&\frac{\bE^{n+1}-\bE^n}{\Delta t} =0,\,
\end{align}
\end{subequations}
 as $\textnormal{\bf Scheme-a}(\Delta t)$. 
Similarly, for Equation (b),  the implicit midpoint method
\begin{subequations}
\label{schemeb}
\begin{align}
&\frac{f_\alpha^{n+1}-f_\alpha^n}{\Delta t}    +   \frac{1}{2} \mu_\alpha (\bE^n+\bE^{n+1})  \cdot \nabla_\bv \frac{f_\alpha^n+f_\alpha^{n+1}}{2} = 0~,\label{schemeb:1}\\[3mm]
&\frac{\bE^{n+1}-\bE^n}{\Delta t}=- \frac{1}{2} (\bJ^{n}+\bJ^{n+1})\, \{+\frac{1}{2} (\bJ_0^{n}+\bJ_0^{n+1})\}, \label{schemeb:2}
\end{align}
\end{subequations}
is denoted as $\textnormal{\bf Scheme-b}(\Delta t)$. Note that here we have abused the notation, and use superscript $n$, $n+1$ to denote the sub steps in computing equations (a), (b)  rather than the whole time step to compute the VA system. 
Finally, we define
$$\textnormal{\bf Scheme-2}(\Delta t):=\textnormal{\bf Scheme-a}(\Delta t/2) \textnormal{\bf Scheme-b}(\Delta t) \textnormal{\bf Scheme-a}(\Delta t/2).$$

Through simple Taylor expansions, we can verify that both schemes  are  second order accurate in time, also they satisfy discrete energy conservation as illustrated in the theorem below.

\begin{thm}
\label{thm:energy}
With periodic boundary conditions in $\Ox$ domain, the schemes above preserve the discrete total energy $TE_n=TE_{n+1}$, where
$$2\,  (TE_n)=\int_\Ox \int_{\Real^n} f_e^n |\bv|^2 d\bv d\bx +\frac{1}{\mu_i} \int_\Ox  \int_{\Real^n} f_i^n |\bv|^2 d\bv d\bx+\int_\Ox |\bE^n|^2 d\bx$$
in {\bf Scheme-1}$(\Delta t)$ and  {\bf Scheme-2}$(\Delta t)$ for both the case of $\bJ_{ext}=0$ and $\bJ_{ext}=\bJ_0$, with $\bE^0_0=0$.
%, and
%\begin{eqnarray*}
%2\,  (TE_n)&=&\int_\Ox \int_\Ov f_e^n |\bv|^2 d\bv d\bx +\frac{1}{\mu_i} \int_\Ov f_i^n |\bv|^2 d\bv d\bx +\int_\Ox \bE^{n+1/2} \cdot \bE^{n-1/2} d\bx\\
%&=&\int_\Ox \int_\Ov f_e^n |\bv|^2 d\bv d\bx +\frac{1}{\mu_i} \int_\Ov f_i^n |\bv|^2 d\bv d\bx+\int_\Ox |\bE^n|^2 d\bx-\frac{\Delta t^2}{4} \int_\Ox |\bJ^n|^2 d\bx
%\end{eqnarray*}
%in {\bf Scheme-2}$(\Delta t)$.
\end{thm}

\emph{Proof.}  When $\bJ_{ext}=0$, the proof  is similar to   \cite{cheng_va} and is omitted.

In the case of $\bJ_{ext}=\bJ_0$, with $\bE_0^0=0$, in {\bf Scheme-1},  we derive
$$
\bE_0^{n+1}=\int_\Ox \bE^{n+1} \, d \bx=\int_\Ox \bE^{n} \, d \bx-\Delta t \int_\Ox (\bJ^{n+1/2}-\bJ_0^{n+1/2} ) \, d \bx=\int_\Ox \bE^{n} \, d \bx=\bE_0^n.
$$
Therefore,  $\bE_0^n=0,\, \forall n$, if $\bE_0^0=0$.

From \eqref{scheme1}, the additional contribution of $\bJ_0$ term to the total energy difference at $t^{n+1}$ compared to $t^n$  is
$$
\Delta t \int_\Ox (\bE^{n+1}+\bE^n)\cdot \bJ_0^{n+1/2} \, d\bx =\Delta t |\Ox|(\bE_0^{n+1}+\bE_0^n)  \cdot \bJ_0^{n+1/2} =0.
$$
Therefore, we established total energy conservation for $\textnormal{\bf Scheme-1}$ in this case. The proof for $ \textnormal{\bf Scheme-2}$ is similar and omitted. 
 $\Box$
 
% From this theorem, we can see that  {\bf Scheme-1} and {\bf Scheme-3} exactly preserve the total energy, while {\bf Scheme-2} achieves near conservation of the total energy, as is the case for symplectic methods for Hamiltonian systems. The numerical energy from {\bf Scheme-2} is a second order modified version of the original total energy. This is natural due to the second order accuracy of the scheme, and this ensures that over the long run, the numerical energy will not deviate much from its actual value.
% 
%

Similar to \cite{cheng_va, zhong_vm},   higher order time discretizations based on {\bf Scheme-2} can be developed. For simplicity and without loss of generality, we do not pursue them in this paper.

\subsection{Fully discrete methods}
\label{sec: fully}
In this section, we will discuss the spatial discretizations and formulate the fully discrete schemes. In particular, we consider two approaches: one being the explicit scheme, the other being the split implicit scheme. Here, we follow our previous work \cite{cheng_va, zhong_vm} and  use discontinuous Galerkin (DG) methods to discretize the $(\bx, \bv)$ variable. The DG methods are shown to have excellent conservation properties, and when applied to one-species VA and VM systems, the methods can be designed to achieve fully discrete energy conservation \cite{cheng_va, zhong_vm}. For the two-species system, the main difference in our schemes compared to \cite{cheng_va} is the inclusion of the additional species, and the inclusion of $\bJ_0$ term. This causes some additional difficulties as outlined in Section \ref{sec:implicit}.

\subsubsection{Preliminaries}

When discretizing the velocity space, it is necessary  to truncate the domain $\bv \in \Real^n$ into a finite computational region. This is a reasonable assumption as long as the computational domain is taken large enough so that the probability distribution functions vanish at the boundary. For the two-species system, due to the intrinsic scale difference between ions and electrons, we shall use different regions for the two species. In particular, we denote $\Ova$, $\alpha=e, i$ to be the truncated velocity domain for electrons and ions. While the typical size of $\Ove \propto O(1) $, since the ion speed is generally slower, the size of $\Ovi $ would be smaller and $\propto O(\frac{V_{T_i}}{V_{T_e}})=O(\sqrt{\frac{T_i/T_e}{m_i/m_e}})$. As for the $\Ox$ domain, without loss of generality, periodic boundary condition is assumed. We remark here that  our methods can be easily adapted to other types of boundary conditions.

Now we are ready to introduce the mesh and underlying piecewise polynomial spaces. We define $\mT_h^x=\{\Kx\}$ and $\mT_\ah^{v}=\{\Kva\}$ be  partitions of $\Ox$ and $\Ova$, $\alpha=e, i$ respectively,  with $\Kx$ and $\Kva$ being  Cartesian elements or simplices. Notice that  we   use the same mesh in $\bx$ domain for the two species due to their coupling in the Amp\`{e}re equation. However, the mesh in $\bv$ domain is  different for the two-species due to the size difference between $\Ove$ and $\Ovi$.

The meshes for the two species are defined as $\mT_\ah=\{\Ka: \Ka=\Kx\times\Kva, \forall \Kx\in\mT_h^x, \forall \Kva\in\mT_\ah^{v}\}$ . Let $\mE_x$ be the set of the edges of $\mT_h^x$ and  $\mE_{v_\alpha}$ be the set of the edges of $\mT_\ah^{v}$;  then the edges of $\mT_\ah$ will be $\mE_{\alpha}=\{\Kx\times e_{v_\alpha}: \forall\Kx\in\mT_h^x, \forall e_{v_\alpha}\in\mE_{v_\alpha}\}\cup \{e_x\times\Kva: \forall e_x\in\mE_x, \forall \Kva\in\mT_\ah^{v}\}$, $\alpha = e, i$. Here we take into account the periodic boundary condition in the $\bx$-direction when defining $\mE_x$ and $\mE_{\alpha}$. 
%Furthermore, $\mE_v=\mE_v^i\cup\mE_v^b$ with $\mE_v^i$ and $\mE_v^b$ being the set of interior and boundary edges of $\mT_h^v$, respectively.

We will make use of  the following discrete spaces: for $\alpha=e, i$,
\begin{subequations}
\begin{align}
\mG_\ah^{k}&=\left\{g\in L^2(\Omega): g|_{K=\Kx\times\Kva}\in P^k(\Kx\times\Kva), \forall \Kx\in\mT_h^x, \forall \Kva\in\mT_\ah^{v} \right\}, \label{eq:sp:f}\\
\mS_\ah^{k}&=\left\{g\in L^2(\Omega): g|_{K=\Kx\times\Kva}\in P^k(\Kx)\times P^k(\Kva), \forall \Kx\in\mT_h^x, \forall \Kva\in\mT_\ah^{v}\right\},\\
\mZ_\ah^k&=\left\{z\in L^2(\Ova): w|_\Kva\in P^k(\Kva), \forall \Kva\in\mT_\ah^{v} \right\}~,
\end{align}
\end{subequations}
and
\begin{subequations}
\begin{align}
\mU_h^k&=\left\{\bU\in [L^2(\Omega_x)]^{d_x}: \bU|_\Kx\in [P^k(\Kx)]^{d_x}, \forall \Kx\in\mT_h^x \right\}~,\\
\mW_h^k&=\left\{w\in L^2(\Omega_x): w|_\Kx\in P^k(\Kx), \forall \Kx\in\mT_h^x \right\}~,
\end{align}
\end{subequations}
where $d_x$ is the number of dimension for the $\bx$ domain,  $P^k(D)$ denotes the set of polynomials of  total degree at most $k$ on $D$.
The discussion about those spaces for Vlasov equations can be found in \cite{cheng_vp, cheng_vm}.

For piecewise  functions defined with respect to $\mT_h^x$ or $\mT_\ah^v$, we further introduce the jumps and averages as follows. For any edge $e=\{K_x^+\cap K_x^-\}\in\mE_x$, with $\bn_x^\pm$ as the outward unit normal to $\partial K_x^\pm$,
$g^\pm=g|_{K_x^\pm}$, and $\bU^\pm=\bU|_{K_x^\pm}$, the jumps across $e$ are defined  as
\begin{equation*}
[g]_.={g^+}{\bn_.^+}+{g^-}{\bn_.^-},\qquad [\bU]_.={\bU^+}\cdot{\bn_.^+}+{\bU^-}\cdot{\bn_.^-}
\end{equation*}
and the averages are
\begin{equation*}
\{g\}_{.}=\frac{1}{2}({g^+}+{g^-}),\qquad \{\bU\}_.=\frac{1}{2}({\bU^+}+{\bU^-}),
\end{equation*}
where $.$ are used to denote $\bx$ or $\bv_\alpha$.

\subsubsection{The explicit method}
\label{sec:explicit}

In this subsection, we will describe the explicit DG methods formulated with time diescretization {\bf Scheme-1}. In particular, we look for $ f_{\alpha,h}^{n+1/2},f_{\alpha,h}^{n+1}\in\mG_\ah^k$, $\alpha = i, e$, 
%$\bE_h^{n+1} \in\mU_h^k$,
such that for any $\psi_{\alpha,1},\psi_{\alpha,2}\in\mG_\ah^k$,
\begin{subequations}
\label{dscheme1}
\begin{align}
&\int_\Ka \frac{f_{\alpha, h}^{n+1/2}-f_{\alpha, h}^n}{\Delta t/2}  \psi_{\alpha,1} d\bx d\bv
- \int_\Ka f_{\alpha, h}^n\bv\cdot\nabla_\bx \psi_{\alpha,1} d\bx d\bv
-\mu_\alpha \int_\Ka f_{\alpha, h}^n \bE_h^n \cdot\nabla_\bv \psi_{\alpha,1} d\bx d\bv\notag\\[3mm]
&\quad+ \int_{\Kva}\int_{\df\Kx} \widehat{\mu_\alpha f_\ah^n \bv\cdot \bn_x} \psi_{\alpha,1} ds_x d\bv + \int_{\Kx} \int_{\df\Kva} \widehat{( \mu_\alpha f_\ah^n \bE_h^n \cdot \bnva)} \psi_{\alpha,1} ds_{v_\alpha} d\bx=0~,\label{dscheme1:1}\\[3mm]
%&\int_\Kx\frac{\bE_h^{n+1}-\bE_h^n}{\Delta t}\cdot \bU_h d \bx=-\int_\Kx \bJ_h^{n+1/2}\cdot \bU_h d \bx,  \quad \textrm{where}\, \,\bJ_h^{n+1/2}= \int_\Ov f_h^{n+1/2} \bv d\bv~,\label{dscheme1:2}\\
& \frac{\bE_h^{n+1}-\bE_h^n}{\Delta t}=-  \bJ_h^{n+1/2}\,\{+  \bJ_{h,0}^{n+1/2}\} \quad \textrm{where}\, \,\bJ_h^{n+1/2}= \int_\Ov( f_{i,h}^{n+1/2} - f_{e,h}^{n+1/2})\bv d\bv~,\label{dscheme1:2}\\[3mm]
&\qquad  \qquad \qquad \qquad \qquad \qquad \qquad   \qquad  \qquad \bJ^{n+1/2}_{h,0}= \frac{1}{|\Ox|}\int_{\Ox} \bJ_h^{n+1/2} d\bx\notag\\[3mm]
& \int_\Ka \frac{f_{\alpha, h}^{n+1}-f_{\alpha, h}^n}{\Delta t}  \psi_{\alpha,2} d\bx d\bv
- \int_\Ka f_{\alpha, h}^{n+1/2}\bv\cdot\nabla_\bx \psi_{\alpha,2} d\bx d\bv
- \frac{1}{2} \mu_\alpha \int_\Ka f_{\alpha, h}^{n+1/2} (\bE_h^n+\bE_h^{n+1}) \cdot\nabla_\bv \psi_{\alpha,2} d\bx d\bv\notag\\[3mm]
&\quad+  \int_{\Kva}\int_{\df\Kx} \widehat{\mu_\alpha f_{\alpha, h}^{n+1/2} \bv\cdot \bn_x} \psi_{\alpha,2} ds_x d\bv + \frac{1}{2} \int_{\Kx} \int_{\df\Kva} \widehat{(\mu_\alpha f_{\alpha, h}^{n+1/2} (\bE_h^n+\bE_h^{n+1}) \cdot \bnva)} \psi_{\alpha,2} ds_{v_\alpha} d\bx=0.\label{dscheme1:3}
\end{align}
\end{subequations}
Here $\bn_x$ and $\bnva$ are outward unit normals of $\df\Kx$ and $\df\Kva$, respectively.  
Following the discussion in \cite{cheng_va}, to deal with filamentation,   we use the dissipative upwind numerical fluxes, i.e.,
\begin{subequations}
\begin{align}
\mu_\alpha\widehat{f_{\alpha, h}^n \bv\cdot \bn_x}:&=\widetilde{\mu_\alpha f_{\alpha, h}^n \bv}\cdot \bn_x=\left(\{\mu_\alpha f_{\alpha, h}^n\bv\}_x+\frac{|\mu_\alpha\bv\cdot\bn_x|}{2}[f_{\alpha, h}^n]_x\right)\cdot\bn_x~,\label{eq:flux:1}\\
\mu_\alpha\widehat{f_{\alpha, h}^n \bE_h^n \cdot \bnva}:&=\widetilde{\mu_\alpha f_{\alpha, h}^n \bE_h^n }\cdot\bnva=\left(\{\mu_\alpha f_{\alpha, h}^n \bE_h^n\}_{v_\alpha}+\frac{|\mu_\alpha\bE_h^n\cdot\bnva|}{2}[f_{\alpha, h}^n]_{v_\alpha}\right)\cdot\bn_{\bv_\alpha}~.
\end{align}
\end{subequations}
%or central flux
%\begin{subequations}
%\begin{align}
%\mu_\alpha\widehat{f_{\alpha, h}^n \bv\cdot \bn_x}:&=\mu_\alpha\widetilde{f_{\alpha, h}^n \bv}\cdot \bn_x= \mu_\alpha\{f_{\alpha, h}^n\bv\}_x \cdot\bn_x~,\label{eq:flux:1}\\
%\mu_\alpha\widehat{f_{\alpha, h}^n \bE_h^n \cdot \bn_\bv}:&=\mu_\alpha\widetilde{f_{\alpha, h}^n \bE_h^n }\cdot\bn_\bv=\mu_\alpha \{f_{\alpha, h}^n \bE_h^n\}_\bv \cdot\bn_\bv~,
%\end{align}
%\end{subequations}
The upwind fluxes in \eqref{dscheme1:3} are defined similarly. 
%It is well known that the upwind flux is more dissipative and the central flux is more dispersive. With the central flux, DG methods for the linear transport equation has sub-optimal order for odd degree polynomials. A numerical comparison of central and upwind fluxes for the VA system is  shown in \cite{cheng_va}. 

\subsubsection{The implicit method}
\label{sec:implicit}

In this subsection, we would like to design  fully discrete implicit DG methods with  $\textnormal{\bf Scheme-2}$. The key idea is to solve each split equation in their respective reduced dimensions. A complete discussion of similar methods for one-species models on general mesh has been included in \cite{cheng_va}. In this paper, the main difficulty for applying such an approach is for the case of $\bJ_{ext}=\bJ_0$,  which causes a coupling in the $\bx$ direction for equation (b).% due to the $\bJ_0$ term. 

For simplicity, below we will describe  the scheme in detail under 1D1V setting.
For one-dimensional problems, we use a mesh that is a tensor product of grids in the $x$ and $v$ directions, and the domain  is  partitioned as follows:
\begin{equation*}
\label{2dcell1} 0=x_{\frac{1}{2}}<x_{\frac{3}{2}}< \ldots
<x_{N_x+\frac{1}{2}}=L , \qquad -V_{c,\alpha}=v_{{\frac{1}{2}},\alpha}<v_{\frac{3}{2},\alpha}<
\ldots <v_{N_{v,\alpha}+\frac{1}{2},\alpha}=V_{c,\alpha}, \quad \alpha=e, i
\end{equation*}
where $V_{c,\alpha}$ is chosen large enough as the cut-off speed for species $\alpha$.   The mesh is defined as
\begin{eqnarray*}
\label{2dcell2} && K_{r,j,\alpha}=[x_{r-\frac{1}{2}},x_{r+\frac{1}{2}}]
\times [v_{j-\frac{1}{2},\alpha},v_{j+\frac{1}{2},\alpha}] , \nonumber \\
&&K_{x,r}=[x_{r-1/2},x_{r+1/2}],
\quad K_{v,j,\alpha}=[v_{j-\frac{1}{2},\alpha},v_{j+\frac{1}{2},\alpha}]\,, \     \quad r=1,\ldots N_x, \quad j=1,\ldots N_{v,\alpha} ,
\end{eqnarray*}
Let  $\Delta x_r= x_{r+1/2}-x_{r-1/2}$, $\Delta v_{j,\alpha}=v_{j+1/2, \alpha}-v_{j-1/2, \alpha}$ be the length of each interval. $x_r^{(l)}, l=1, \ldots, k+1$ be the $(k+1)$ Gauss quadrature points on $K_{x,r}$ and $v_{\alpha, j}^{(m)}, m=1, \ldots, k+1$ be the $(k+1)$ Gauss quadrature points on $K_{v,j,\alpha}$. Now we are ready to describe our scheme.
\bigskip

\underline{Algorithm {\bf Scheme-a}$(\Delta t)$}
\medskip

To solve from $t^n$ to $t^\na$
 \[
  \textrm{ (a)} \left\{
  \begin{array}{l}
\partial_t f_\alpha  +  v\, \partial_x f_\alpha   = 0, \quad \alpha=e, i\\
\partial_t E = 0,
  \end{array} \right.
\]
\begin{enumerate}
\item For each species $\alpha$ and  $j=1,\ldots N_{v,\alpha} , m=1, \ldots, k+1$, we seek $g_{\alpha,j}^{(m)}(x) \in \mW_h^k$, such that
\begin{eqnarray}
 &&   \int_{K_{x,r}} \frac{g_{\alpha,j}^{(m)}(x)-f_{\alpha,h}^n(x, v_{\alpha, j}^{(m)} )}{\Delta t} \varphi_{\alpha,h} \, dx - \int_{K_{x,r}} v_{\alpha, j}^{(m)}\frac{g_{\alpha,j}^{(m)}(x)+f_{\alpha,h}^n(x, v_j^{(m)} )}{2} (\varphi_{\alpha,h})_x \, dx \label{schemeas} \\
 &&+ \widehat{v_{\alpha, j}^{(m)} \frac{g_{\alpha,j}^{(m)}(x_{r+\frac{1}{2}})+f_{\alpha,h}^n(x_{r+\frac{1}{2}}, v_{\alpha, j}^{(m)} )}{2} } (\varphi_{\alpha,h})_{r+\frac{1}{2}}^- \,
    -\widehat{v_{\alpha, j}^{(m)} \frac{g_{\alpha,j}^{(m)}(x_{r-\frac{1}{2}})+f_{\alpha,h}^n(x_{r-\frac{1}{2}}, v_{\alpha, j}^{(m)} )}{2} } (\varphi_{\alpha,h})_{r-\frac{1}{2}}^+ =0 \notag
\end{eqnarray}
holds for any test function $\varphi_{\alpha,h}(x) \in \mW_h^k$, where the flux terms in \eqref{schemeas} are chosen as the upwind flux, similar to Section \ref{sec:explicit}.
\item Let $f_{\alpha,h}^{n+1}$ be the unique polynomial in $\mS_\ah^k$, such that $f_{\alpha,h}^\na(x_r^{(l)}, v_{\alpha, j}^{(m)})=g_{\alpha,j}^{(m)}(x_r^{(l)}), \, \forall r, j, l, m$.
\end{enumerate}

\bigskip

\underline{Algorithm {\bf Scheme-b}$(\Delta t)$}
\medskip

\textbf{Case 1:} $\bJ_{ext}=0$. This case is similar to the discussion of \cite{cheng_va}.

To solve from $t^n$ to $t^\na$
\[
  \textrm{(b)} \left\{
  \begin{array}{l}
\partial_t f_\alpha  +  \mu_\alpha E \,\partial_vf_{\alpha}  = 0, \quad \alpha=e, i \\
\partial_t E = -J,
  \end{array} \right.
\]

\begin{enumerate}
\item For each species $\alpha$, $r=1,\ldots N_x ,  l=1, \ldots, k+1$, we seek $g_{\alpha,r}^{(l)}(v) \in \mZ_\ah^k$ and $E_r^{(l)}$, such that
\begin{eqnarray}
  &\displaystyle\int_{K_{v,j,\alpha}} \frac{g_{\alpha,r}^{(l)}(v)-f_{\alpha,h}^n(x_r^{(l)}, v)}{\Delta t} \varphi_{\alpha,h} \,  dv - \mu_\alpha \int_{K_{v,j,\alpha}} \frac{E_h^n(x_r^{(l)})+E_r^{(l)}}{2}  \frac{g_{\alpha,r}^{(l)}(v)+f_{\alpha,h}^n(x_r^{(l)}, v)}{2}  (\varphi_{\alpha,h})_v \, dv \notag\\[2mm]
   & + \mu_\alpha\displaystyle \frac{E_h^n(x_r^{(l)})+E_r^{(l)}}{2} \widehat{\frac{g_{\alpha,r}^{(l)}(v_{j+\frac{1}{2},\alpha})+f_{\alpha,h}^n(x_r^{(l)}, v_{j+\frac{1}{2}, \alpha})}{2}} (\varphi_{\alpha,h})^-_{ j+\frac{1}{2}} \notag\\[2mm]
    &-\displaystyle \mu_\alpha \frac{E_h^n(x_r^{(l)})+E_r^{(l)}}{2} \widehat{\frac{g_{\alpha,r}^{(l)}(v_{j-\frac{1}{2},\alpha})+f_{\alpha,h}^n(x_r^{(l)}, v_{j-\frac{1}{2},\alpha})}{2}} (\varphi_{\alpha,h})^+_{ j-\frac{1}{2}} =0 \label{scheme3ss}\\[4mm]
   &\displaystyle\frac{E_r^{(l)}-E_h^n(x_r^{(l)})}{\Delta t}=  -\frac{1}{2}(J_h^n(x_r^{(l)})+J_r^{(l)}), \; \notag
\end{eqnarray}
holds for any test function $\varphi_{\alpha,h}(v) \in \mZ_\ah^k$,  where $J_h^n(x)=\int_\Ovi f_{i,h}^n(x, v) v dv- \int_\Ove f_{e,h}^n(x,v) v dv,$ $J_r^{(l)}=\int_\Ovi g_{i,r}^{(l)}v dv- \int_\Ove g_{e,r}^{(l)} v dv$, and   the flux terms in \eqref{scheme3ss} are chosen as the upwind flux, similar to Section \ref{sec:explicit}.
\item Let $f_{\alpha,h}^{n+1}$ be the unique polynomial in $\mS_\ah^k$, such that $f_{\alpha,h}^\na(x_r^{(l)}, v_{\alpha,j}^{(m)})=g_{\alpha,r}^{(l)}(v_{\alpha,j}^{(m)}),\, \forall r, j, l, m$.
Let $E_h^{n+1}$ be the unique polynomial in $\mW_h^k$, such that $E_h^\na(x_r^{(l)})=E_r^{(l)}, \, \forall r, l$.
\end{enumerate}

To solve \eqref{scheme3ss}, a
Jacobian-free Newton-Krylov solver \cite{knoll2004jacobian}  (KINSOL) is necessary. Notice  we need to set a tolerance parameter $\epsilon_{tol}$  in KINSOL, and that may cause some slight deviation in the conservation.

\bigskip

\textbf{Case 2:} $\bJ_{ext}=\bJ_0$  with $\bE_0(t=0)=0$. This case is different due to the $x$ coupling from $J_0$ term.
To solve from $t^n$ to $t^\na$
\[
  \textrm{(b)} \left\{
  \begin{array}{l}
\partial_t f_\alpha  +   \mu_\alpha E \,\partial_vf_{\alpha}  = 0, \quad \alpha=e, i \\
\partial_t E = -J+J_0,
  \end{array} \right.
\]
a direct generalization of  \eqref{scheme3ss} would cause a nonlinearly coupled system in the whole $(x, v)$ space because $J_0$ involves all elements in $x$. To resolve this issue, we employ a simple Gauss-Seidel   iteration as outlined below.
\begin{enumerate}
\item Initialize with $f_\alpha^{n+1,0}=f_\alpha^n$, $\alpha=e, i$.
\item Iterate on $k=0, 1, \ldots$, solve
\begin{subequations}
\label{schemegs}
\begin{align}
&\frac{E^{n+1,k+1}-E^n}{\Delta t}=- \frac{1}{2} (J^{n}+J^{n+1,k})+\frac{1}{2} (J_0^{n}+J_0^{n+1,k}), \label{schemegs:1}\\[2mm]
&\frac{f_\alpha^{n+1,k+1}-f_\alpha^n}{\Delta t}    +  \frac{1}{2}\mu_\alpha(E^n+E^{n+1,k+1})  \cdot \nabla_\bv \frac{f_\alpha^n+f_\alpha^{n+1,k+1}}{2} = 0~,\label{schemegs:2}
\end{align}
\end{subequations}
until convergence, i.e. $||f_\alpha^{n+1,K+1}-f_\alpha^{n+1,K}||_\infty<\epsilon_{tol}$, where $\epsilon_{tol}$ is a preset tolerance parameter.
\item Set $f_\alpha^{n+1}=f_\alpha^{n+1,K+1}$, $E^{n+1}=E^{n+1, K+1}$.
\end{enumerate}
We notice that the \eqref{schemegs:1} can be implemented explicitly, i.e. for $\forall \,r, l$,
$$
\frac{E^{n+1,k+1}_h(x_r^{(l)})-E_h^n(x_r^{(l)})}{\Delta t}=  -\frac{1}{2}(J_h^n(x_r^{(l)})+J_h^{n+1,k}(x_r^{(l)}))+\frac{1}{2}(J_{h,0}^n+J_{h,0}^{n+1,k})
$$
where 
$J_h^n(x)=\int_\Ovi f_{i,h}^n(x, v) dv- \int_\Ove f_{e,h}^n(x,v) v dv,$ $J_h^{n+1,k}(x)=\int_\Ovi f_{i,h}^{n+1,k}(x, v) dv- \int_\Ove f_{e,h}^{n+1,k}(x,v) v dv,$ and $J_{h,0}^n=\frac{1}{L}\int_\Ox J_h^n(x) dx,$ $J_{h,0}^{n+1,k}=\frac{1}{L}\int_\Ox J_h^{n+1,k}(x) dx.$ This would determine uniquely $E_h^{n+1,k+1} \in \mW_h^k$.

The linear systems resulting from \eqref{schemegs:2} can be evaluated at each Gauss quadrature nodes in $x$ direction, i.e. for each species $\alpha$, $\forall \, r, l$, we seek $g_{\alpha,r}^{(l)}(v) \in \mZ_\ah^k$ and $E_r^{(l)}$, such that
\begin{eqnarray}
  &\hspace{-10mm}\displaystyle\int_{K_{v,j,\alpha}} \frac{g_{\alpha,r}^{(l)}(v)-f_{\alpha,h}^n(x_r^{(l)}, v)}{\Delta t} \varphi_{\alpha,h} \,  dv - \mu_\alpha \int_{K_{v,j,\alpha}} \frac{E_h^n(x_r^{(l)})+E^{n+1,k+1}_h(x_r^{(l)})}{2}  \frac{g_{\alpha,r}^{(l)}(v)+f_{\alpha,h}^n(x_r^{(l)}, v)}{2}  (\varphi_{\alpha,h})_v \, dv \notag\\
   &\hspace{-6mm} \displaystyle+\mu_\alpha \frac{E_h^n(x_r^{(l)})+E^{n+1,k+1}_h(x_r^{(l)})}{2} \widehat{\frac{g_{\alpha,r}^{(l)}(v_{j+\frac{1}{2},\alpha})+f_{\alpha,h}^n(x_r^{(l)}, v_{j+\frac{1}{2},\alpha})}{2}} (\varphi_{\alpha,h})^-_{ j+\frac{1}{2}} \notag\\
    &\displaystyle-\mu_\alpha\frac{E_h^n(x_r^{(l)})+E^{n+1,k+1}_h(x_r^{(l)})}{2} \widehat{\frac{g_{\alpha,r}^{(l)}(v_{j-\frac{1}{2},\alpha})+f_{\alpha,h}^n(x_r^{(l)}, v_{j-\frac{1}{2},\alpha})}{2}} (\varphi_{\alpha,h})^+_{ j-\frac{1}{2}} =0 \label{schemegsd}
\end{eqnarray}
holds for any test function $\varphi_{\alpha,h}(v) \in \mZ_\ah^k$, where the flux terms in \eqref{schemegsd} are chosen as the upwind flux, similar to Section \ref{sec:explicit}. Then we let $f_{\alpha,h}^{n+1,k+1}$ be the unique polynomial in $\mS_\ah^k$, such that $f_{\alpha,h}^{n+1,k+1}(x_r^{(l)},  v_{\alpha,j}^{(m)})=g_{\alpha,r}^{(l)}( v_{\alpha,j}^{(m)}),\, \forall \, \alpha, r, l, j, m$.

\bigskip

Finally, we recall
$\textnormal{\bf Scheme-2}(\Delta t)=\textnormal{\bf Scheme-a}(\Delta t/2) \textnormal{\bf Scheme-b}(\Delta t) \textnormal{\bf Scheme-a}(\Delta t/2)$ and this completes the description of the fully implicit method.
%The main difference compared to the one-species case is the coupling of an additional species in \eqref{scheme3ss}.

\subsubsection{Properties of the fully discrete methods}
In this subsection, we summarize the conservation properties of the fully discrete methods for the two-species VA system. The proof of the theorems below is similar to \cite{cheng_va} by utilizing Theorem \ref{thm:energy} and the properties  of Gauss quadrature formulas, and thus is omitted. 

\begin{thm}[{Total particle number conservation}]
\label{thm:mass}
The DG schemes described in Sections \ref{sec:explicit}, \ref{sec:implicit}
  preserve the total particle number of the system, i.e.
$$\int_\Ox \int_\Ova f_{\alpha,h}^\na dv dx =\int_\Ox \int_\Ova f_{\alpha,h}^n  dv dx, \quad \alpha=e, i.$$
\end{thm}

\begin{thm}[{Total energy conservation}]
\label{thm:fullenergy}
If $k \geq 2$, the   DG schemes described in Sections \ref{sec:explicit}, \ref{sec:implicit} preserve the discrete total energy $TE_n=TE_{n+1}$, where
$$2 \,(TE_n)=\int_\Ox \int_\Ove f_{e,h}^n |v|^2 dv dx+\frac{1}{\mu_i}\int_\Ox \int_\Ovi f_{i,h}^n |v|^2 dv dx +\int_\Ox |E_h^n|^2 dx.$$
\end{thm}

\begin{thm}[{$L^2$ stability}]
The fully implicit DG scheme described in Section \ref{sec:implicit} is $L^2$ stable, i.e.
$$\int_\Ox \int_\Ova |f_{\alpha,h}^\na|^2 dv dx \leq \int_\Ox \int_\Ova |f_{\alpha,h}^n|^2  dv dx, \quad \alpha=e, i.$$
\end{thm}

We notice that the theorems above do not take into account the deviation of $f_\alpha$ from zero at $\partial \Ova$ and the tolerance parameters in the implicit solves. Those factors are the only  possible error sources in the numerical computations for particle number and energy conservation.

\section{Numerical Results }
\label{sec:numerical}

In this section, we demonstrate the performance of our methods in the 1D1V setting.  For simplicity, we use uniform meshes in $x$ and $v$ directions, while we note that nonuniform mesh can also be used under the DG framework.
We use quadratic polynomial spaces and test $\textnormal{\bf Scheme-1}$ with space $\mG_\ah^2$, and $\textnormal{\bf Scheme-2}$  with space $\mS_\ah^2$, respectively.

The time step $\Delta t$ is chosen according to
 \begin{align*}
 \Delta t_e = CFL \left(V_{c,e}\frac{N_x}{L}+\mu_e E_{\max}\frac{N_{v,e}}{V_{c,e}}\right)^{-1}, \quad  \Delta t_i = CFL \left(V_{c,i}\frac{N_x}{L}+\mu_i E_{\max}\frac{N_{v,i}}{V_{c,i}}\right)^{-1}
 \end{align*}
 $$\Delta t = \min(\Delta t_e, \Delta t_i).$$
Notice that the explicit scheme has to satisfy the CFL restriction for stability, while the implicit method can allow large CFL numbers in the computation. 

Two numerical examples are considered in this section: a two-species model with $\bJ_{ext}=0$ which has the limit of Landau damping when the mass ratio $\mu_i \rightarrow 0$, and the CDIAW with $\bJ_{ext}=\bJ_0$ and $\bE_0(t=0)=0$.

For  $\textnormal{\bf Scheme-2}$, we use KINSOL from SUNDIALS \cite{hindmarsh2005sundials} to solve the nonlinear algebraic systems \eqref{scheme3ss} and  \eqref{schemegs:2} resulting from the discretization of equation (b), and we set the tolerance number to be $\varepsilon_{tol}=10^{-12}$. 
The tolerance $\epsilon_{tol}$ is set to be $10^{-11}$ in the  Gauss-Seidel  iteration solving the system \eqref{schemegs}.

\subsection{Testing the one-species limit}
In this subsection, we consider the two-species VA system with a fixed temperature ratio $T_e/T_i=2$ and varying mass ratio to test the one-species limit of our methods. In particular, by letting the initial condition be
%test the one-species limit using the Landau damping example, where the initial conditions are given by
\begin{subequations}
    \label{landauinit}
  \begin{align}
      f_{e}(x,v,0)&=\left(1+A \cos(\kappa x)\right)\frac{1}{\sqrt{2\pi}}e^{-v^2/2},\\
        f_{i}(x,v,0)& = \frac{1}{\sqrt{2\pi\gamma}}e^{-v^2/2\gamma},
   \end{align}
   \end{subequations}
   where  $A = 0.5, \; \kappa = 0.5$ and  $\gamma=(T_im_e/T_em_i)^{1/2}$, we will recover the one-species Landau damping in the limit of $\mu_i=m_e/m_i\rightarrow0$. This initial condition  corresponds to ions in a uniform equilibrium state with a slightly perturbed electron distribution. 
   
The computational domain for $x$ is set to be  $[0,\, L]$, with $L = 4\pi$. The domain of velocity $v$ for electron and ion is chosen to be $[-V_{c,e},V_{c,e}]$ with $V_{c,e}=8$ and $[-V_{c,i},V_{c,i}]$ with $V_{c,i}=\sqrt\gamma V_{c,e}$, respectively, such that $f\simeq 0$ on the boundaries. We use a mesh of uniform $N_x$ cells in the $x$ direction, and $N_v=N_{v,e}=N_{v,i}$ cells in the $v$ direction.
Our numerical examples are performed with two sets of mass ratios,  $\mu_i^{-1}=m_i/m_e=25$ (the reduced mass ratio) and $\mu_i^{-1}=m_i/m_e=1836$ (the real mass ratio). The real mass ratio corresponds to heavy ions that are essentially immobile.

\begin{figure}[!htbp]
\centering
\subfigure[$\textnormal{\bf Scheme-1}$. $\mu_i=1/25$.]{\includegraphics[width=0.45\textwidth]{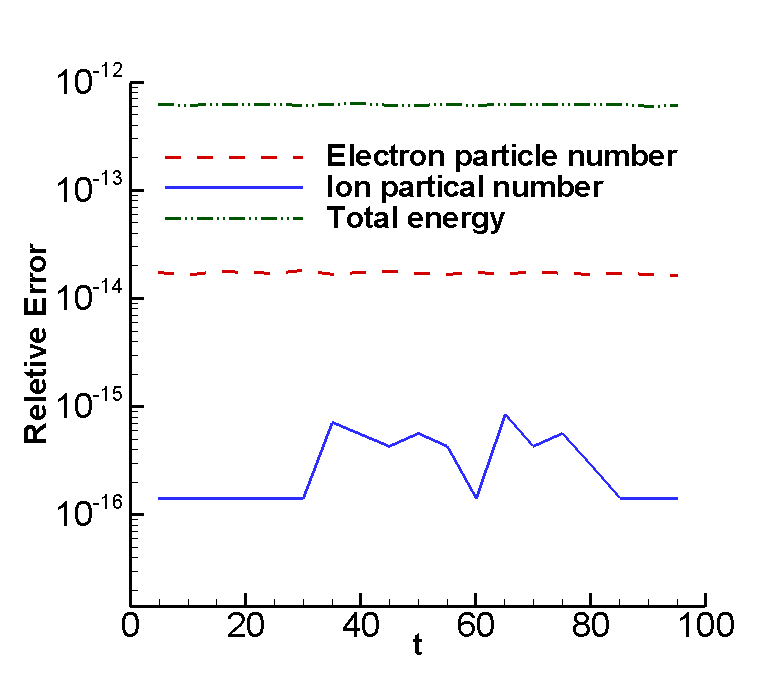}}
\subfigure[$\textnormal{\bf Scheme-1}$. $\mu_i=1/1836$.]{\includegraphics[width=0.45\textwidth]{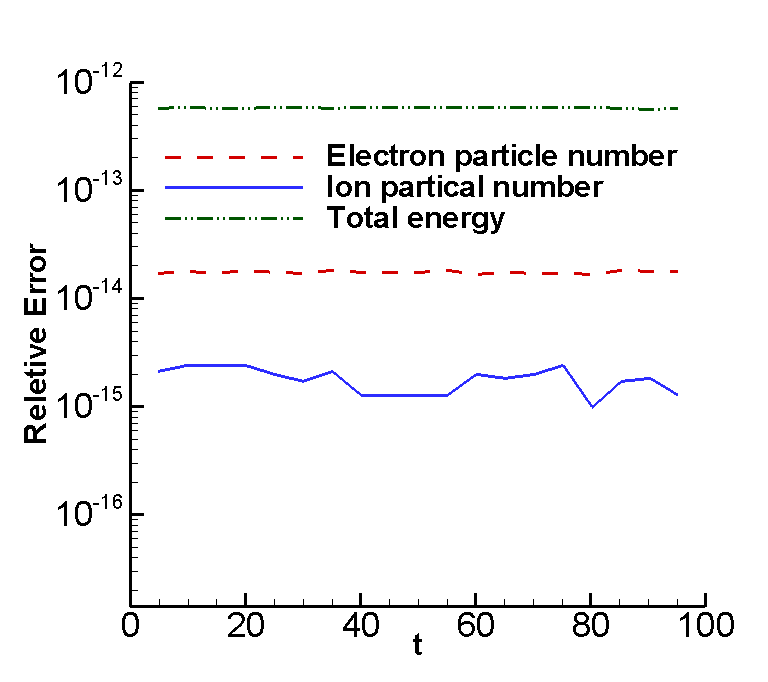}}\\
\subfigure[$\textnormal{\bf Scheme-2}$. $\mu_i=1/25$.]{\includegraphics[width=0.45\textwidth]{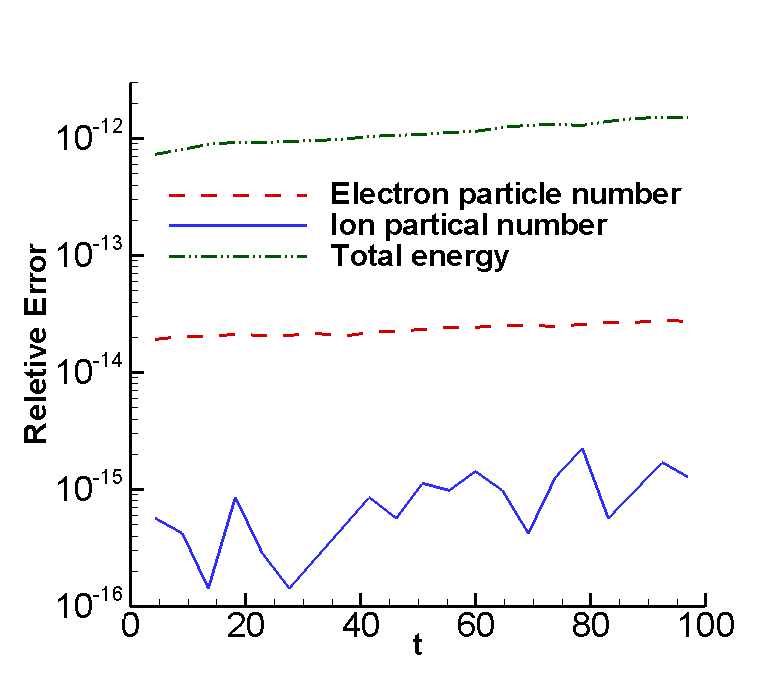}}
\subfigure[$\textnormal{\bf Scheme-2}$. $\mu_i=1/1836$.]{\includegraphics[width=0.45\textwidth]{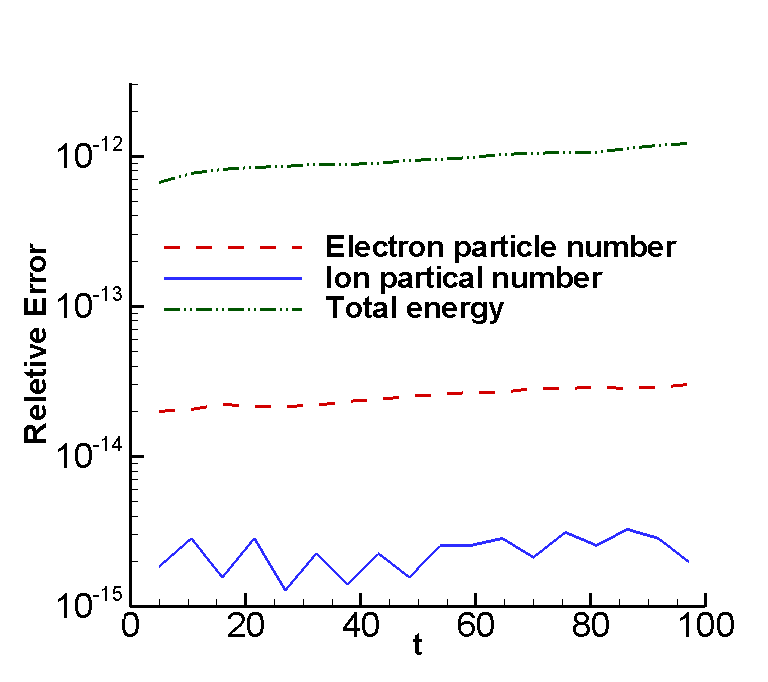}}
\caption{Evolution of absolute value of relative error in total particle number and total energy. $N_x=100, N_v=200.$}
\label{figure_landauconserve}
\end{figure}

\begin{figure}[!htbp]
\centering
%\subfigure[$\textnormal{\bf Scheme-1}$. $\mu_2=1/25$.]{\includegraphics[width=0.45\textwidth]{landauexplicit25.eps}}
%\subfigure[$\textnormal{\bf Scheme-1}$. $\mu_2=1/1836$.]{\includegraphics[width=0.45\textwidth]{landauexplicitreal.eps}}\\
\subfigure[$\textnormal{\bf Scheme-2}$. $\mu_i=1/25$.]{\includegraphics[width=0.4\textwidth]{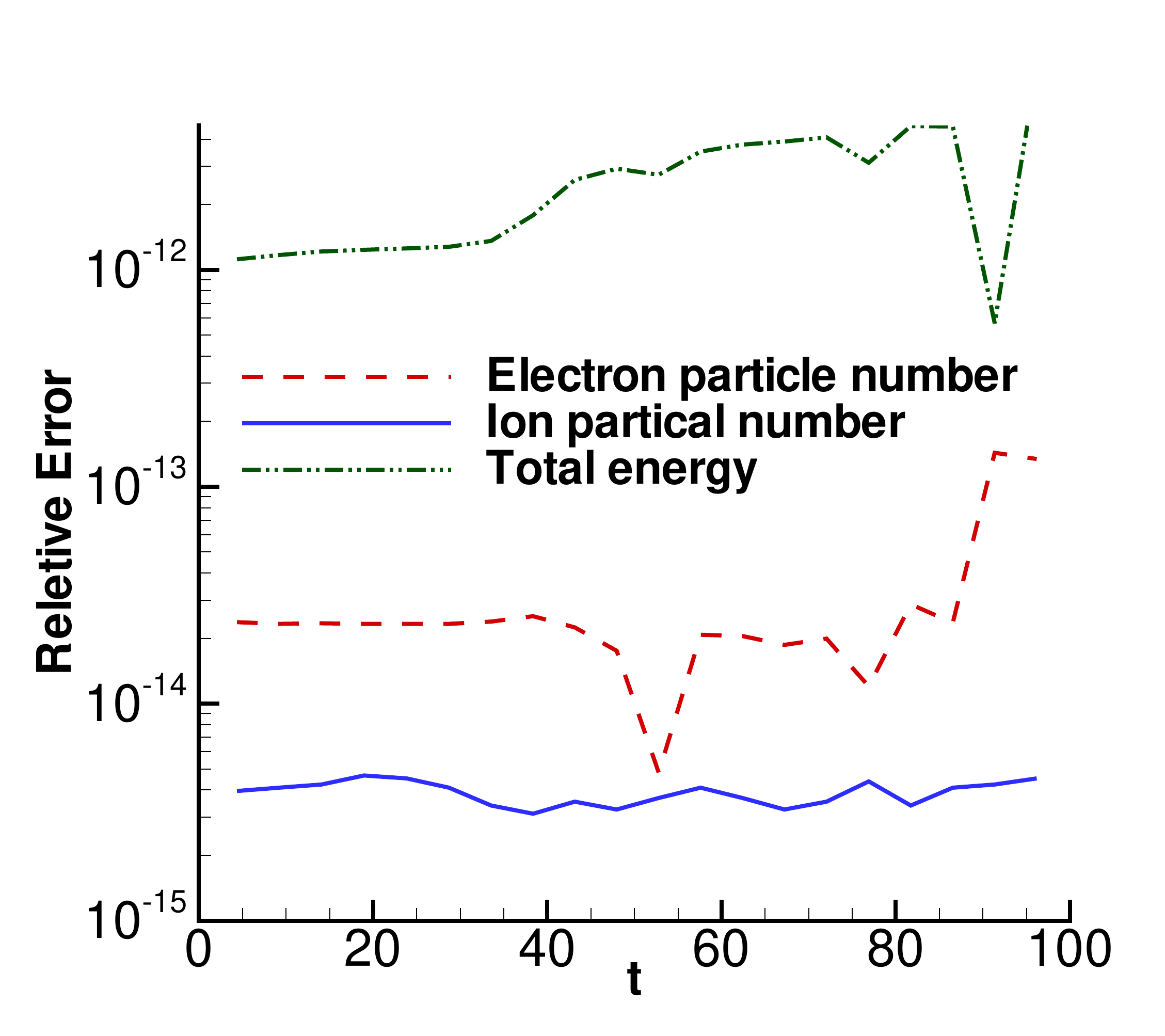}}
\subfigure[$\textnormal{\bf Scheme-2}$. $\mu_i=1/1836$.]{\includegraphics[width=0.4\textwidth]{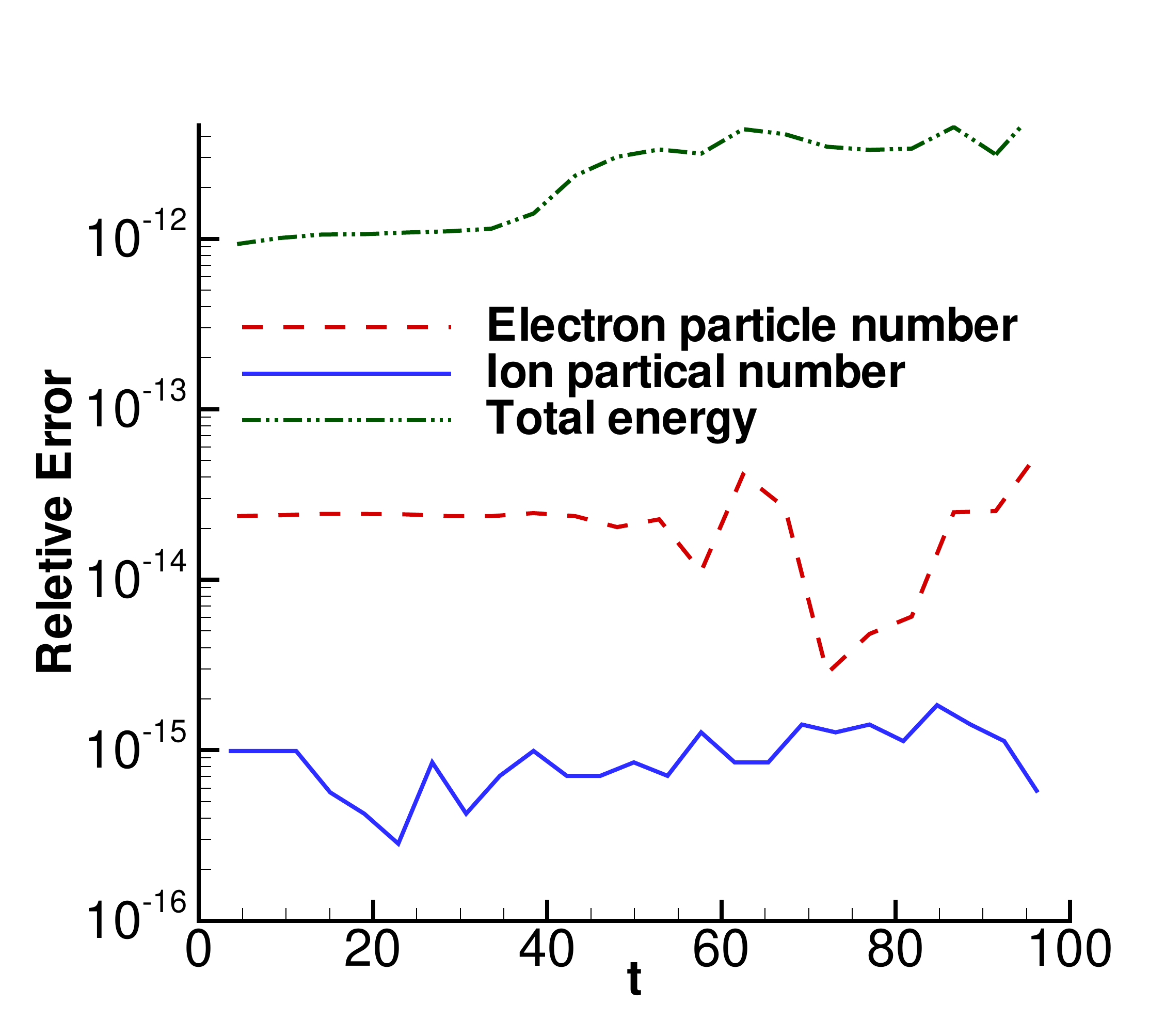}}
\caption{Evolution of absolute value of relative error in total particle number and total energy with a coarse mesh. $CFL = 5$. $N_x=40, N_v=80.$}
\label{figure_landauconservemesh40}
\end{figure}

We first verify the conservation of total particle number and total energy for our schemes. 
 Figure \ref{figure_landauconserve} shows the absolute value of relative error of the total particle number and total energy for   $\textnormal{\bf Scheme-1}$ ($CFL = 0.13,$ typical time step size $\Delta t \approx 0.002$) and 
$\textnormal{\bf Scheme-2}$ ($CFL = 5,$ typical time step size $ \Delta t \approx 0.077$) with two sets of mass ratios with $N_x=100$ and $N_v=200$. We can see that all errors stay small, below $10^{-11}$ for the whole duration of the simulation.
 In Figure \ref{figure_landauconservemesh40}, we use a coarse mesh ($N_x=40$, $N_v=80$, $CFL=5$) to plot the errors in the conserved quantities to demonstrate that the conservation properties of our schemes are mesh independent. We use  $\textnormal{\bf Scheme-2} $ to demonstrate the behavior. Upon comparison with the results from finer mesh  in Figures \ref{figure_landauconserve}, we conclude that the mesh size has no impact on the conservation of total particle number and total energy as predicted by  Theorems \ref{thm:mass} and \ref{thm:fullenergy}. This verify the conservation properties are independent of mesh sizes, and we are allowed to use even under-resolved mesh to achieve high accuracy in particle number and energy conservations. 
 
\begin{figure}[!htbp]
\centering
\subfigure[$\log FM_1$]{\includegraphics[width=0.4\textwidth]{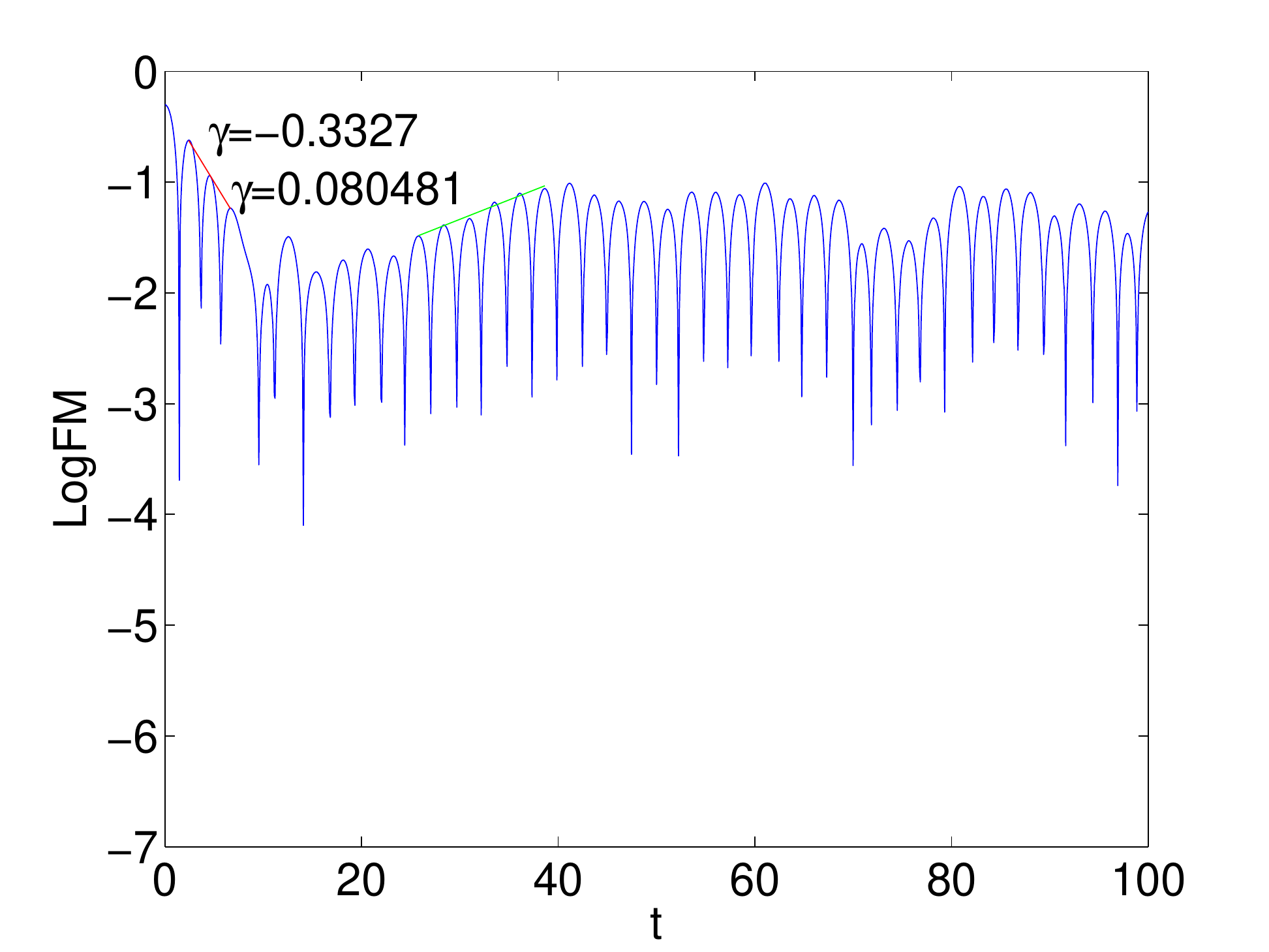}}
\subfigure[$\log FM_2$]{\includegraphics[width=0.4\textwidth]{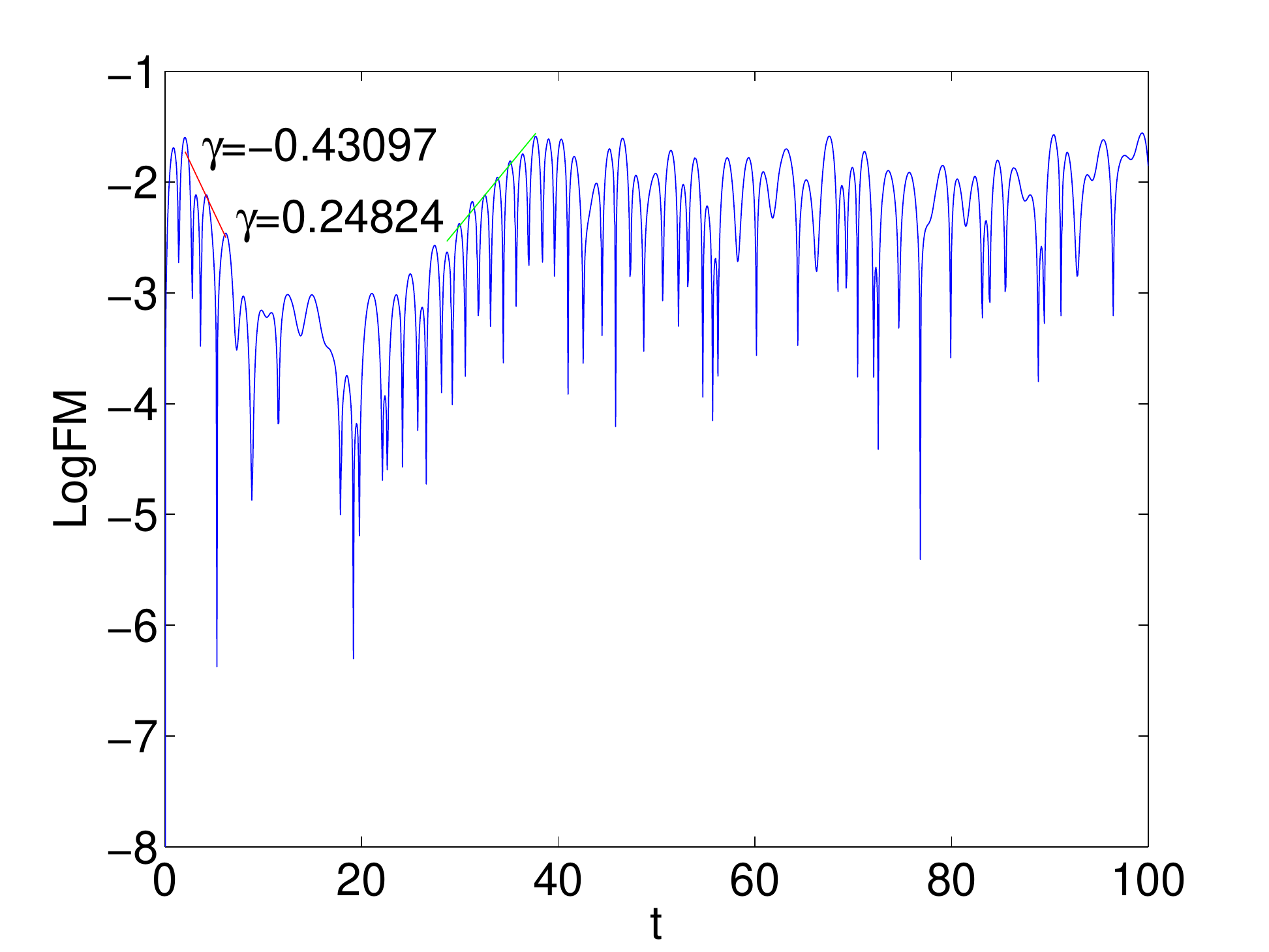}}\\
\subfigure[$\log FM_3$]{\includegraphics[width=0.4\textwidth]{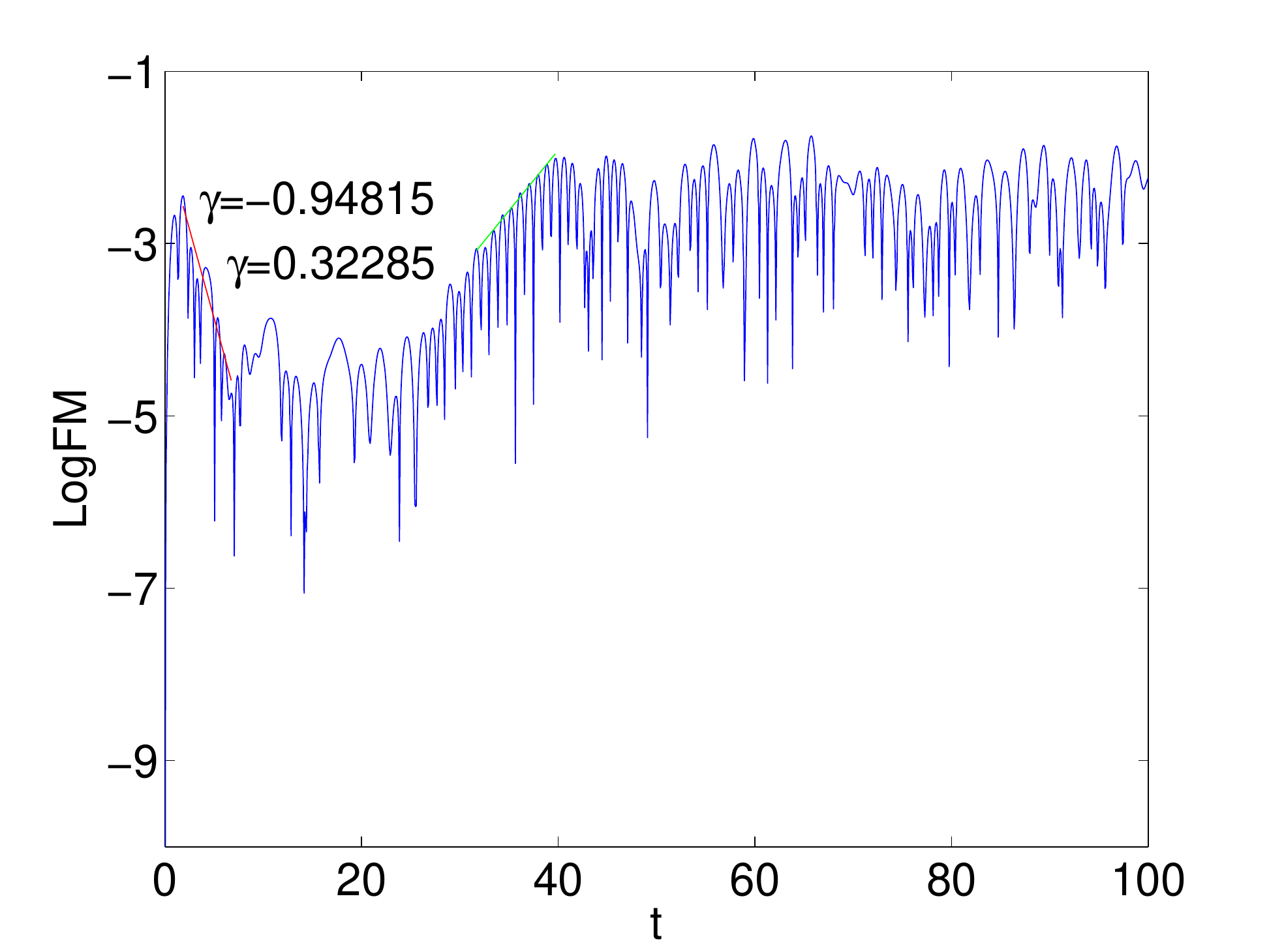}}
\subfigure[$\log FM_4$]{\includegraphics[width=0.4\textwidth]{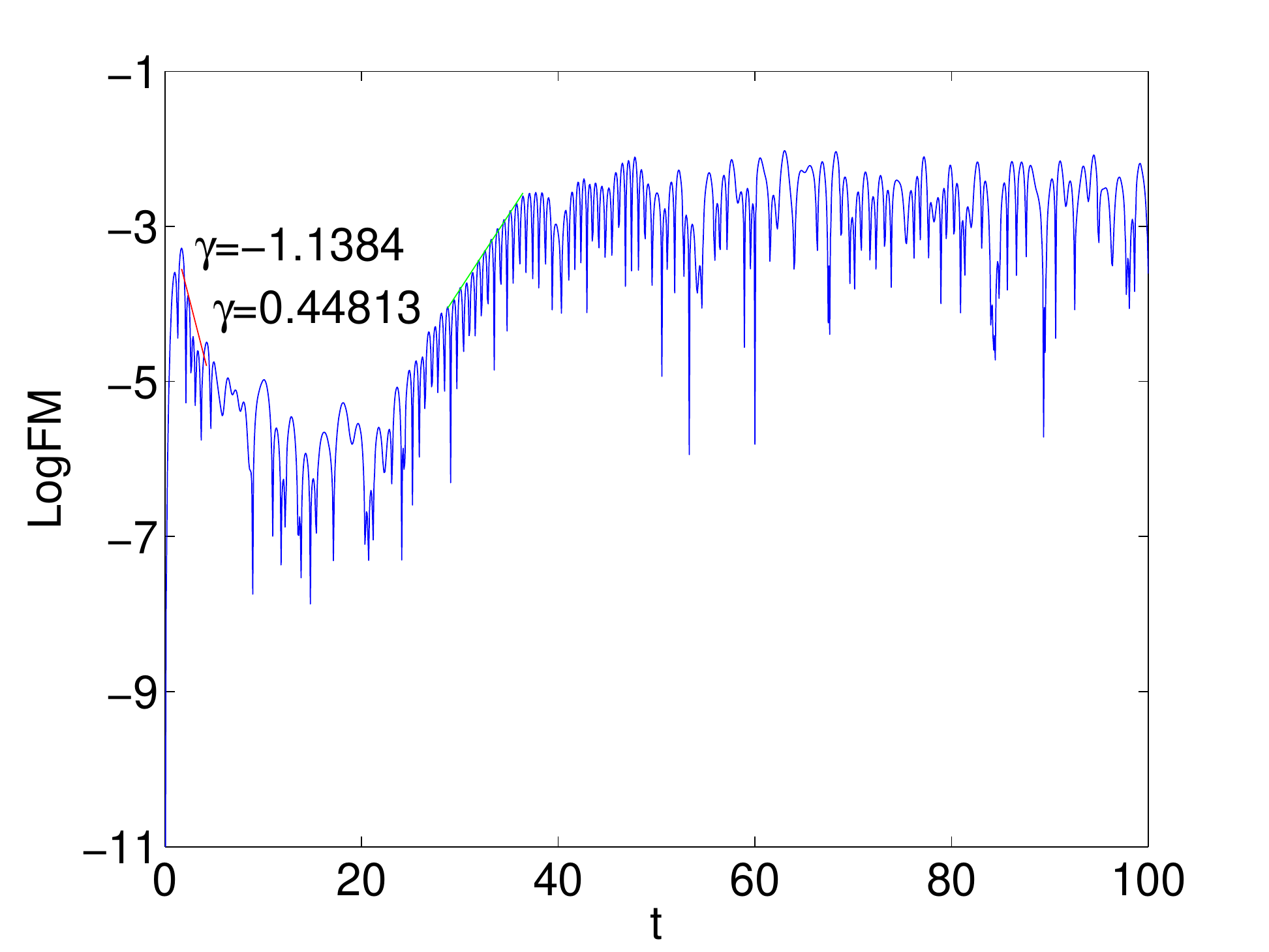}}
\caption{Log Fourier modes of Landau damping. $\textnormal{\bf Scheme-2}$. $\mu_2=1/25$. $CFL = 5$ (typical time step size $\Delta t \approx 0.077)$.  $N_x=100, N_v=200.$}
\label{figure_logfm25}
\end{figure}

\begin{figure}[!htbp]
\centering
\subfigure[$\log FM_1$]{\includegraphics[width=0.45\textwidth]{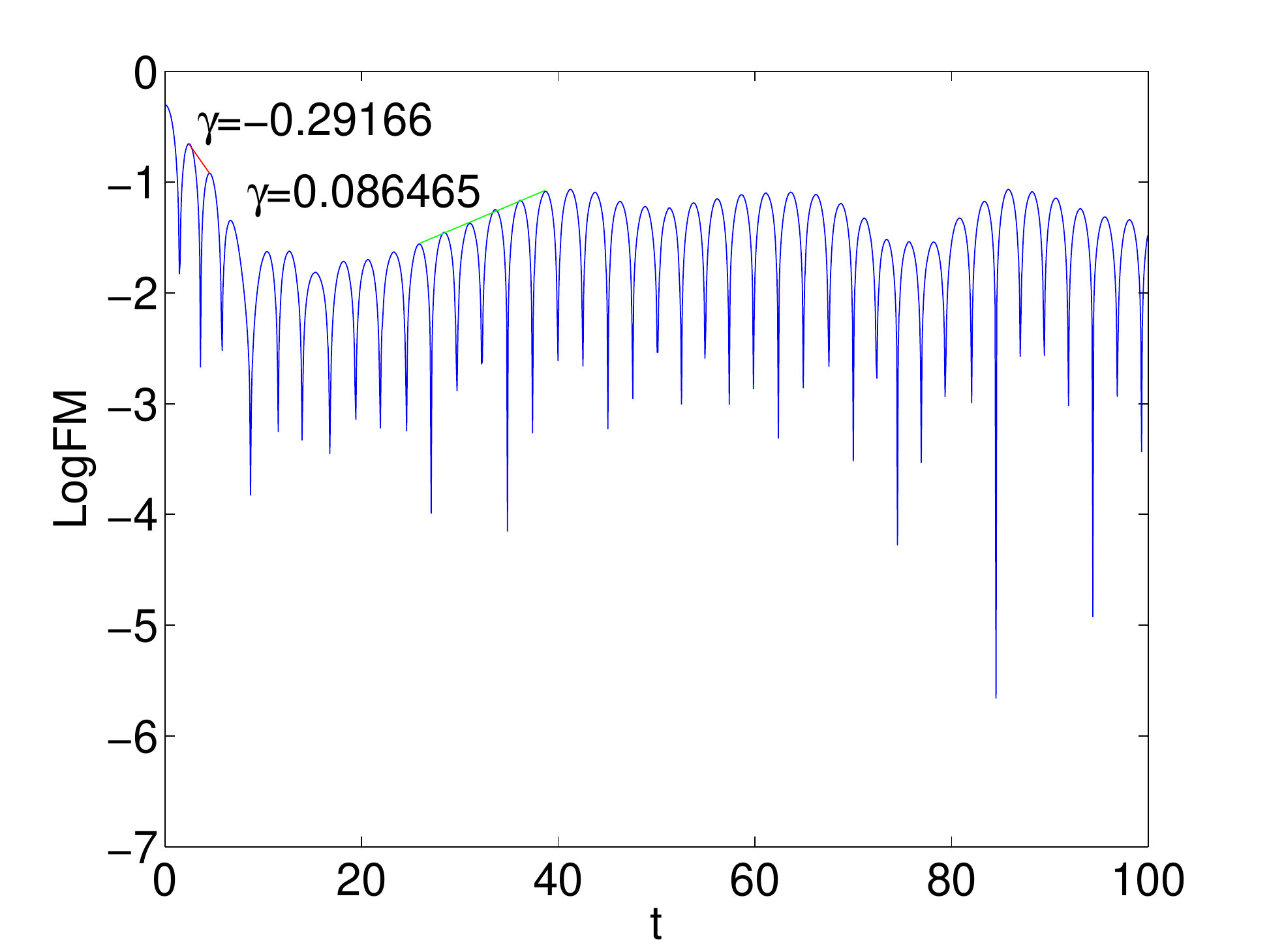}}
\subfigure[$\log FM_2$]{\includegraphics[width=0.45\textwidth]{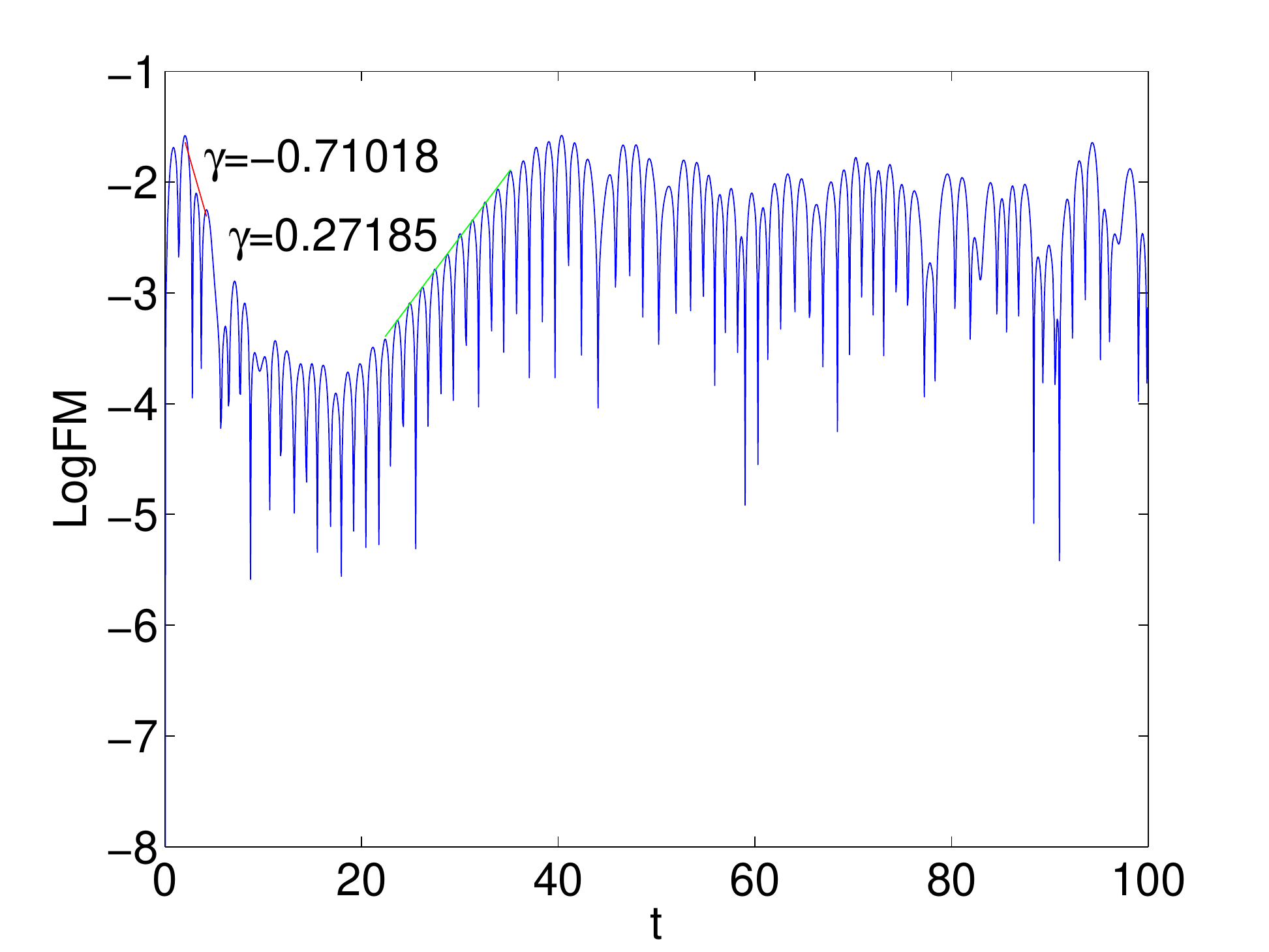}}\\
\subfigure[$\log FM_3$]{\includegraphics[width=0.45\textwidth]{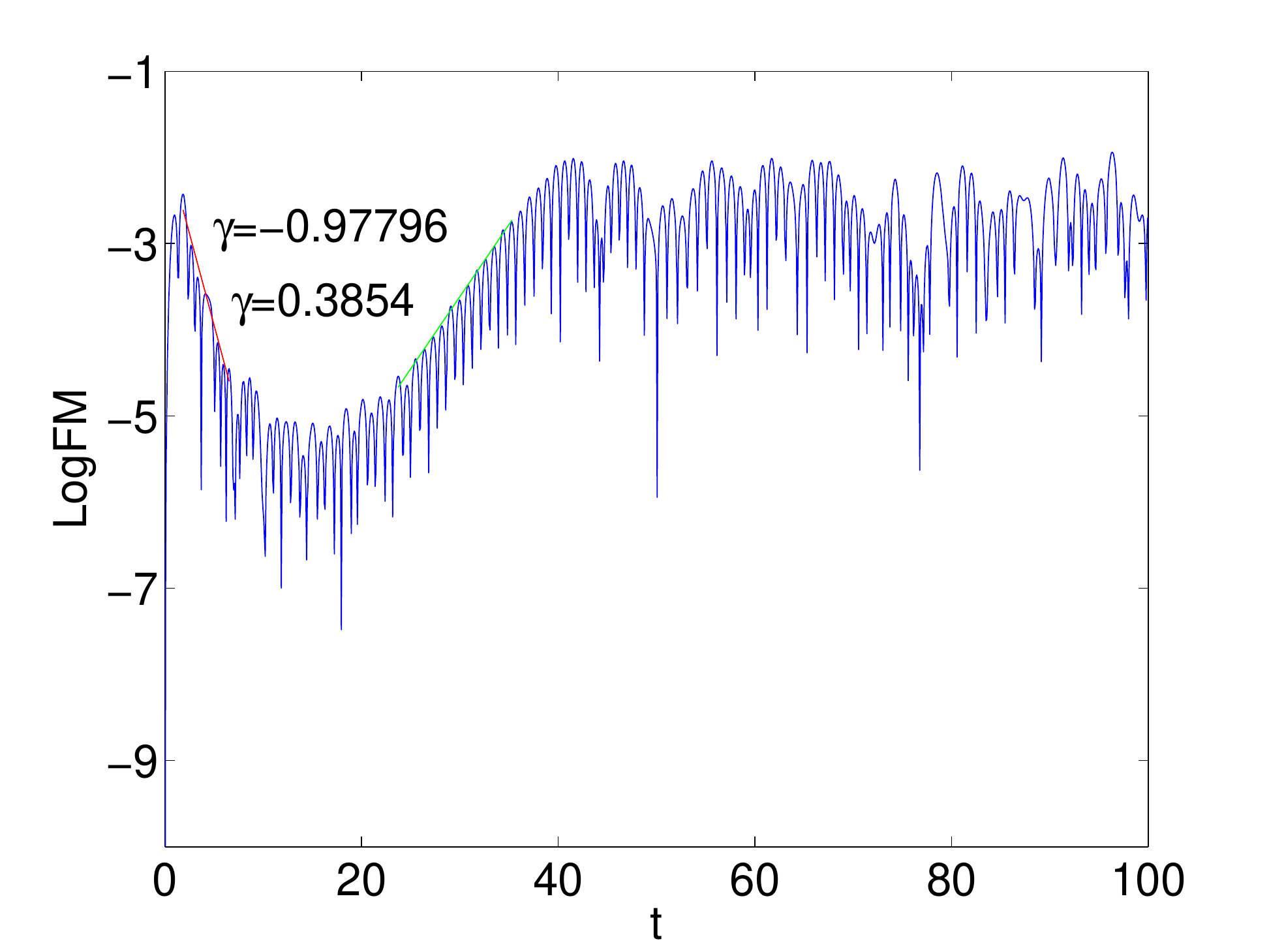}}
\subfigure[$\log FM_4$]{\includegraphics[width=0.45\textwidth]{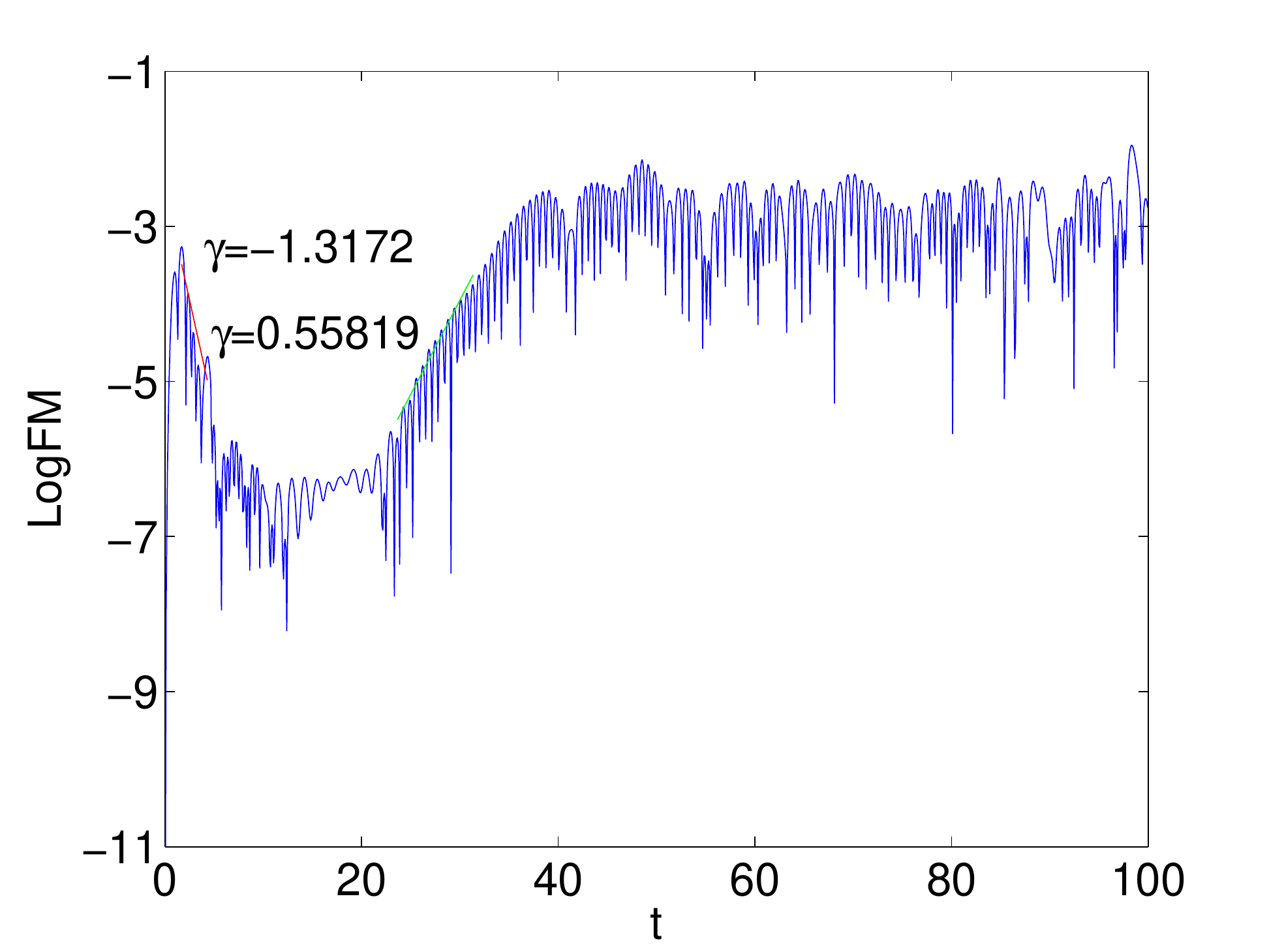}}
\caption{Log Fourier modes of Landau damping. $\textnormal{\bf Scheme-2}$. $\mu_2=1/1836$. $CFL = 5$ (typical time step size $\Delta t \approx 0.077)$.  $N_x=100, N_v=200.$}
\label{figure_logfmreal}
\end{figure}

Next, we compare the difference  of numerical simulations between different mass ratios. We plot the first four   Log Fourier modes for the electric field.
The $n$-th Log Fourier mode for the electric field $E(x,t)$  is defined as
$$
logF\!M_n(t)=\log_{10} \left(\frac{1}{L} \sqrt{\left|\int_0^L E(x, t) \sin(\kappa nx) \, dx \right|^2 +
\left|\int_0^L E(x, t) \cos(\kappa nx)\,  dx \right|^2} \right).$$
They are crucial quantities to investigate the qualitative behavior of the solution \cite{chengknorr_76}. 
By comparing Figures \ref{figure_logfm25} and \ref{figure_logfmreal}, we see that both mass ratios demonstrate similar qualitative behavior for the four modes. Upon a detailed comparison to 
 the one-species Landau damping result \cite{cheng_va},  the real mass ratio clearly yields decay and growth rates that are much closer to the one-species Landau damping, showing the convergence of the model to the one-species limit when $\mu_i \rightarrow 0$.
 
\begin{figure}[!htbp]
\centering
\subfigure[$\log FM_1$]{\includegraphics[width=0.45\textwidth]{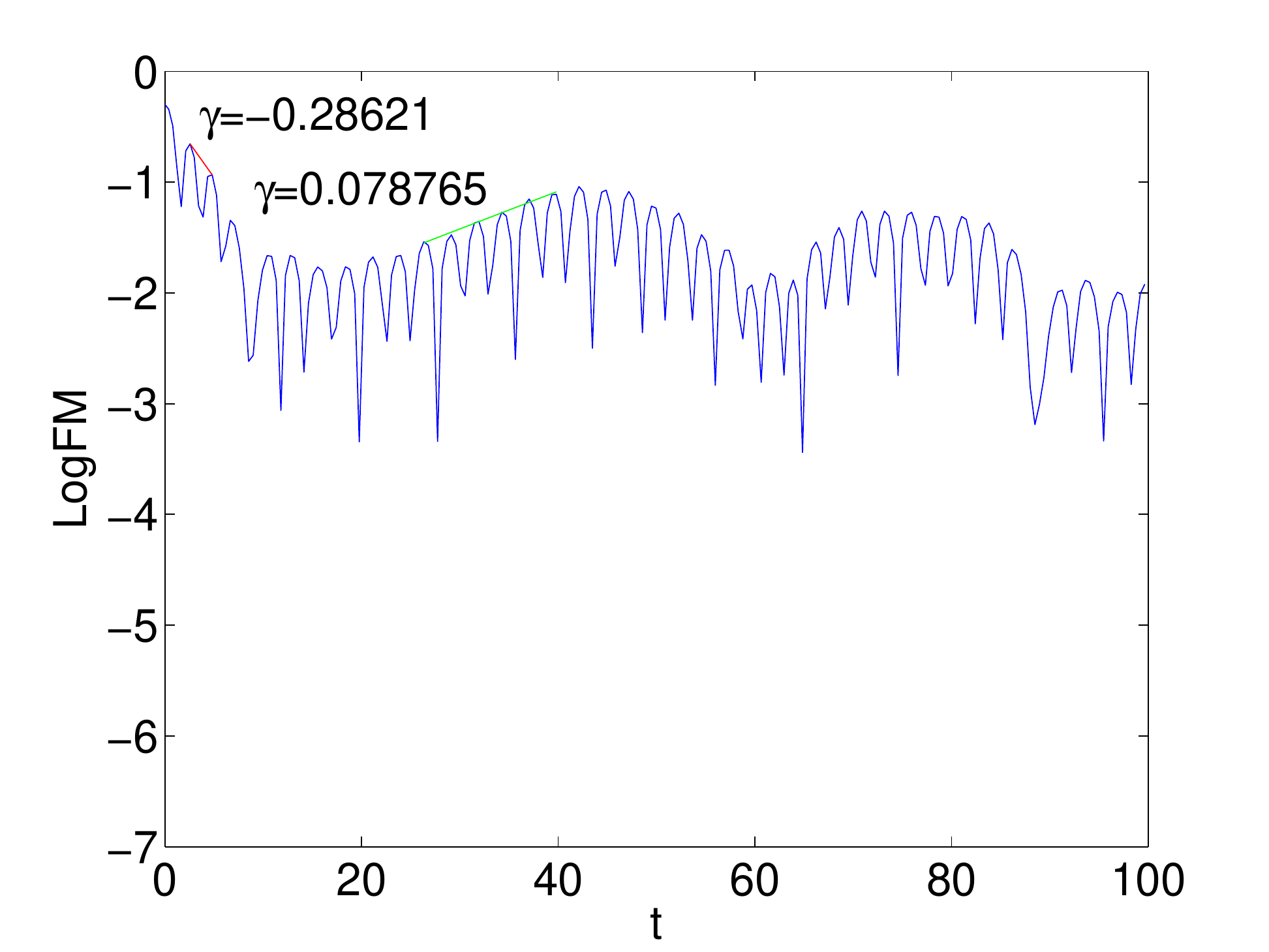}}
\subfigure[$\log FM_2$]{\includegraphics[width=0.45\textwidth]{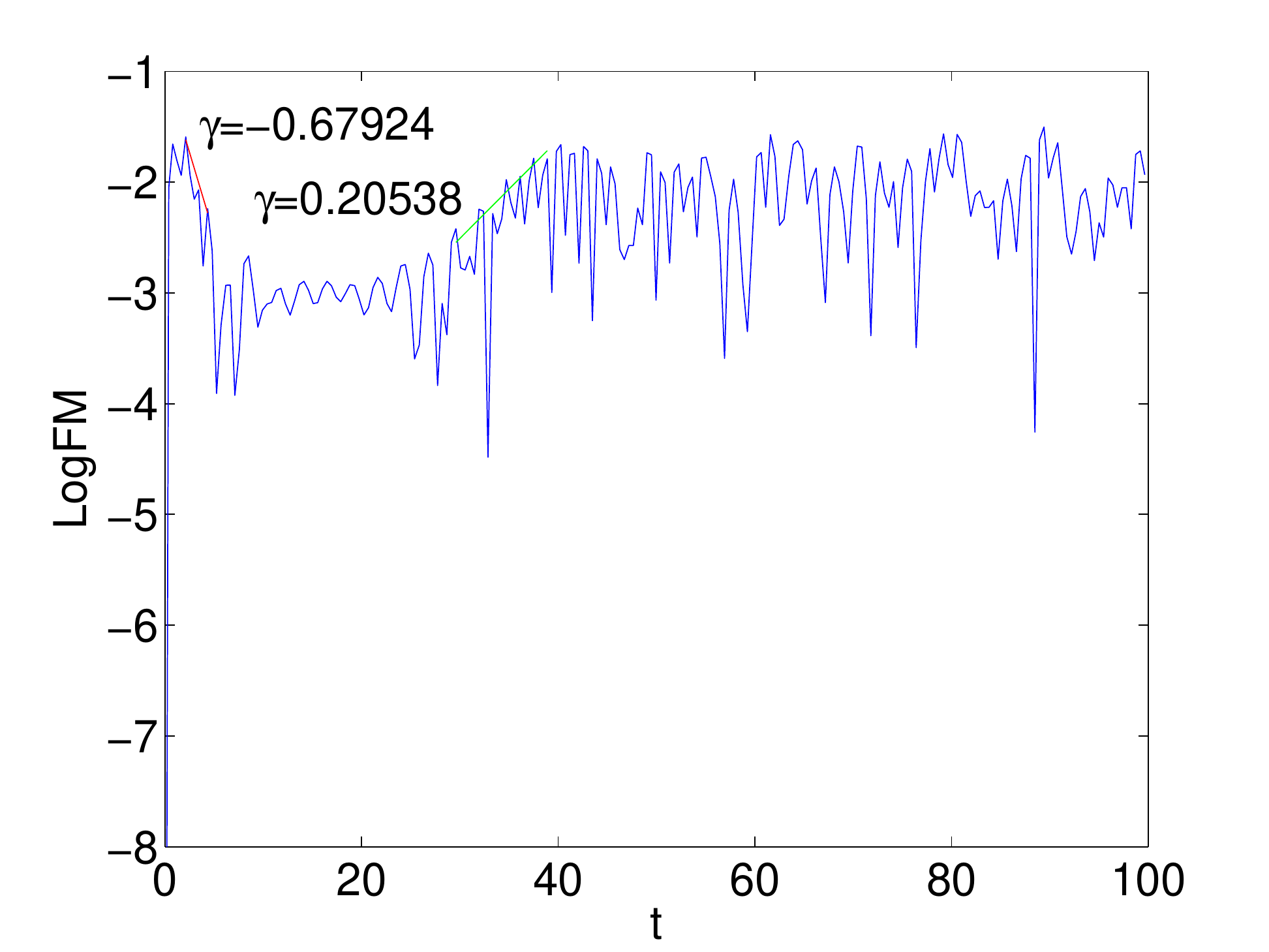}}\\
\subfigure[$\log FM_3$]{\includegraphics[width=0.45\textwidth]{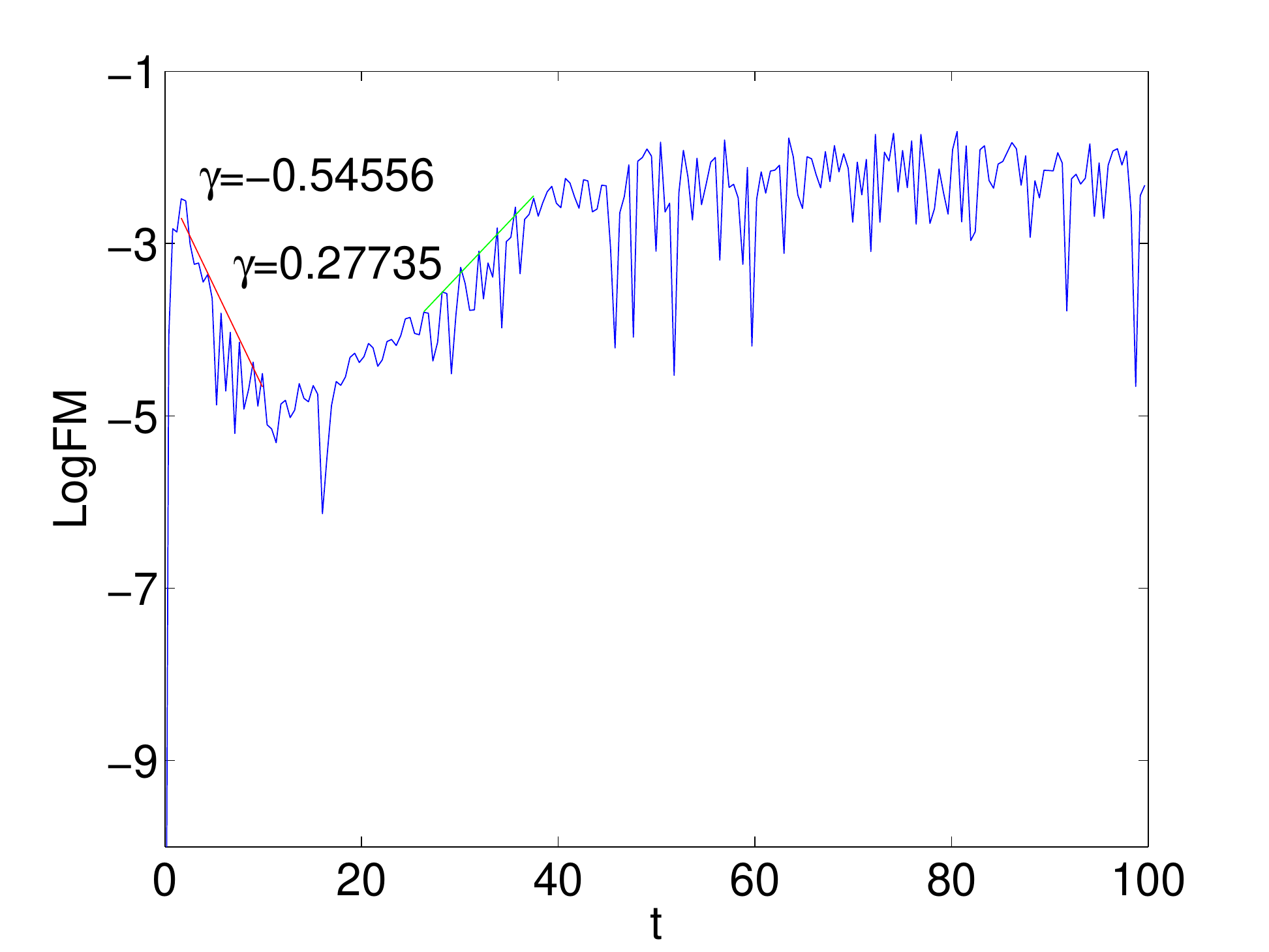}}
\subfigure[$\log FM_4$]{\includegraphics[width=0.45\textwidth]{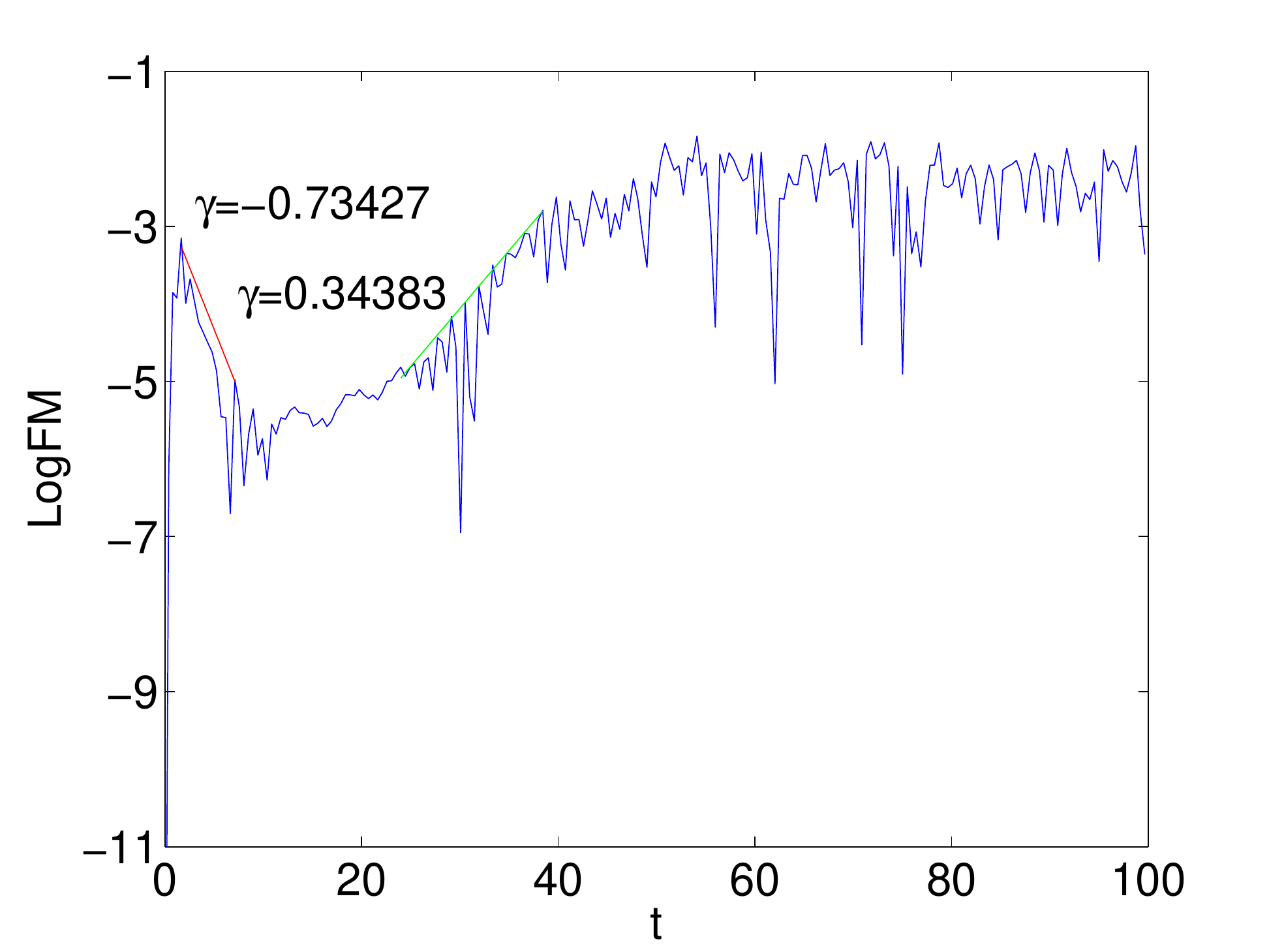}}
\caption{Log Fourier modes of Landau damping. $\textnormal{\bf Scheme-2}$. $\mu_2=1/1836$. $CFL= 30$ (typical time step size $\Delta t \approx 0.46)$, $N_x=100, N_v=200.$}
\label{figure_logfmrealcfl30}
\end{figure}

\begin{figure}[!htbp]
\centering
\subfigure[$\log FM_1$]{\includegraphics[width=0.45\textwidth]{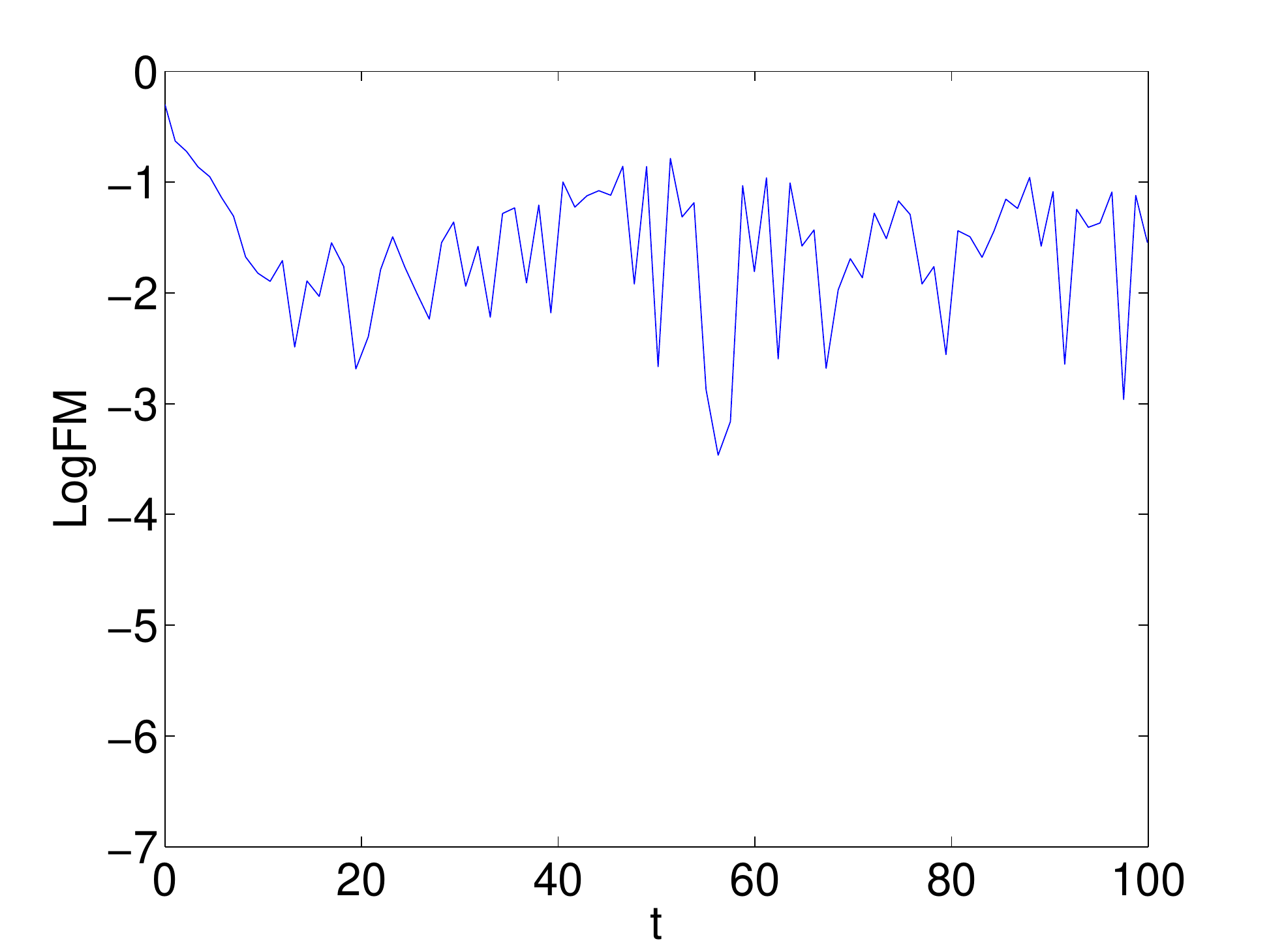}}
\subfigure[$\log FM_2$]{\includegraphics[width=0.45\textwidth]{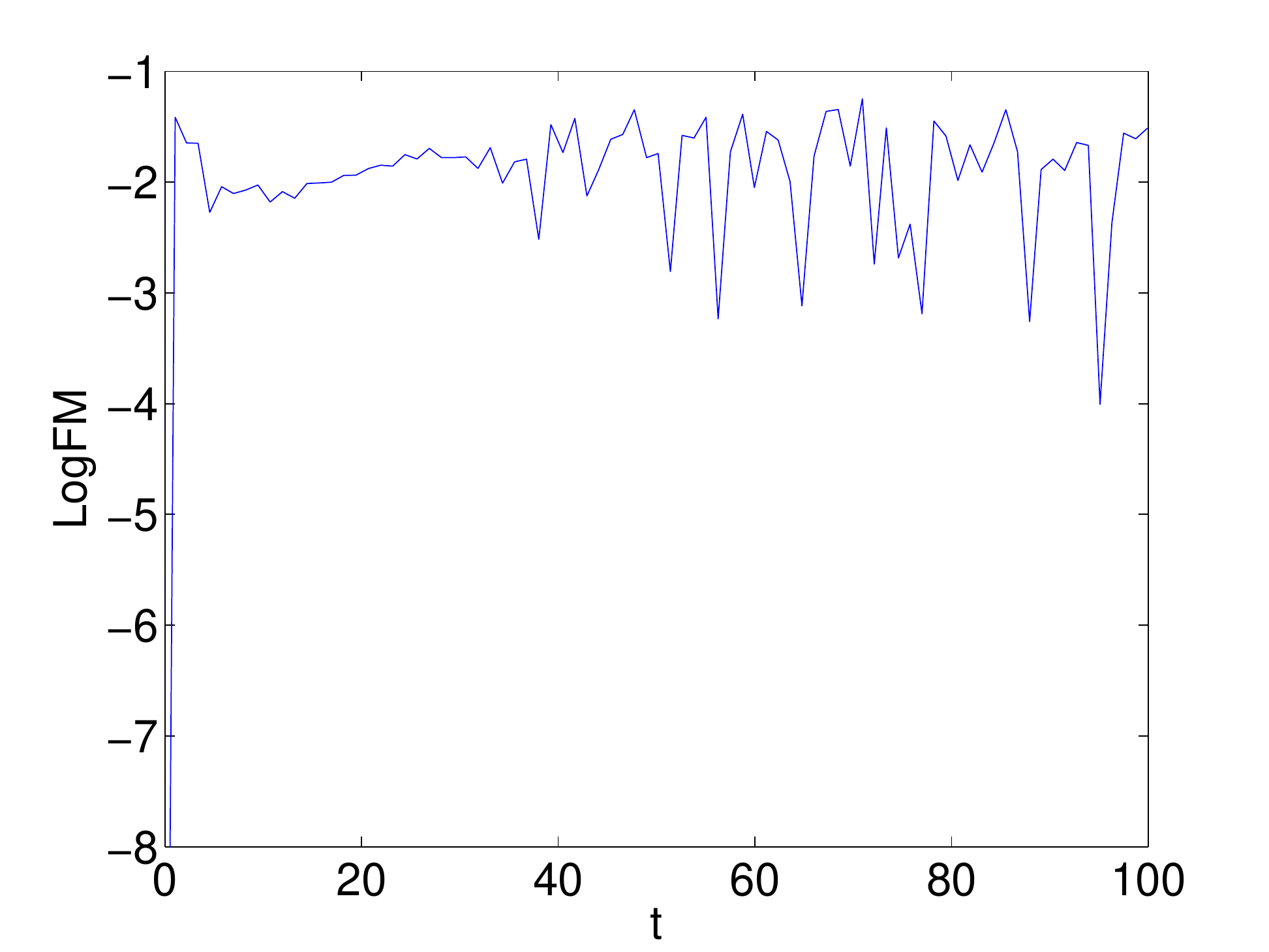}}\\
\subfigure[$\log FM_3$]{\includegraphics[width=0.45\textwidth]{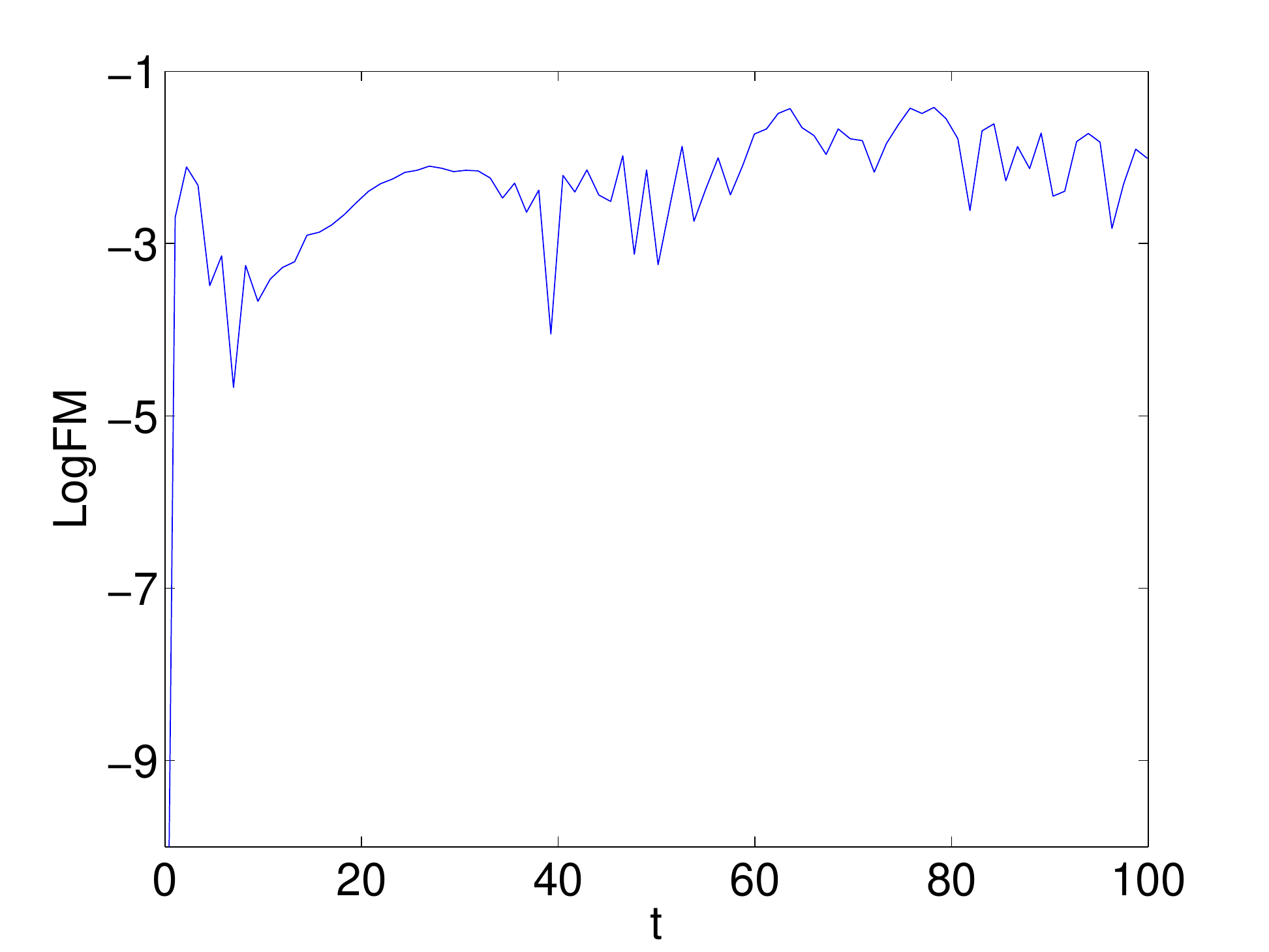}}
\subfigure[$\log FM_4$]{\includegraphics[width=0.45\textwidth]{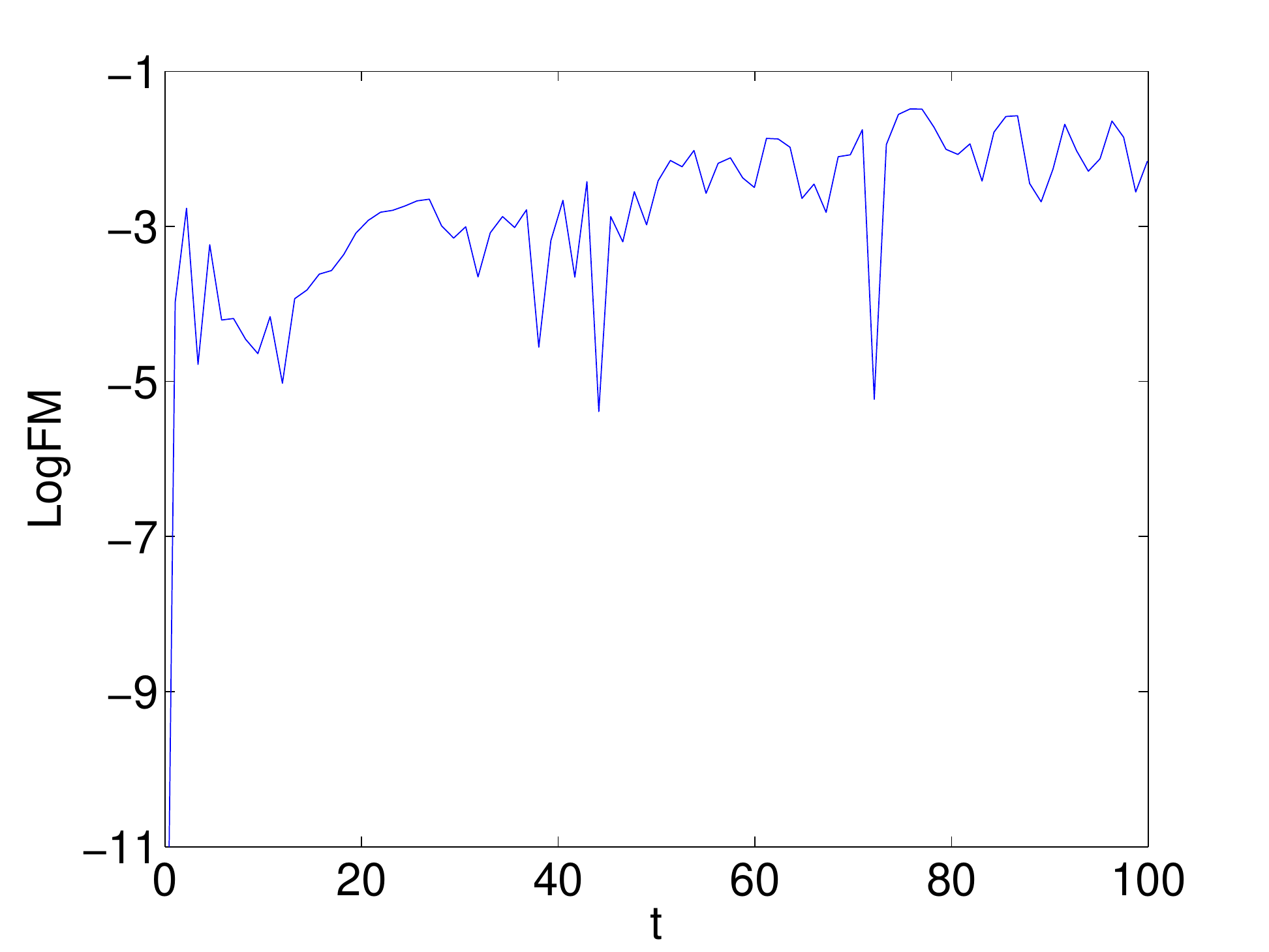}}
\caption{Log Fourier modes of Landau damping. $\textnormal{\bf Scheme-2}$. $\mu_2=1/1836$. $CFL= 80$ (typical time step size $\Delta t \approx 1.2)$, $N_x=100, N_v=200.$}
\label{figure_logfmrealcfl80}
\end{figure}

\begin{figure}[!htbp]
\centering
\subfigure[$\log FM_1$]{\includegraphics[width=0.45\textwidth]{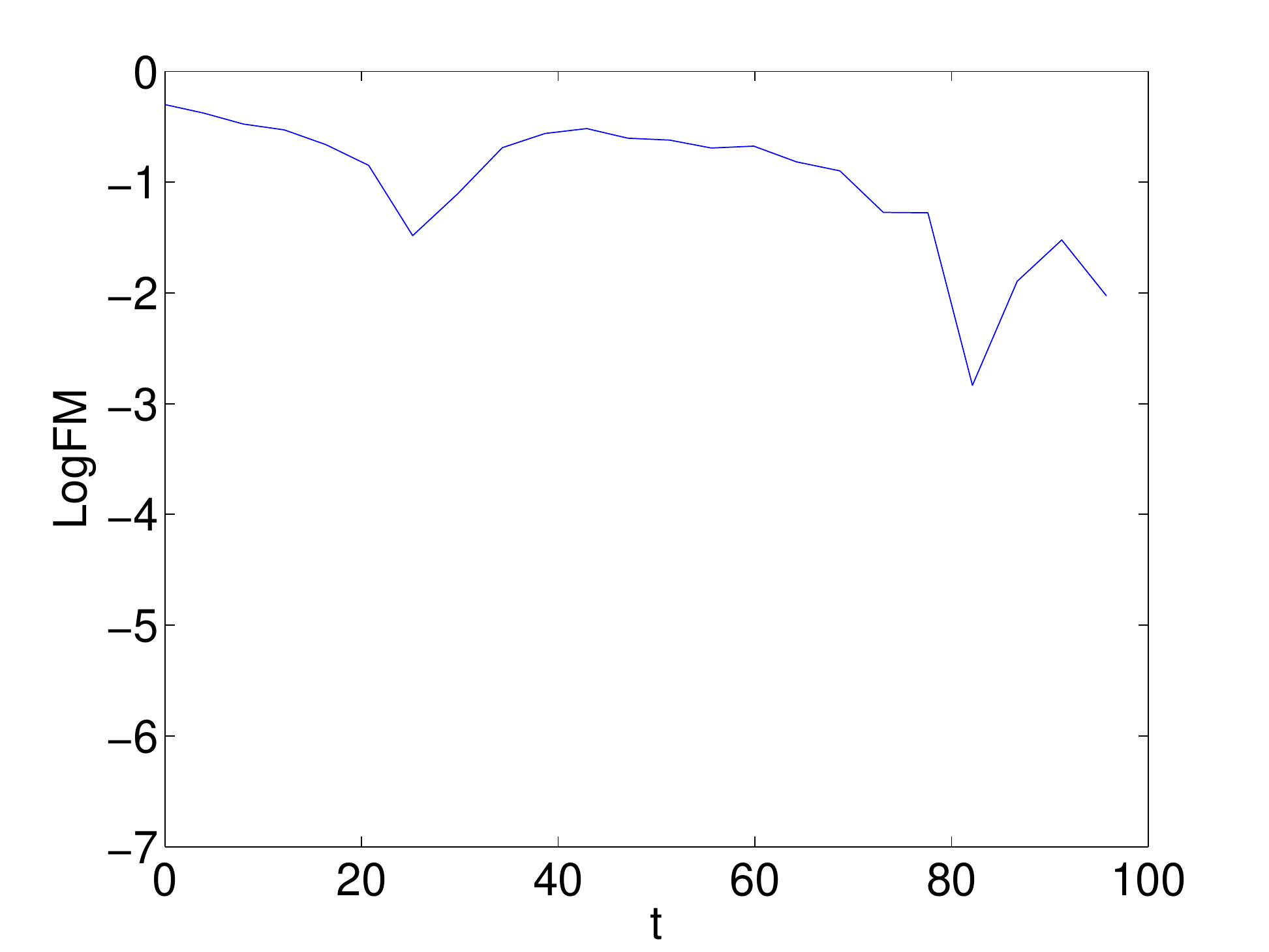}}
\subfigure[$\log FM_2$]{\includegraphics[width=0.45\textwidth]{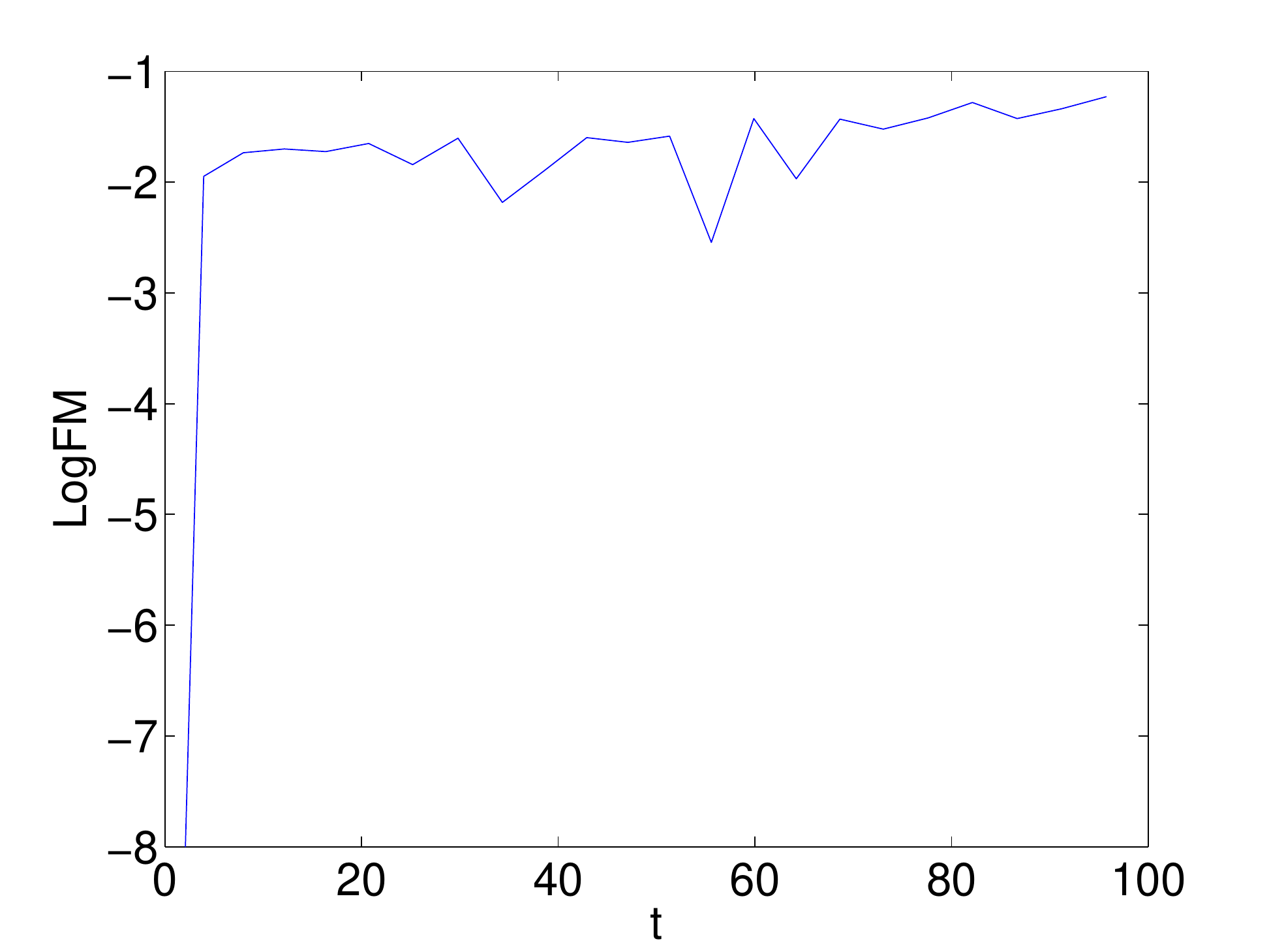}}\\
\subfigure[$\log FM_3$]{\includegraphics[width=0.45\textwidth]{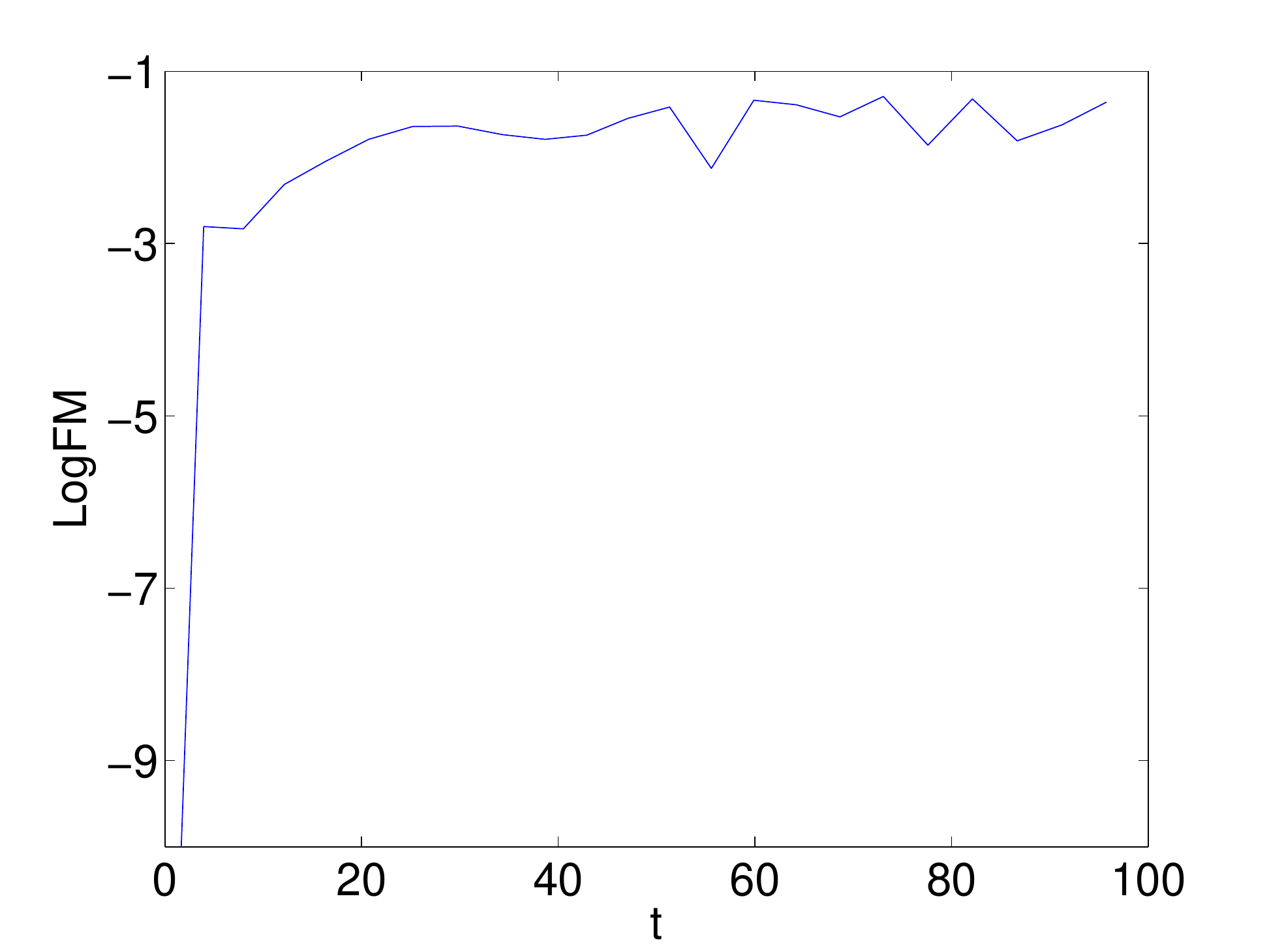}}
\subfigure[$\log FM_4$]{\includegraphics[width=0.45\textwidth]{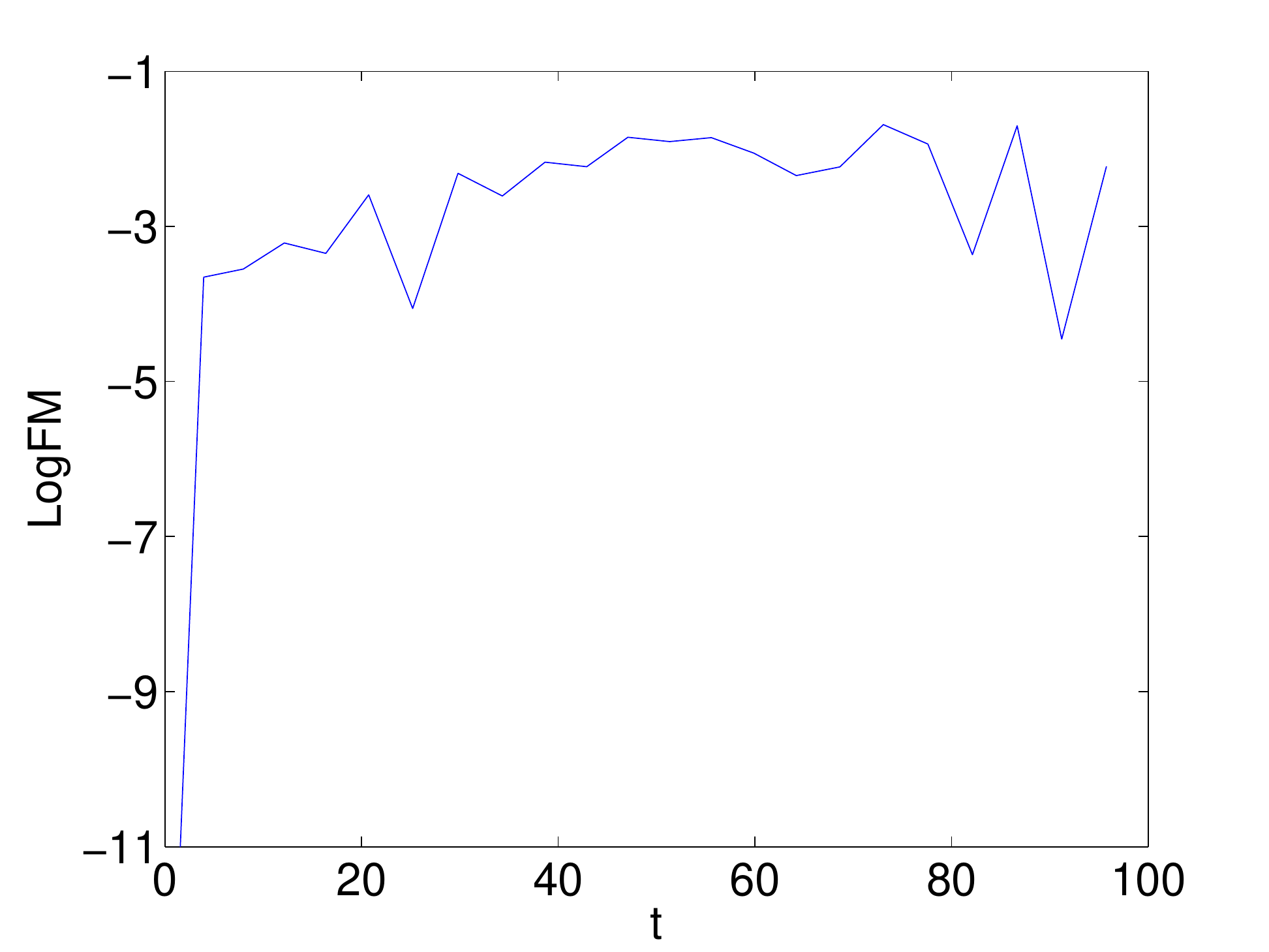}}
\caption{Log Fourier modes of Landau damping. $\textnormal{\bf Scheme-2}$. $\mu_2=1/1836$. $CFL= 300$ (typical time step size $\Delta t \approx 4.3)$, $N_x=100, N_v=200.$}
\label{figure_logfmrealcfl300}
\end{figure}

  Finally, we investigate the influence of the time step size $\Delta t$ on the behavior of the solutions. In  Figures \ref{figure_logfmreal}, \ref{figure_logfmrealcfl30}, \ref{figure_logfmrealcfl80} and \ref{figure_logfmrealcfl300}, we   plot the results using the real mass ratio with  varying CFL numbers by the implicit scheme. Although the four simulations are all numerically stable, we clearly observe that the numerical runs with $\Delta t>1$ (i.e. $w_{pe} \Delta t>1$ in the unscaled variables) in Figures \ref{figure_logfmrealcfl80} and \ref{figure_logfmrealcfl300} fail to capture the subtle electron kinetic effects. Naturally, larger time steps will filter out high frequency in time, and if the electron kinetic effects are important, it would be necessary to use time step size smaller than $w_{pe}^{-1}$. 
  %On the other hand, $w_{pi}=w_{pe} \sqrt{\mu_i}$


%\begin{figure}[!htbp]
%\subfigure[$\log FM_1$]{\includegraphics[width=0.45\textwidth]{landauimplicitzerologfm1.eps}}
%\subfigure[$\log FM_2$]{\includegraphics[width=0.45\textwidth]{landauimplicitzerologfm2.eps}}\\
%\subfigure[$\log FM_3$]{\includegraphics[width=0.45\textwidth]{landauimplicitzerologfm3.eps}}
%\subfigure[$\log FM_4$]{\includegraphics[width=0.45\textwidth]{landauimplicitzerologfm4.eps}}
%\caption{Log Fourier modes of Landau damping. $\textnormal{\bf Scheme-3}$. $\mu_2=0$. $CFL = 5$.  $100\times 200$ mesh.}
%\end{figure}

%\clearpage
\subsection{Current-driven ion-acoustic instability}
\label{sec:cdiaw}
In this subsection, we perform a detailed numerical study of the current-driven ion-acoustic instability. 
 Ion-acoustic waves are natural wave modes in unmagnetized plasmas. 
Current-driven ion-acoustic instability are generated by giving the electrons an initial uniform drift velocity $v_{de}$ relative to the ions such that $v_{de}>v_{crit}$ where 
$v_{crit}$ is the threshold for the ion-acoustic instability. This test example corresponds to the case of the two-species VA system with $\bJ_{ext}=\bJ_0$, $\bE_0(0)=0$.

The initial conditions of the distribution functions  are given by
\begin{subequations}
\begin{align}
f_{e}(x,v,0)&=\left(1+\sum_{n=1}^{N_{\max}}E_{tf}\kappa_n\cos(\kappa_n x+\varphi_n)\right)\frac{1}{\sqrt{2\pi}}e^{-(v-v_{de})^2/2} \label{feinit}\\
f_{i}(x,v,0)& = \frac{1}{\sqrt{2\pi\gamma}}e^{-v^2/2\gamma}
\end{align}
\end{subequations}
where $N_{\max}$ is the number of modes permitted in the simulation, $\varphi_n$ is a random phase, $E_{tf}$ is the thermal fluctuation level,
$v_{de}$ is the uniform drift velocity for the electrons and $\gamma=(T_im_e/T_em_i)^{1/2}$. The initial condition for the electric field 
$$
E(x,0)=-\sum_{n=1}^{N_{\max}}E_{tf}\sin(\kappa_n x+\varphi_n)
$$
is obtained by  the Poisson equation and clearly satisfies $E_0(0)=0$.

In the numerical runs, we use  simulation parameters $S_1$  as listed in Table \ref{table_parameter}, which are the rescaled version of the parameters used in Table 1  of \cite{ionacoustic}. In particular, we focus on the reduced  mass ratio with $\mu_i=m_e/m_i=1/25$ instead of the real mass ratio. 
As demonstrated in   \cite{ionacoustic}, the reduced mass ratio yields qualitatively similar results as the real mass ratio, but enables faster computations. 
%It is advantageous to use a reduced mass ratio when modeling current-driven ion-acousticwaves.
This is because the real mass ratio would require a large number of electron velocity grid points  in order to accommodate the relatively small range of resonant phase velocities. By using the reduced mass ratio, the run time of the simulation is kept to a reasonable level for us to perform 100 simulations with random phase perturbations.
%By increasing the number of electron velocity grid points, the time taken to perform all the calculations required to evolve the Vlasov equation by one time step is also increased, and the total run time of the simulation becomes too long to be practical. By using the reduced mass ratio, the run timeof the simulation is kept to a reasonable level.

In Table \ref{table_parameter},
$\lambda_{\min},\,\lambda_{\max}$ are the smallest and largest wavelengths of linearly unstable
ion-acoustic wave modes calculated by solving numerically the linear dispersion relation. The domain for $x$ is set to be  $[0,\, L]$, where $L=\lambda_{\max}$. 
Let $\kappa_0=2\pi/L$ denotes the wave number. $\kappa_n=n\kappa_0$ for $1\leq n\leq N_{\max}$ with
$N_{\max}=\lambda_{\max}/\lambda_{\min}=53$.
The domain of velocity $v$ for electron and ion is chosen to be $[-V_{c,e},V_{c,e}]$ and $[-V_{c,i},V_{c,i}]$, respectively, such that $f\simeq 0$ on the boundaries.
$v_{ph,\min},\,v_{ph,\max}$ are the smallest and largest phase velocity of the unstable
wave modes calculated from the solution to the linear dispersion. The resolution of the velocity grids is controlled by the phase velocities of the growing wave modes, such that 
there are three velocity points in the linear unstable region ($v_{ph,\min}<v<v_{ph,\max}$) so that interaction with the unstable wave modes is possible \cite{ionacoustic}. 
Under this scaling, $E_{tf}=6.76\times 10^{-5}$ in \eqref{feinit}. Under our scaling, all the variables ploted in the figures of the following sections can be read as the values of quantities of the reference \cite{ionacoustic} after multiplied by the corresponding factor listed in the last column of Table \ref{table_parameter}.

\begin{table}[!htbp]
\centering
\caption{Summary of simulation parameters (rescaled)}
\label{table_parameter}
%\begin{tabular}{ c | >{\centering}m{3cm}| m{1cm}}
%$\begin{tabu} spread 0.2in {X[$c] X[$c] | X[$c] X[$c] }
$\begin{tabu} [m] {c| c |c |c}

\hline
%1&2&3\\(
\text{Parameters} & S_1\text{(reduced mass ratio)} &\text{Variables plotted} & \text{Scaled factor}\\
\hline
m_i/m_e &25 & t & \omega_{pe} \\
T_e/T_i & 2 & x & 3.97 \\ 
\lambda_{\min} & 7.98 &\theta_e^m&\sqrt{2}\\
\lambda_{\max} & 426.60 &f_e  & 11.81\\
v_{ph,\min}   &0.23 & f_i& 11.81\\
v_{ph,\max} & 0.29 &\eta& 7.58\times 10^5 \\
V_{c,e} & 10.30 &E & 0.504\\
V_{c,i} & 2.87&  \kappa &0.252 \\
N_x & 500 & & \\
N_{v,e},\,N_{v,e}& 890& & \\
%v_{de} &1.7&v & \theta_e^m&\sqrt{2}\\
%\lambda_{\min} & 7.98 &f_e  & 11.81\\
%\lambda_{\max} & 426.60 &f_i& 11.81\\
%v_{ph,\min}   &0.23 & \eta& 7.58\times 10^5\\
%v_{ph,\max} & 0.29 &E & 0.504 \\
%V_{c,e} & 10.30 &\kappa &0.252 \\
%V_{c,i} & 2.87& & \\
%N_x & 500 & & \\
%N_{v,e},\,N_{v,e}& 890& & \\
\hline
\end{tabu}$
\end{table}

\begin{figure}[!htbp]
 \centering
\subfigure[$\textnormal{\bf Scheme-1}$. $\mu_i=1/25$.]{\includegraphics[width=0.45\textwidth]{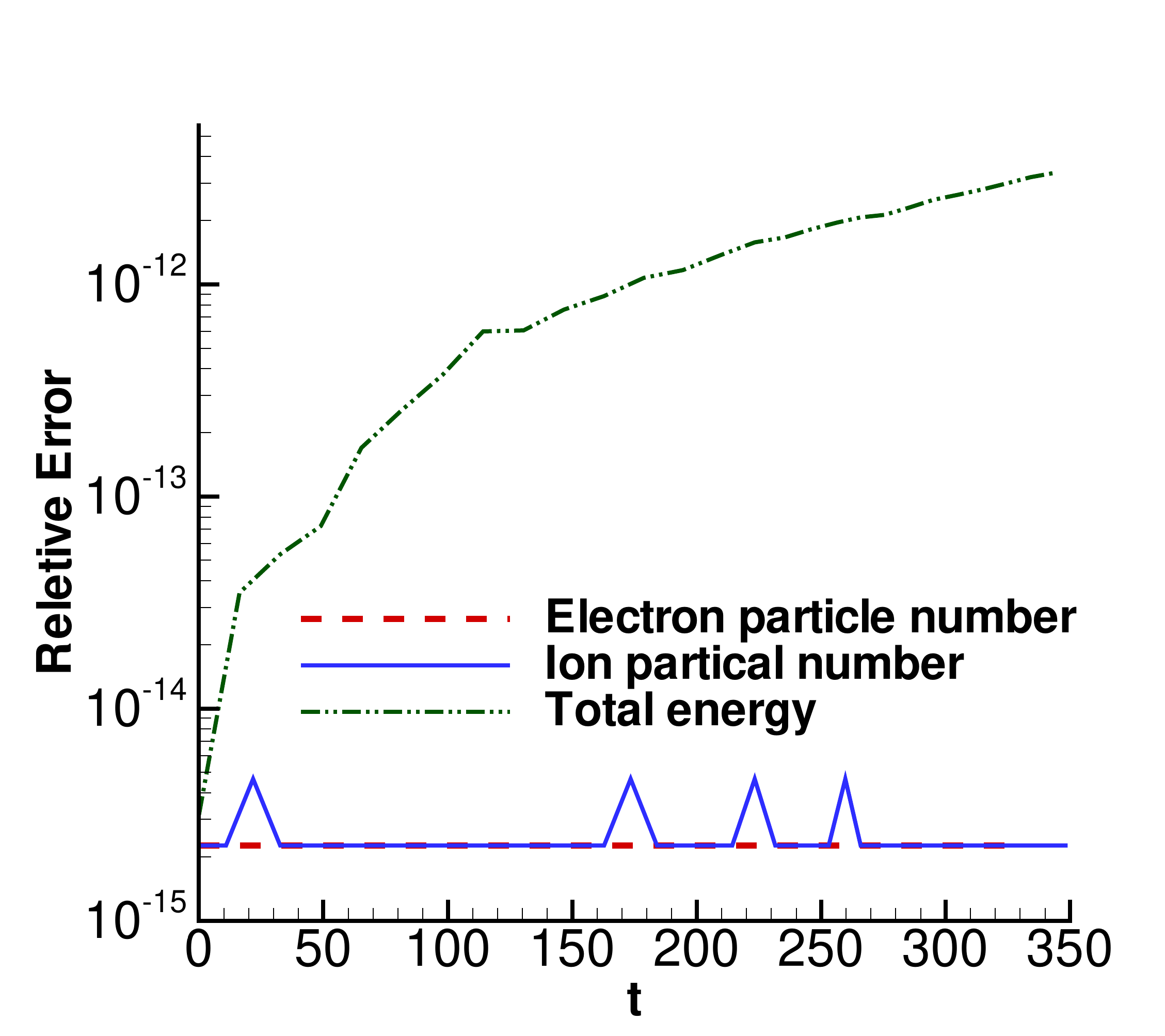}}
\subfigure[$\textnormal{\bf Scheme-2}$. $\mu_i=1/25$.]{\includegraphics[width=0.45\textwidth]{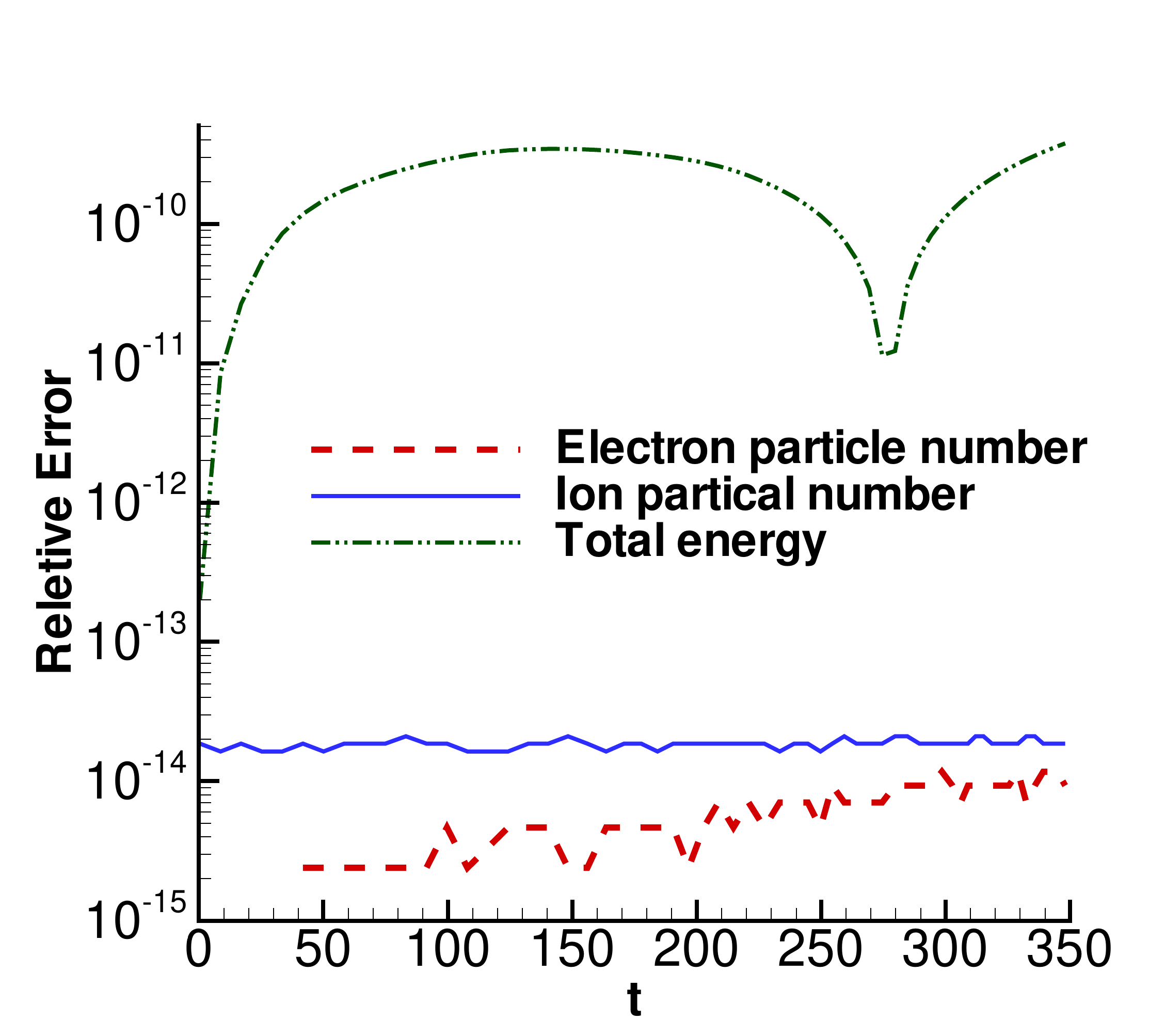}}
%\subfigure[$\textnormal{\bf Scheme-1}$. $\mu_i=1/1836$.]{\includegraphics[width=0.45\textwidth]{landauexplicitreal.eps}}\\
%\subfigure[$\textnormal{\bf Scheme-2}$. $\mu_i=1/25$.]{\includegraphics[width=0.45\textwidth]{landauimplicit25.eps}}
%\subfigure[$\textnormal{\bf Scheme-2}$. $\mu_i=1/1836$.]{\includegraphics[width=0.45\textwidth]{landauimplicitreal.eps}}
\caption{Evolution of absolute value of relative error in total particle number and total energy. }
\label{figure_iaconserve}
\end{figure}

Our first test in this subsection is to verify the conservation properties of the proposed methods $\textnormal{\bf Scheme-1}$ and $\textnormal{\bf Scheme-2}$. In the numerical simulations, we use $CFL = 0.13$ (typical time step size $ \Delta t \approx 0.011$) for $\textnormal{\bf Scheme-1}$ and $ CFL = 5$ (typical time step size $ \Delta t \approx 0.41$) for 
$\textnormal{\bf Scheme-2}$.
 Figure \ref{figure_iaconserve} shows  the absolute value of the relative error of the total particle number and total energy for $\textnormal{\bf Scheme-1}$ and $\textnormal{\bf Scheme-2}$ with simulations parameters $S_1$. We  observe that the relative errors   stay small, below $ 10^{-11}$ for $\textnormal{\bf Scheme-1}$
 and below $10^{-9}$ for $\textnormal{\bf Scheme-2}$.
The errors of total energy for $\textnormal{\bf Scheme-2}$ are slightly  larger mainly due to the error in  the Gauss-Seidel  iteration relating to the preset  tolerance parameter $\epsilon = 10^{-11}$.

%\subsubsection{Anomalous resistivity and other collections of data}

One of the important quantity to consider is 
the anomalous resistivity $\eta$  defined as  \cite{davidson}
%\begin{align}
%\eta = \frac{m_e}{n_0e^2}\left(-\frac{1}{p_e}\frac{\partial p_e}{\partial t}\right)= \frac{m_e}{n_0e^2}\left(-\frac{1}{\bJ}\frac{\partial\bJ}{\partial t}\right)
%\end{align}
\begin{align}
\eta = \frac{m_e}{n_0e^2}\left(-\frac{1}{\bJ_0}\frac{\partial\bJ_0}{\partial t}\right)
\end{align}
%where $p_e=m_e\int\int v f_{0e} ddxdv$ is the electron momentum(measured in the rest frame of ions). $f_{0e}$ is the spatially-averaged distribution functions.
 The resistivity  $\eta$ 
in our simulations is scaled by $m_e\omega_{pe}/n_0e^2$ and calculated at each time step
using a first order backward finite difference method for the time derivative. More discussions about the calculation of resistivity can be found in
 \cite{thesision}.

\begin{figure}
\centering
\subfigure[$v_{de}=0.17$]{\includegraphics[width=0.45\textwidth]{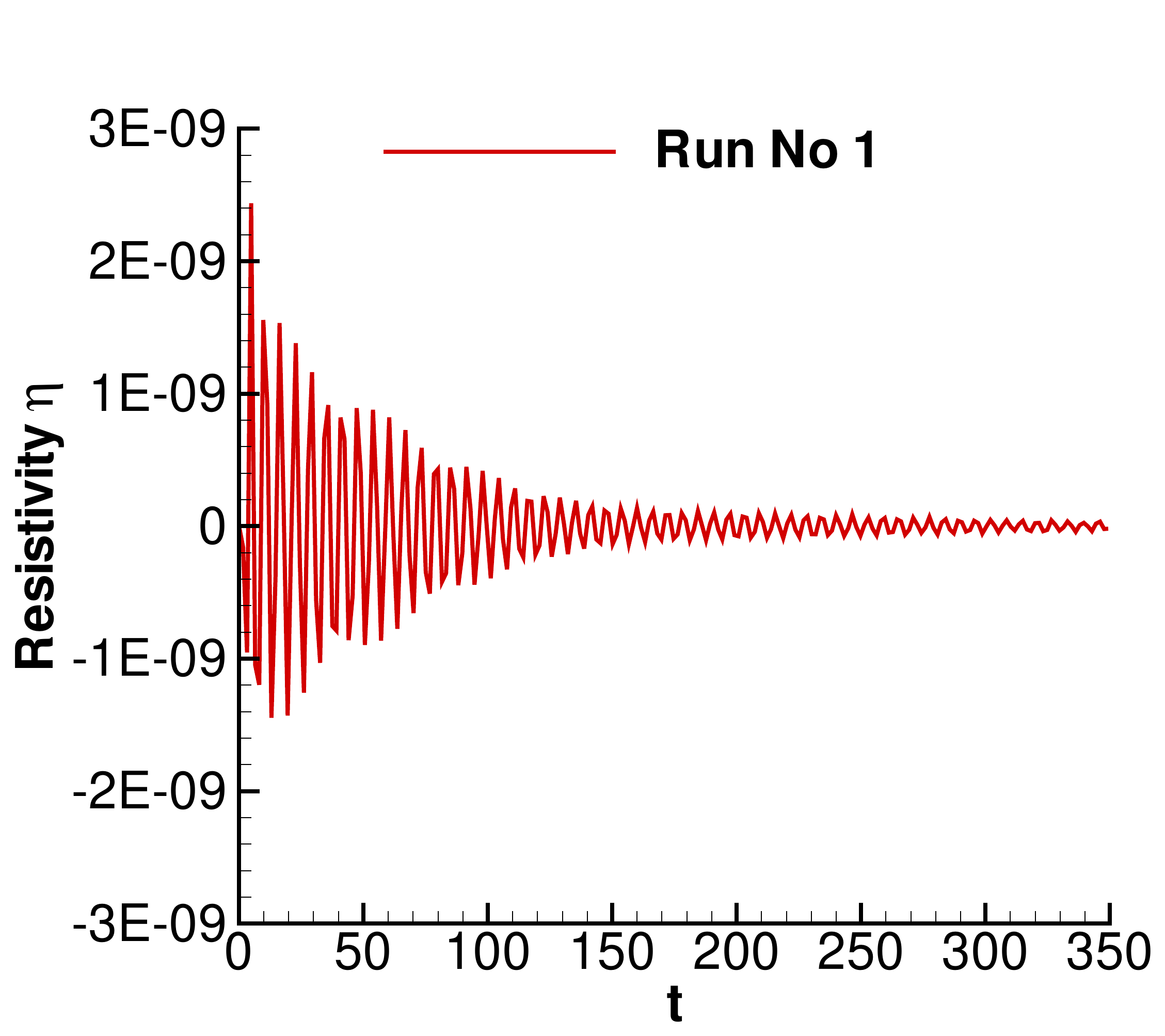}
\label{iaexplicitresistivity25:2}}
\hspace{1cm}
\subfigure[$v_{de}=1.7$]{\includegraphics[width=0.45\textwidth]{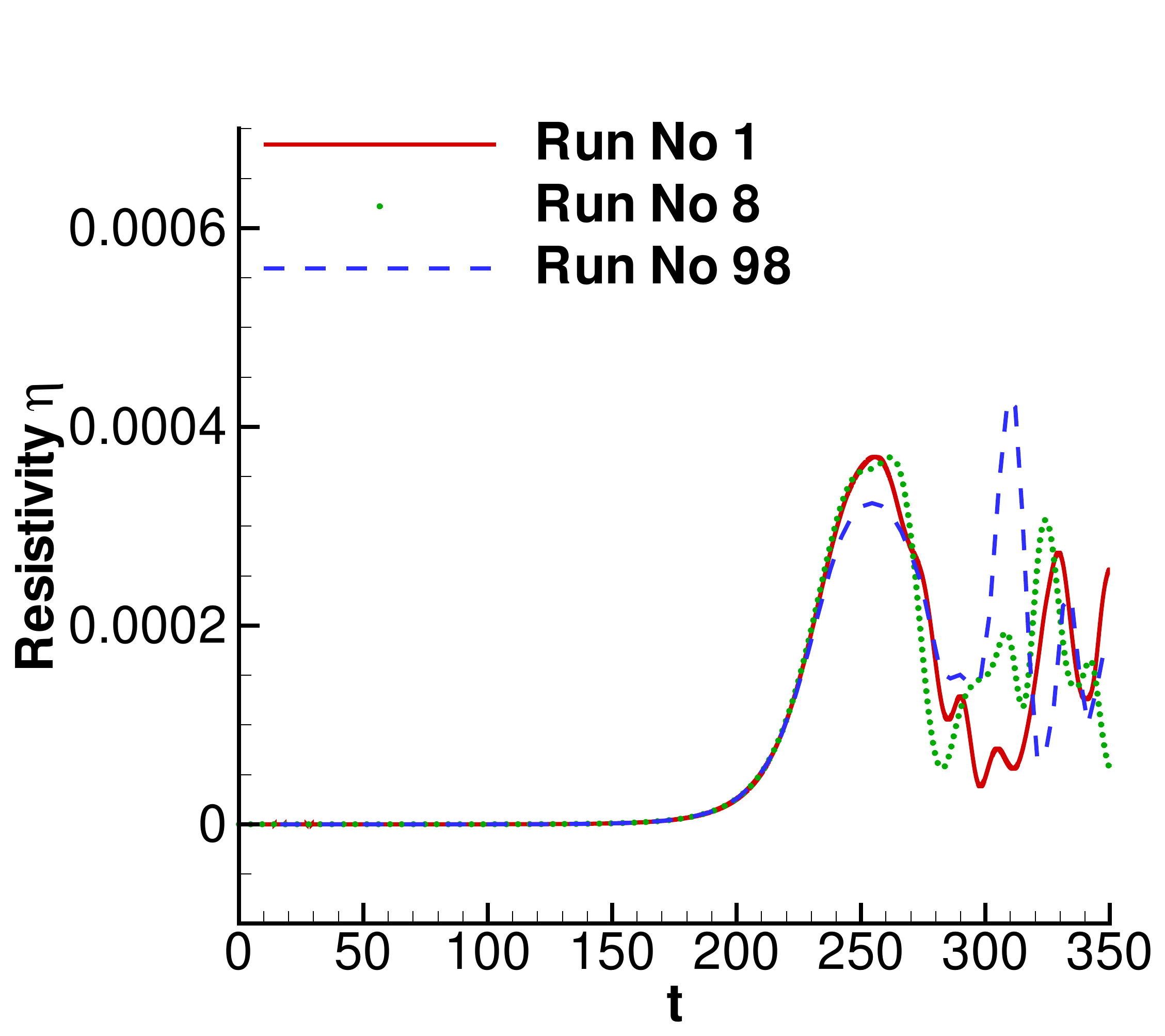}
\label{iaexplicitresistivity25:1}}
\caption{Time evolution of anomalous resistivity. $\textnormal{\bf Scheme-1}$.}
\label{iaexplicitresistivity25}
\end{figure}

We first perform a numerical test with $v_{de}=0.17$. This drift velocity is not large enough to trigger instability \cite{thesision} and the wave eventually got damped as illustrated in Figure \ref{iaexplicitresistivity25:2}. In the rest of the paper, we focus on the case of $v_{de}=1.7$, which is a drift velocity chosen large enough to result in the ion-acoustic instability.
Similar to \cite{ionacoustic}, to study the impact of the initial random phases $\varphi_n$, we perform an ensemble of 100 simulations  with random $\varphi_n$
%with different  sets of initial white noise  
  using the explicit method $\textnormal{\bf Scheme-1}$. For each of the 100 simulations, the 
phases of the initial white noise were randomly picked out of the uniform distribution on $[0,2\pi]$.

\begin{figure}[!htbp]
\centering
\subfigure[Electric filed mode] {\includegraphics[width = 0.45\textwidth]{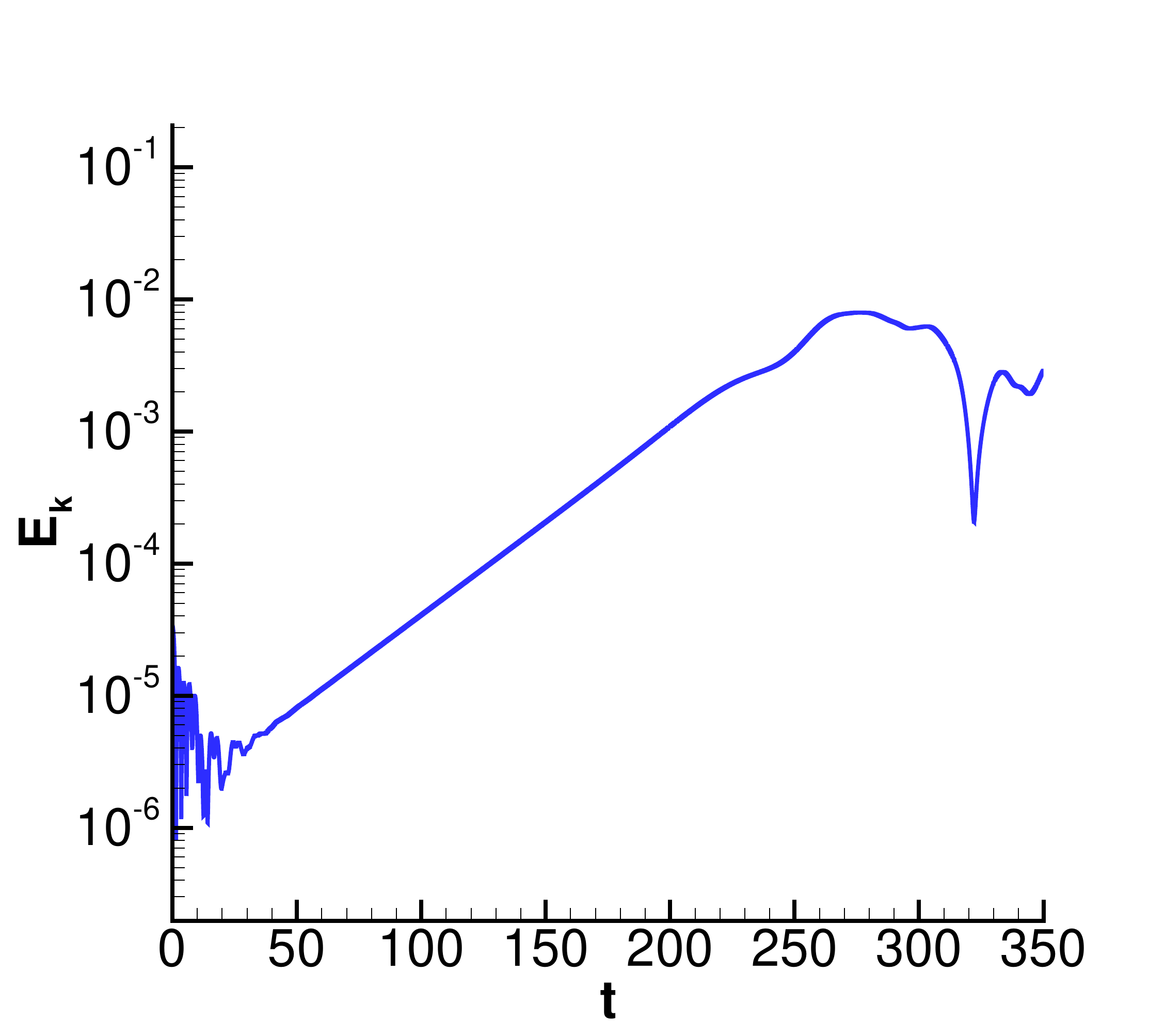}
\label{iaexplicitEfm25}}
%\subfigure[] {\includegraphics[width = 0.45\textwidth]{iaexplicitelectric25.eps}}
\subfigure[ Spatially averaged  distribution functions ] {\includegraphics[width = 0.45\textwidth]{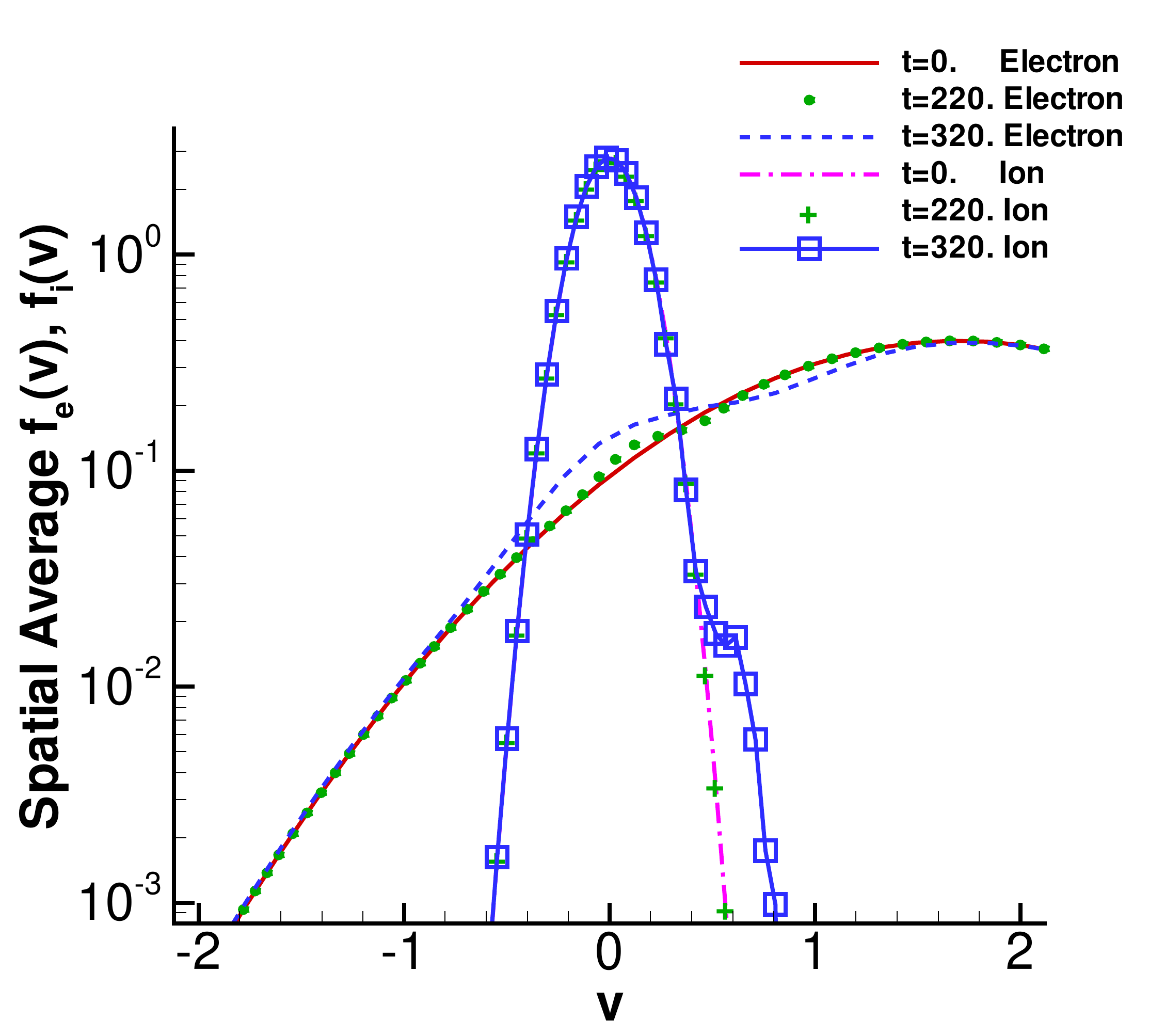}
\label{iaexplicitspatial25}}
\caption{  $\textnormal{\bf Scheme-1}$.  Run No. 1. (a): Time evolution of fastest growing electric filed mode ($\kappa = 0.501$).  (b): Spatially averaged electron and ion distribution functions.}
\end{figure}

Figure \ref{iaexplicitresistivity25:1} shows the time evolution of anomalous resistivity for three representative  simulations Run No. 1, 8, 98.  All the simulation runs  of ion-acoustic waves show similar resistivity evolution, while the exact values of  the anomalous resistivity differ from one simulation to the another, due to the initial fluctuations in the random phases. As in
\cite{ionacoustic}, we also mark different periods of the evolution of the resistivity using four regimes.  This can be explained  by comparing with Figure \ref{iaexplicitEfm25} where the fastest-growing mode of Run No. 1 is plotted. 
\begin{itemize}
\item  $0 \lesssim t \lesssim 50$: the initial onset. The resistivity is negligible during this period, when the wave amplitudes are small and the initial fluctuations are due to the ballistic ``free streaming" solutions to the Vlasov equation \cite{bookion, thesision}.
\item $50 \lesssim t \lesssim 200$: the linear regime. The anomalous resistivity remains close to zero and the only growing modes in the system are the linear modes.  During this time period, the fastest-growing mode grows exponentially in time as predicted  by the linear theory.%, showing in Figure \ref{iaexplicitEfm25}.
\item $200 \lesssim t \lesssim 240$: the quasi-linear regime. The resistivity rises rapidly for $t\gtrsim200$ to a peak at $t\sim 240$ and the fastest-growing mode starts deviating from exponential growth and saturates.
\item $t \gtrsim 240$: the nonlinear regime. The resistivity is relatively stable, oscillating  about the resistivity level reached at the end of the quasi-linear regime. \end{itemize}

\begin{figure}[!htbp]
\centering
%\subfigure[$t=0$. Electron distribution]{\includegraphics[width=0.45\textwidth]{iaexplicitecontourt=0_1.png}}
%\subfigure[$t=0$. Ion distribution]{\includegraphics[width=0.45\textwidth]{iaexpliciticontourt=0_1.png}}
\includegraphics[width=0.75\textwidth]{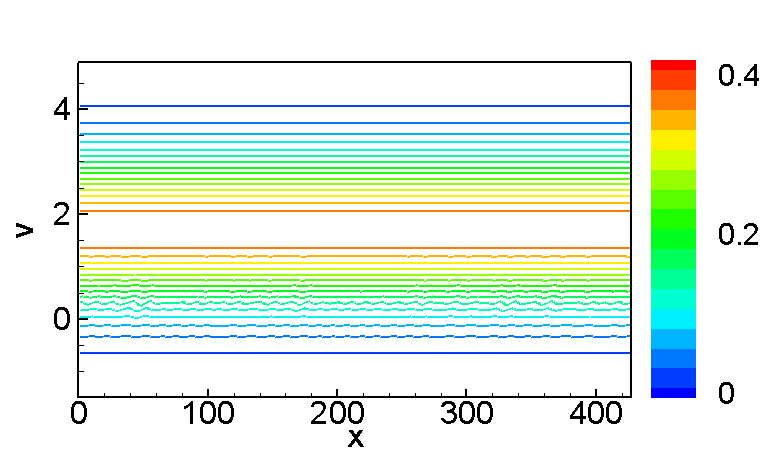}
%\subfigure[$t=130$. Ion distribution]{\includegraphics[width=0.45\textwidth]{iaexpliciticontourt=130_1.png}}\
\includegraphics[width=0.75\textwidth]{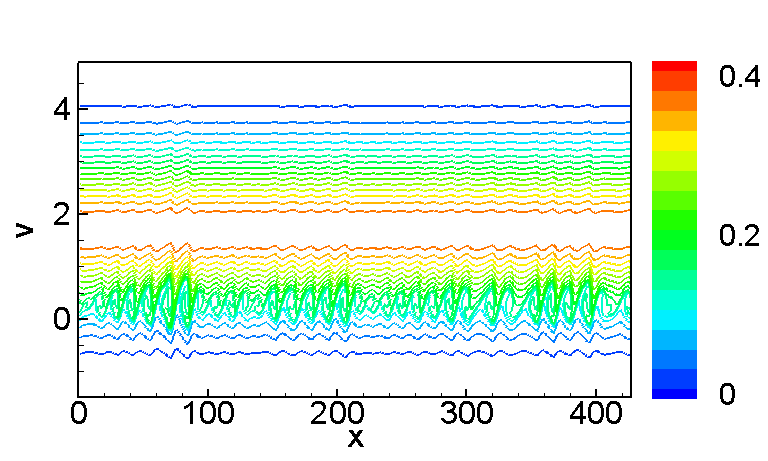}
%\subfigure[$t=220$. Ion distribution]{\includegraphics[width=0.45\textwidth]{iaexpliciticontourt=220_1.png}}
\includegraphics[width=0.75\textwidth]{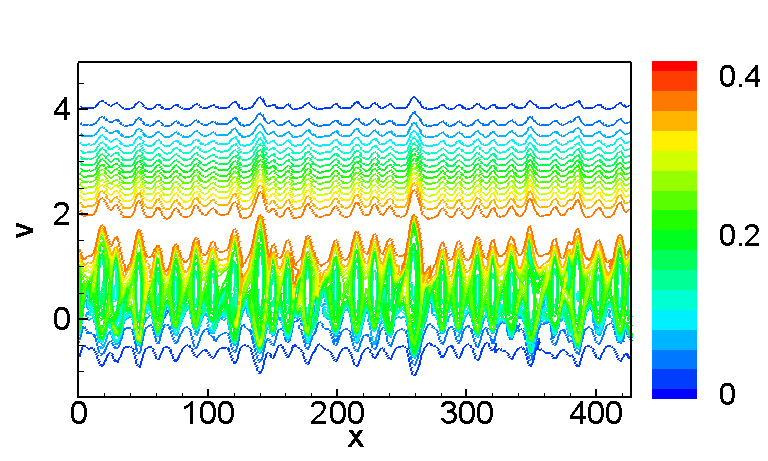}
%\subfigure[$t=320$. Ion distribution]{\includegraphics[width=0.45\textwidth]{iaexpliciticontourt=320_1.png}}
\caption{Contour plots of the electron distributions at  time $t=130$ (top), $t=220$ (middle) and $t=320$ (bottom). $\textnormal{\bf Scheme-1}$. Run No. 1.}
\label{figure_contour1}
\end{figure}

\begin{figure}[!htbp]
\centering
%\subfigure[$t=0$. Electron distribution]{\includegraphics[width=0.45\textwidth]{iaexplicitecontourt=0_1.png}}
%\subfigure[$t=0$. Ion distribution]{\includegraphics[width=0.45\textwidth]{iaexpliciticontourt=0_1.png}}
\includegraphics[width=0.75\textwidth]{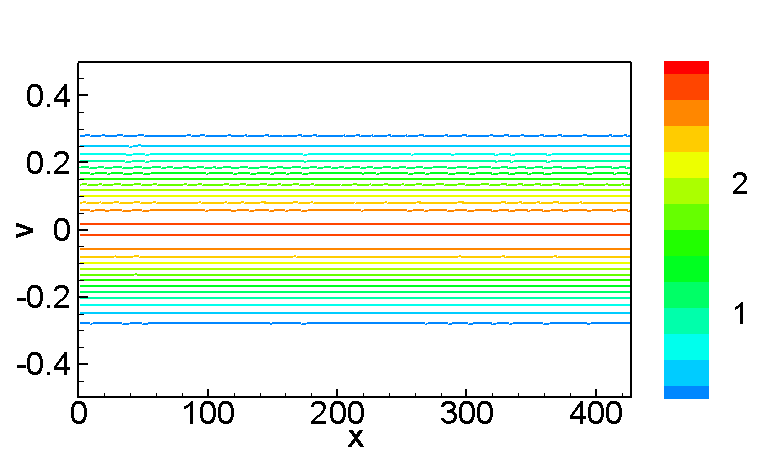}
%\subfigure[$t=130$. Ion distribution]{\includegraphics[width=0.45\textwidth]{iaexpliciticontourt=130_1.png}}\
\includegraphics[width=0.75\textwidth]{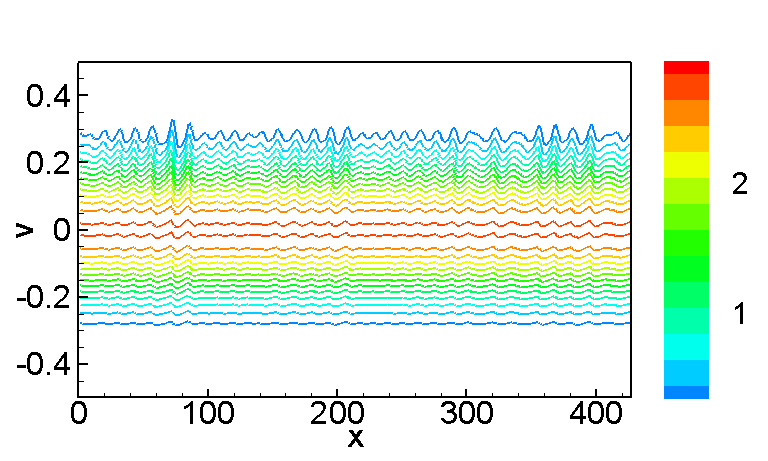}
%\subfigure[$t=220$. Ion distribution]{\includegraphics[width=0.45\textwidth]{iaexpliciticontourt=220_1.png}}
\includegraphics[width=0.75\textwidth]{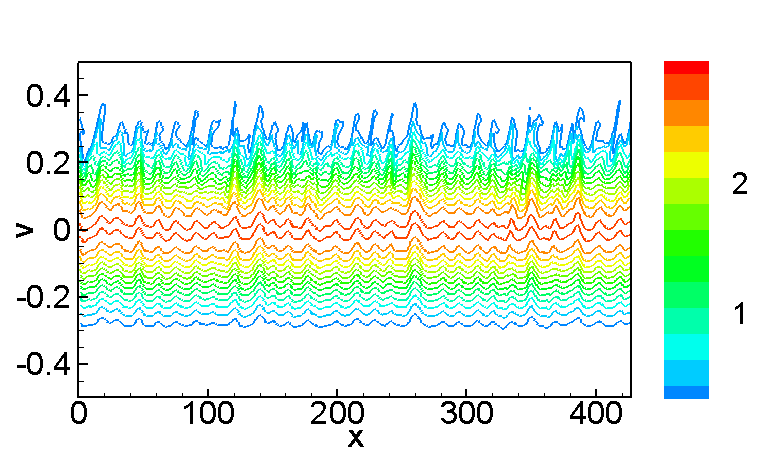}
%\subfigure[$t=320$. Ion distribution]{\includegraphics[width=0.45\textwidth]{iaexpliciticontourt=320_1.png}}
\caption{Contour plots of the ion distributions at  time $t=130$ (top), $t=220$ (middle) and $t=320$ (bottom). $\textnormal{\bf Scheme-1}$. Run No. 1.}
\label{figure_contour2}
\end{figure}

In Figure \ref{iaexplicitspatial25}, we plot the spatially averaged distributions functions of ions and electrons at selected time for Run No 1. We zoom in velocity space and $f_{\alpha}$  in order to see the details.  The two distributions are shown at the beginning to the simulation ($t=0$), near quasi-linear saturation $t=220$, and $t=320$ which is during the nonlinear regime. We obtain qualitatively similar results to \cite{ionacoustic}. Namely, the development of the plateau formation for the distributions functions since the quasi-linear regime indicates the momentum exchange between electrons and ions via the ion-acoustic waves.
% The arrival at quasi-linear saturation can also be identified by the plateau formation at the range of resonant velocities in the spatially averaged electron distribution function \cite{Petkaki2003}. At $t=320$, the plateau in the spatially averaged electron distribution function has spread to lower and higher velocities outside the lear resonant region. The spatially averaged ion distribution shows a high-energy tail, which suggests that momentum exchange has taken place between electrons and ions via the ion-acoustic waves. 
In Figure \ref{figure_contour1} and \ref{figure_contour2}, we plot the probability distribution functions of electron and ion at $t=130, 220, 320$ for Run No. 1. We have zoomed in velocity space in order to see the detail structure of the solutions. At $t=130$ (a time frame in the linear regime ), the electron and ion distributions still stay rather close to the initial distributions. At $t=220$ (a time frame in the quasi-linear regime), deviations from the initial configurations are visible. In particular for the electron distributions, small trapping regions start to form. This is more prominent at $t=320$ (a time frame in the nonlinear regime) and several large trapping islands are displayed in the electron distribution functions. Such observations are consistent with the depiction in Figure \ref{iaexplicitspatial25}.

%To demonstrate that  the variation of the anomalous resistivity in the nonlinear regime is also do to electron and ion bounce motion, we show in 
%  the electron and ion distribution functions and we have zoomed in velocity space in order to see the details. 
%At wpet = 990 a very small trapping
%region has started to form at the resonant electron velocity
%region. By wpet = 1485 the trapping region has expanded in
%velocity space and several trapping islands are present in
%phase space diagram. Coalescence of some islands is visible
%too. Using either the average wavelength of the island
%structures seen in Figure 9 or the wave number of the
%largest amplitude wave mode, the characteristic electron
%trapping wave number is k  0.09  0.14 m1. At wpet =
%990 the ion distribution function is undisturbed yet by the
%presence of the growing electric field modes. By wpet =
%1485 ion trapping regions are beginning to form with a
%similar wave number.

\begin{figure}[!htbp]
\centering
\subfigure[$t=0$.]{\includegraphics[width=0.45\textwidth]{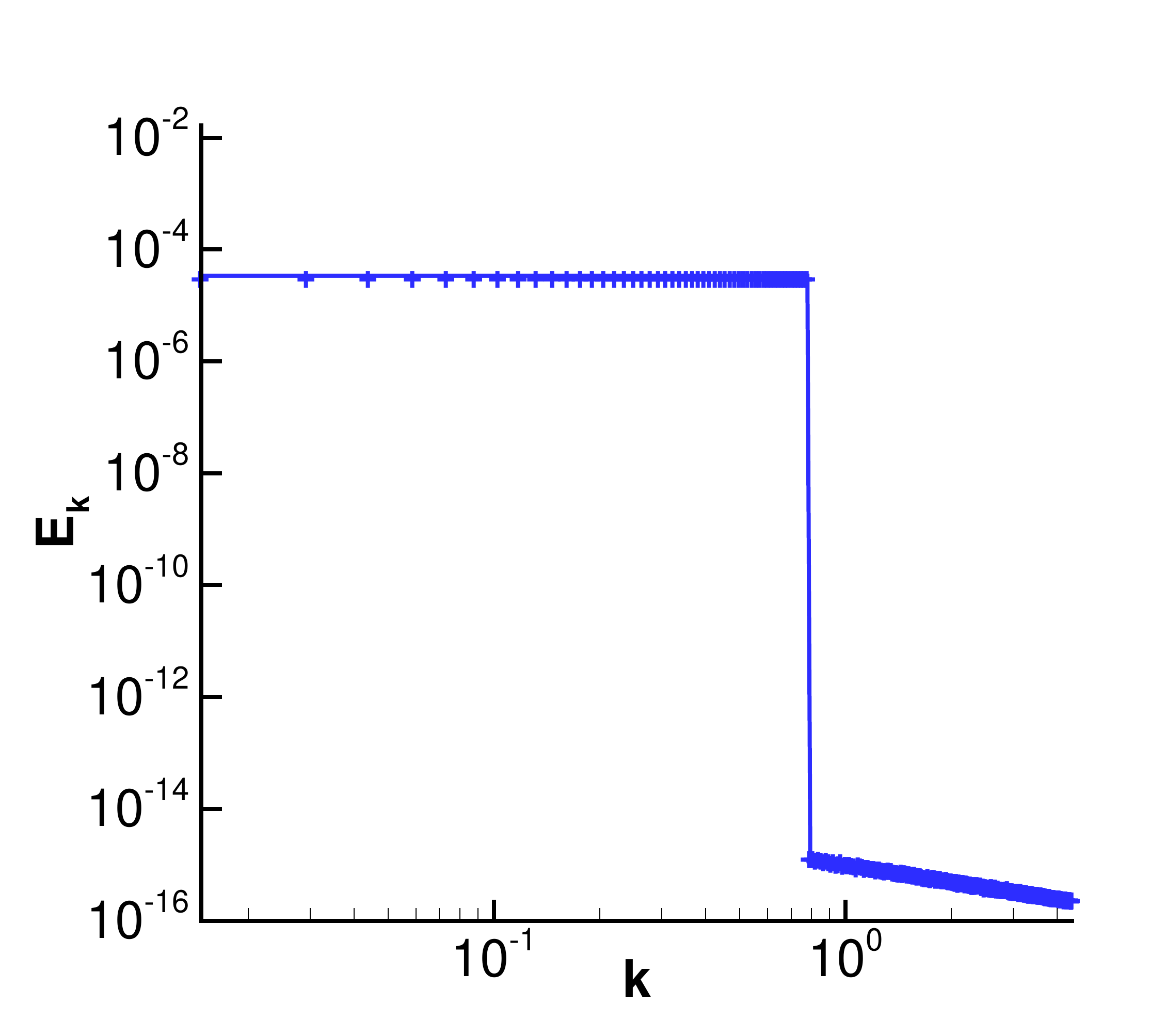}}
\subfigure[$t=130$.]{\includegraphics[width=0.45\textwidth]{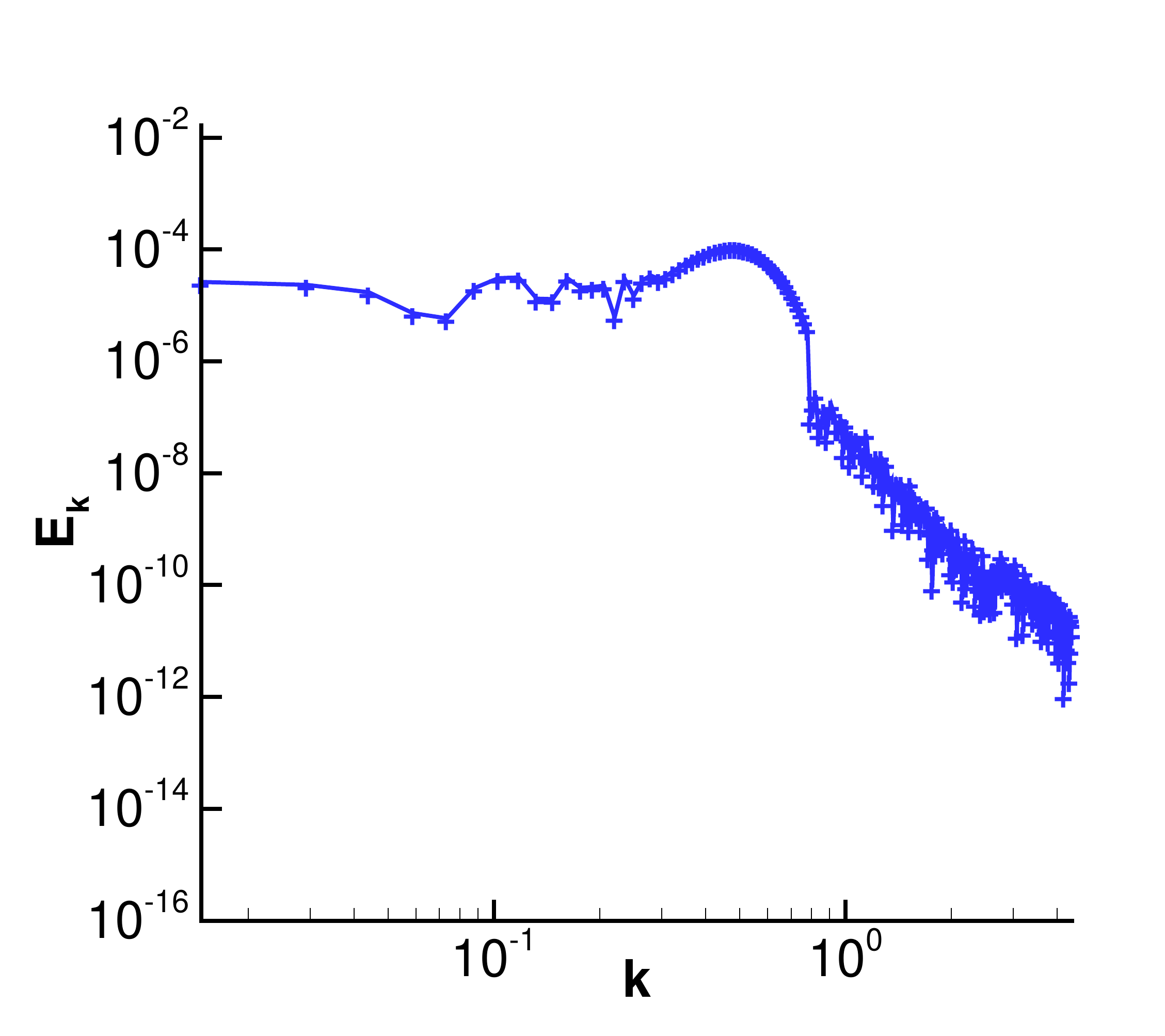}}\\
\subfigure[$t=220$.]{\includegraphics[width=0.45\textwidth]{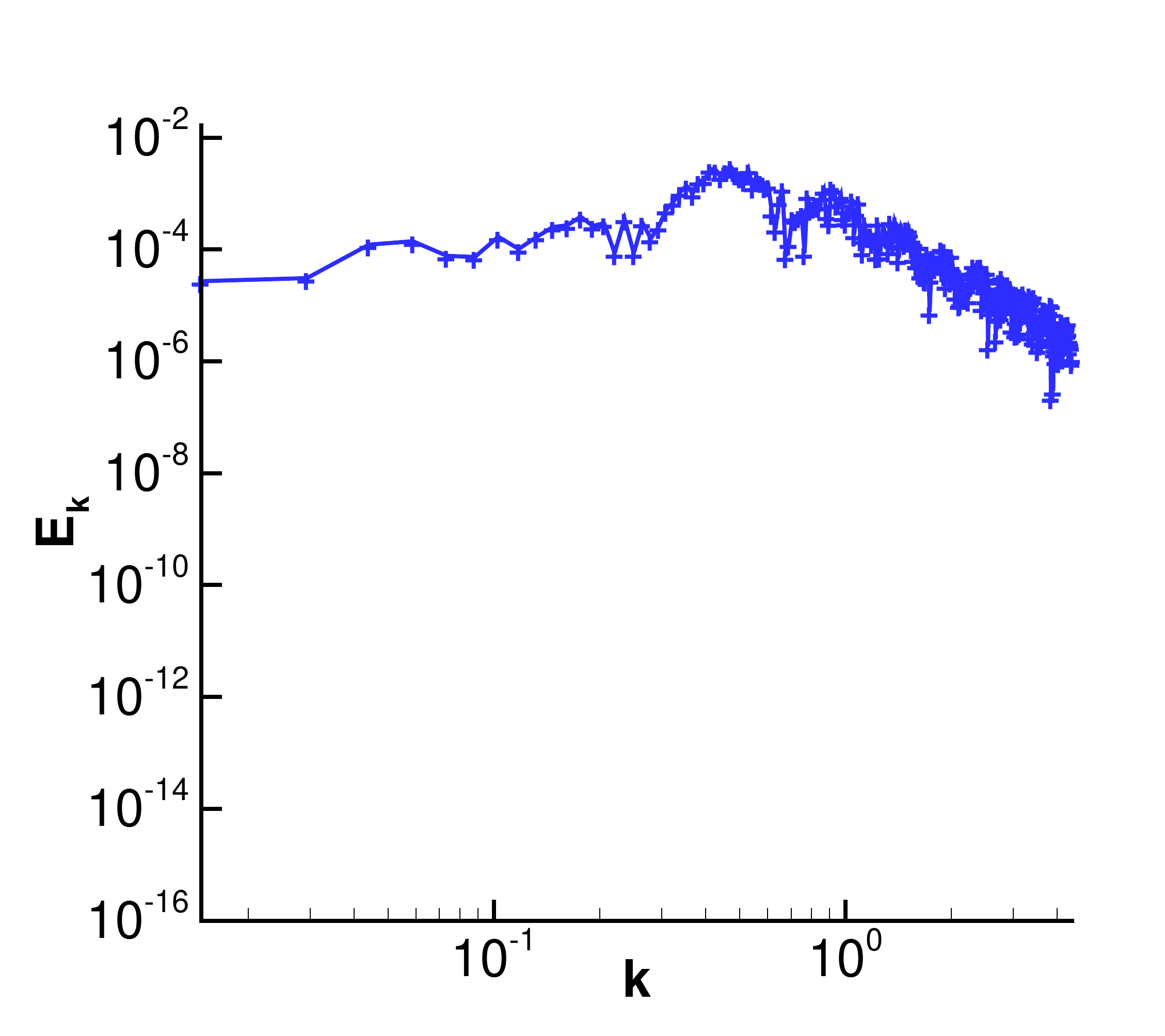}}
\subfigure[$t=320$.]{\includegraphics[width=0.45\textwidth]{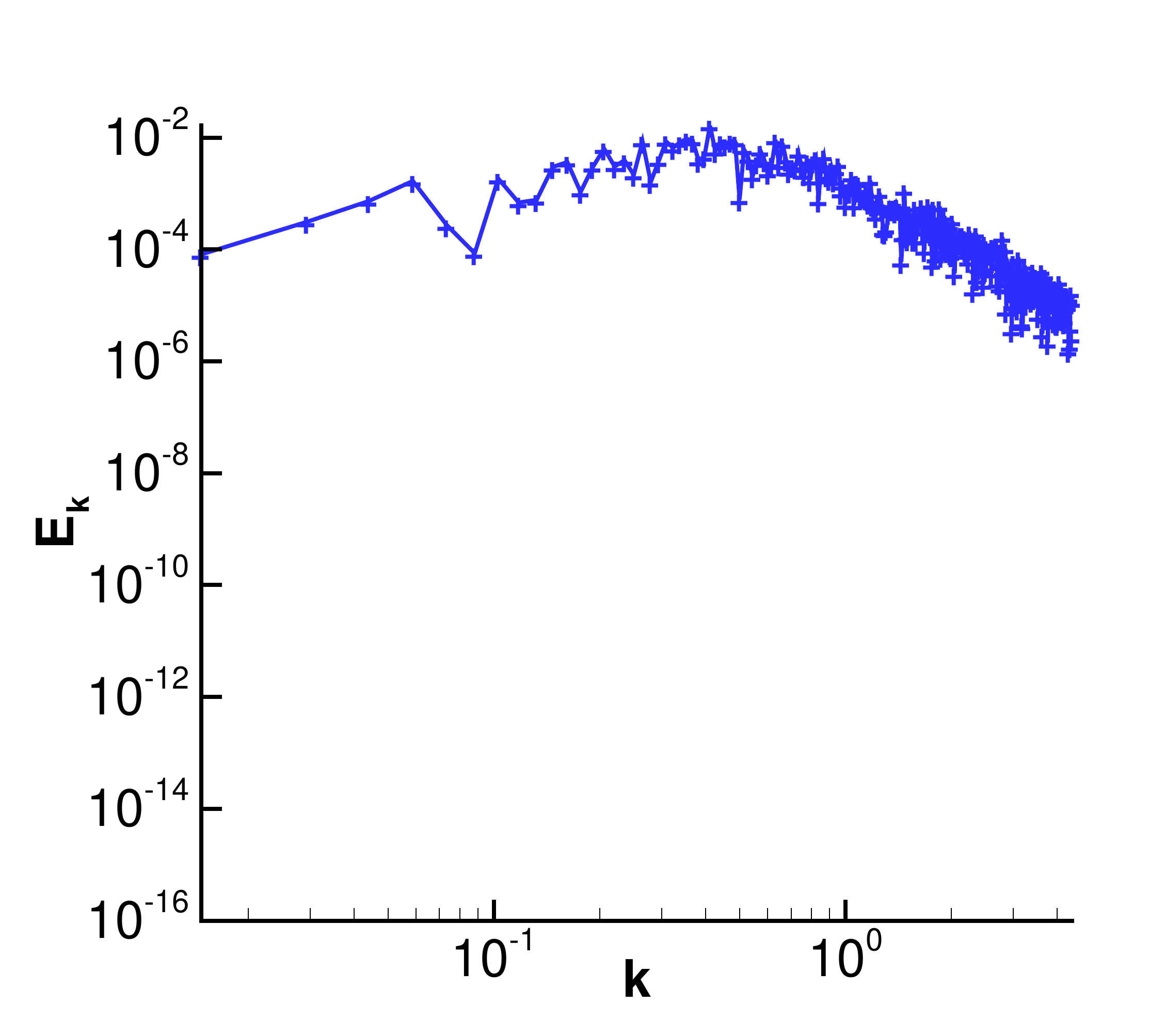}}
\caption{Electric field spectrum in $k$-space($E_k-k$) at four different times. $\textnormal{\bf Scheme-1}$. Run No. 1.}
\label{iaexplicitelectricspectrum}
\end{figure}

 To further show the details of the solutions at those given times, we  plot the electric spectrum in $\kappa$-space in Figure  \ref{iaexplicitelectricspectrum}.  %at the corresponding time frames. %during the beginning, linear, quasi-linear, nonlinear regimes of our simulation Run No. 1. 
 At t=0, we show the electric field spectrum associated with the initial condition, notice for $k>\kappa_{N_{max}}$, the nonzero values are due to the double precision accuracy used in the numerical simulations. For all later times, the spectrum demonstrate similar results as in \cite{ionacoustic}.
At time $t=130$, the wave power is above the initial noise field.  %{\color{red} compare with linear???? need graph?}. 
Since the quasilinear regime, we can observe that the wave energy cascades into wave modes outside the linear resonant region, and generally all modes have increased in power. The development of the power law electric field spectrum is evidence of nonlinear wave-wave coupling.
%At $t=220$, wave energy cascades into wave modes outside the linear resonant region, both at harmonics of the linear modes and not, and generally all modes have increased in power. At $t=320$, the electric field spectrum evolves toward a fully developed forward Iand possibly also inverse) cascade with a power law inertial range from the injection frequency of the fastest-growing linear wave mode to the highest (lowest) $\kappa$-mode. At $t=350$, we stop the simulation since we begin to observe aliasing of the highest $\kappa$-modes to lower $\kappa$-modes. The development of the power law electric field spectrum is evidence of nonlinear wave-wave coupling.

\begin{figure}[!htbp]
%\subfigure[] {\includegraphics[width = 0.45\textwidth]{iaexplicitEfm25.eps}}\\
\subfigure[] {\includegraphics[width = 0.49\textwidth]{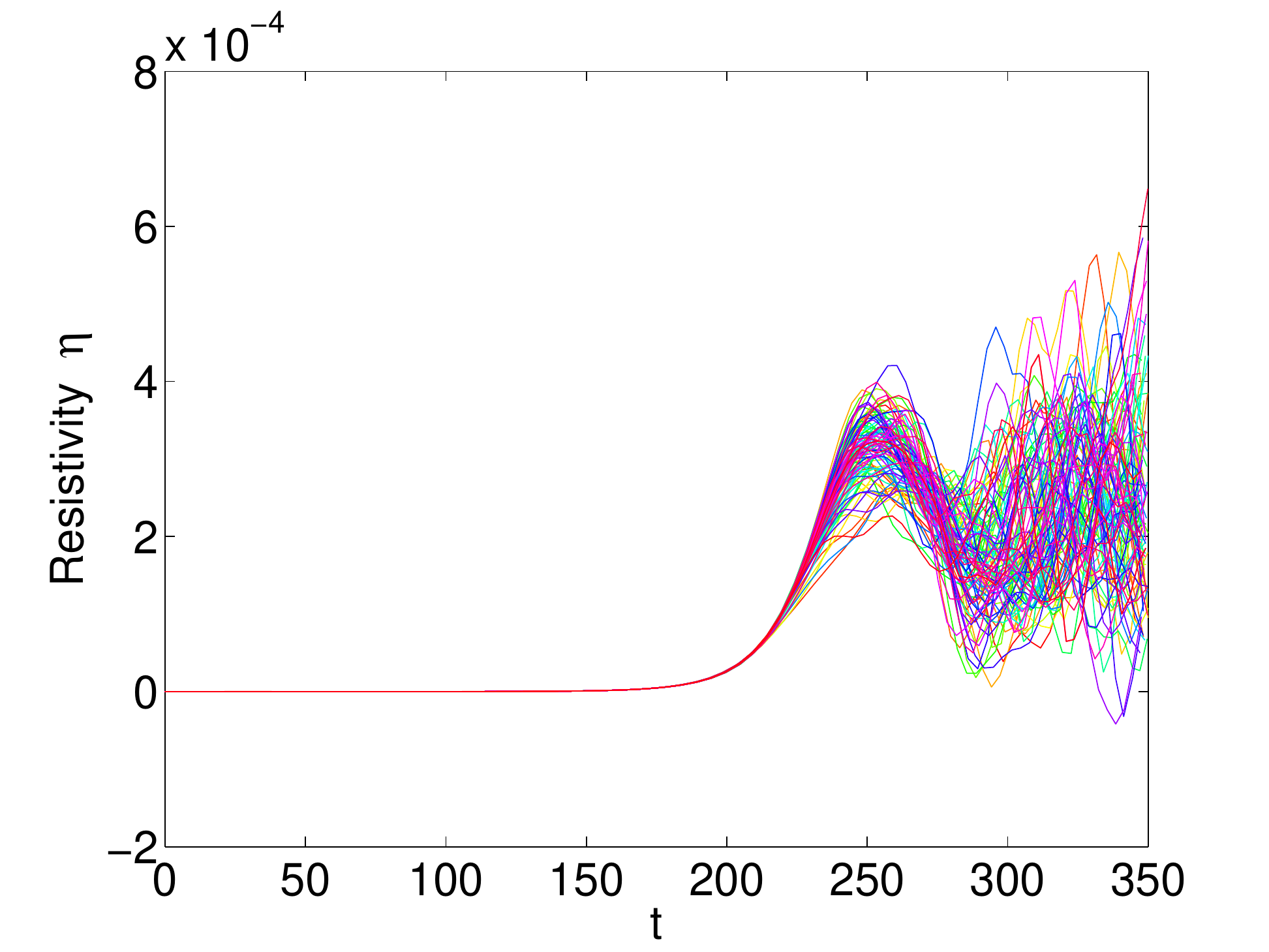}
\label{iaexplicitresistivity104}}
%\subfigure[] {\includegraphics[width = 0.45\textwidth]{iaexplicitspatial25.eps}}
\subfigure[] {\includegraphics[width = 0.49 \textwidth]{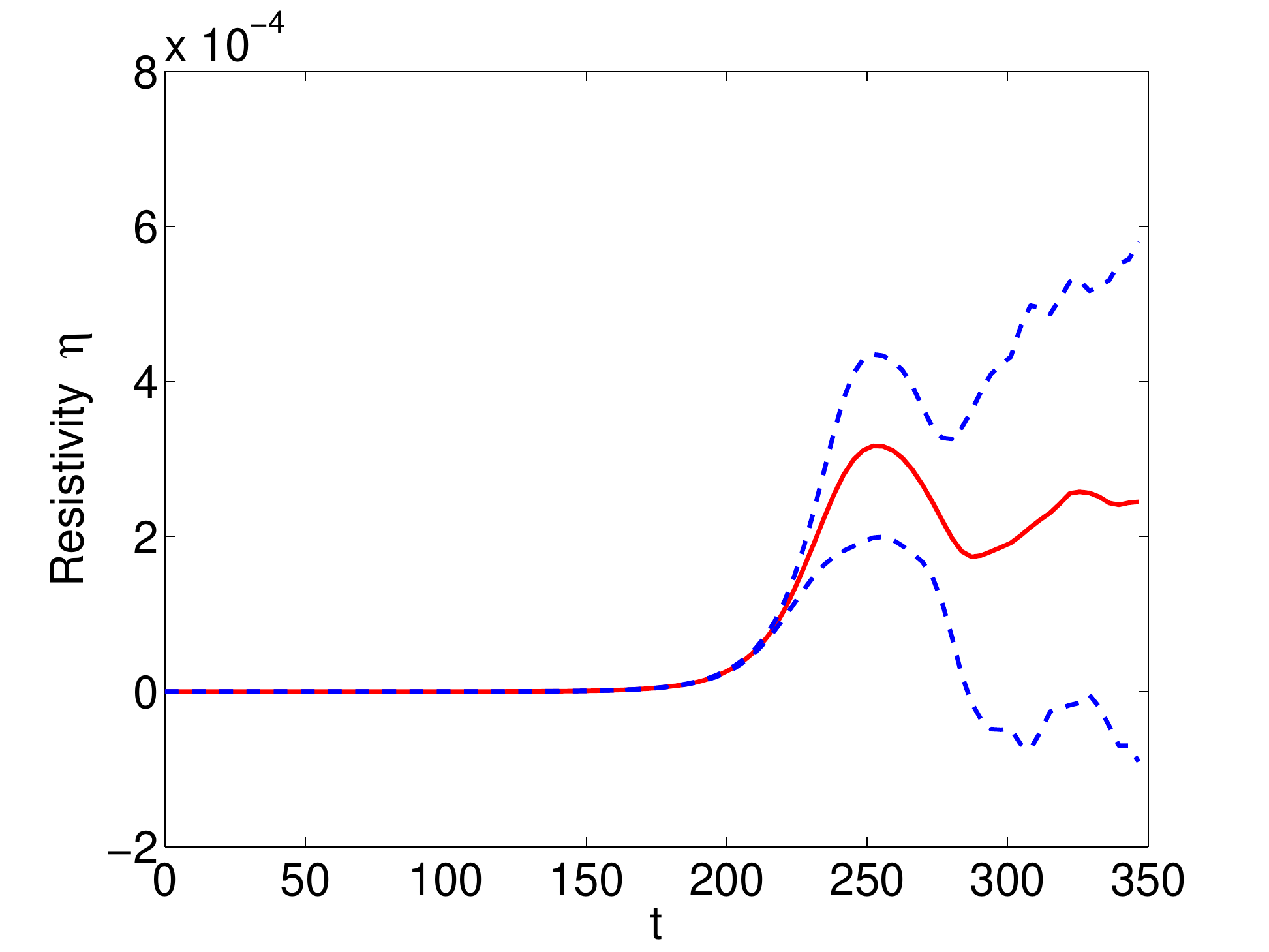}
\label{iaexplicitresistivitymeandev}}
\caption{(a): Overplotted time evolutions of  ion-acoustic resistivity from 100 simulations. (b) Mean value (solid line) and $\pm 3\sigma$ curves (dashed line) of the ensemble  ion-acoustic resistivities.}
\end{figure}

To verify the impact of the initial random phase field, next we perform statistical study of the 100 VA simulations similar to \cite{ionacoustic}.
 We notice that due to the CFL conditions of $\textnormal{\bf Scheme-1}$, each simulation has slightly different time step sizes. To benchmark the meaningful quantitative statistical analysis with \cite{ionacoustic}, we interpolate all the simulations with piecewise cubic Hermite interpolation. The time step and resistivity values  discussed in the following context  are the values after interpolation. The overall conclusion is similar to \cite{ionacoustic}. Figure \ref{iaexplicitresistivity104} over plot the time evolution of all 100 anomalous resistivities,  which proves the similarity and diversity in the behavior of the resistivity.  In Figure \ref{iaexplicitresistivitymeandev}, we plot the mean value of the resistivity calculated by averaging the value of the 100 resistivity values at each time step, and the $\pm$ three standard deviations from the mean. Comparing 
\ref{iaexplicitresistivity104} and \ref{iaexplicitresistivitymeandev}, we see that the range of resistivity values is well confined in $\pm 3\sigma$ of the mean, as would be expected by a Gaussian distribution. 

\begin{figure}[!htbp]
\centering
%\subfigure[] {\includegraphics[width = 0.45\textwidth]{iaexplicitEfm25.eps}}\\
\subfigure[$t=122-126$] {\includegraphics[width = 0.49\textwidth]{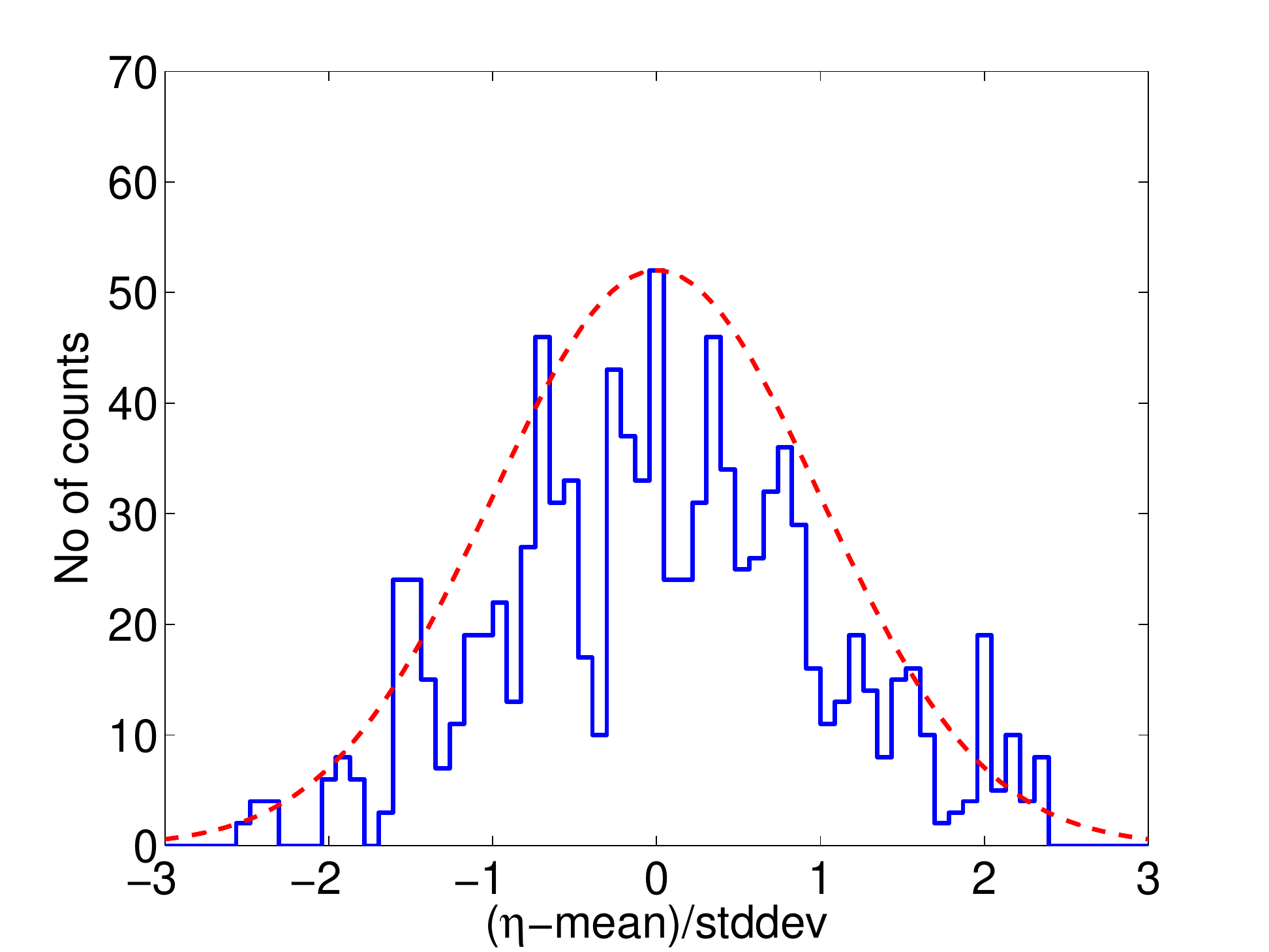}
\label{iaexplicithist:1}}
\subfigure[$t=220-224$]{\includegraphics[width = 0.49\textwidth]{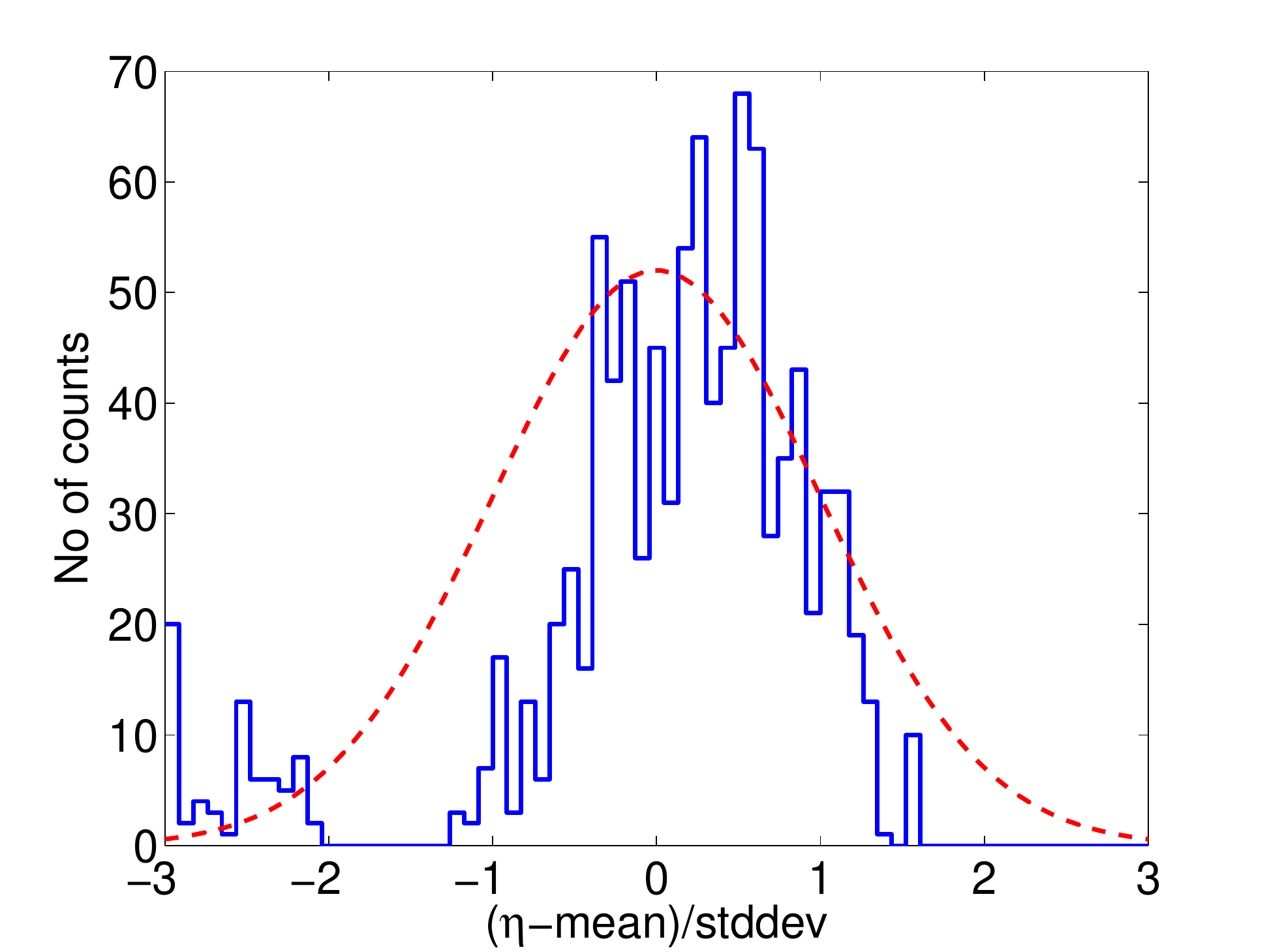}
\label{iaexplicithist:2}}\\
\subfigure[$t=283-287$] {\includegraphics[width = 0.49\textwidth]{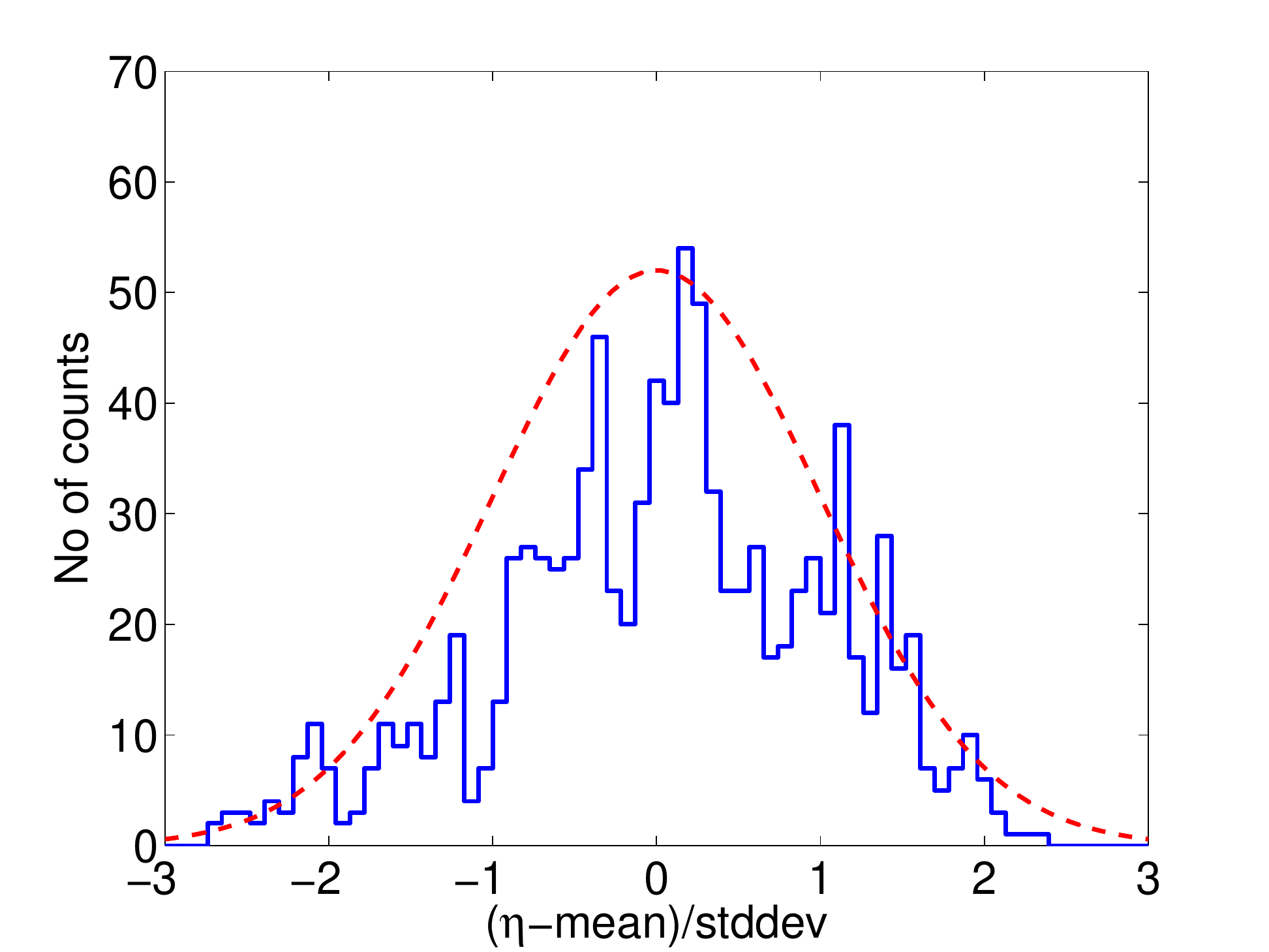}
\label{iaexplicithist:3}}\\
%\subfigure[] {\includegraphics[width = 0.45\textwidth]{iaexplicitspatial25.eps}}
\caption{Probability distribution of standardized  resistivity values ($(\eta(t) - \overline{\eta(t)})/\sigma $) during three time periods:
(a) linear regime, (b) quasi-linear regime, and (c) nonlinear regime.  Dashed line is the plot of Gaussian of mean 0 and standard deviation 1.}
\label{iaexplicithist}
\end{figure}

 To investigate this further, we plot histograms of the  probability distribution of the standardized resistivity values in Figure \ref{iaexplicithist} for three time periods from three different regimes  in the evolution of the instability to study how well they fit into the Gaussians. The standardized value of $\eta$ at $t$ is $(\eta(t)-\overline{\eta(t)})/\sigma(t)$, where $\overline{\eta(t)}$ is the ensemble mean value of $\eta$ at $t$, and $\sigma(t)$ is the standard deviation of $\eta$ at $t$. 
Each histogram comprises of all the standardized resistivity values from 10 consecutive time steps. 
Figure \ref{iaexplicithist:1} shows  the  probability distribution of the standardized resistivity values for $t=122-126$ during the linear regime of the instability. The distribution appears to fit a Gaussian reasonably well. 
Figure \ref{iaexplicithist:2} shows  the  probability distribution of the standardized resistivity values for $t=220-240$ during the quasi-linear regime of the instability. The distribution is sharply peaked and with a left tail longer than the right, which seems deviate from a Gaussian.
Figure \ref{iaexplicithist:3} shows  the  probability distribution of the standardized resistivity values for $t=283-287$ during the nonlinear regime of the instability. The distribution appears to be symmetric and approximately Gaussian.

%It is impractical to perform 104 ensemble runs for much higher mass ratio. However we performed one run with the real mass ration $S_2$ in order to demonstrate that the
%main physical processes are not changed. 
%Piecewise cubic Hermite interpolation
%plots of eta
%
%plot of Ek(t) ({\color{red} k =0.532})
%
%plots of f spaave
%
%plots of f
%
%all our results agrees well with \cite{ionacoustic}?????????/Users/zhongxh/Dropbox/paper/VlasovAmpere/Two species/paper/papervats.pdf
\begin{figure}[!htbp]
%\subfigure[] {\includegraphics[width = 0.45\textwidth]{iaexplicitEfm25.eps}}\\
\subfigure[Skewness] {\includegraphics[width = 0.49\textwidth]{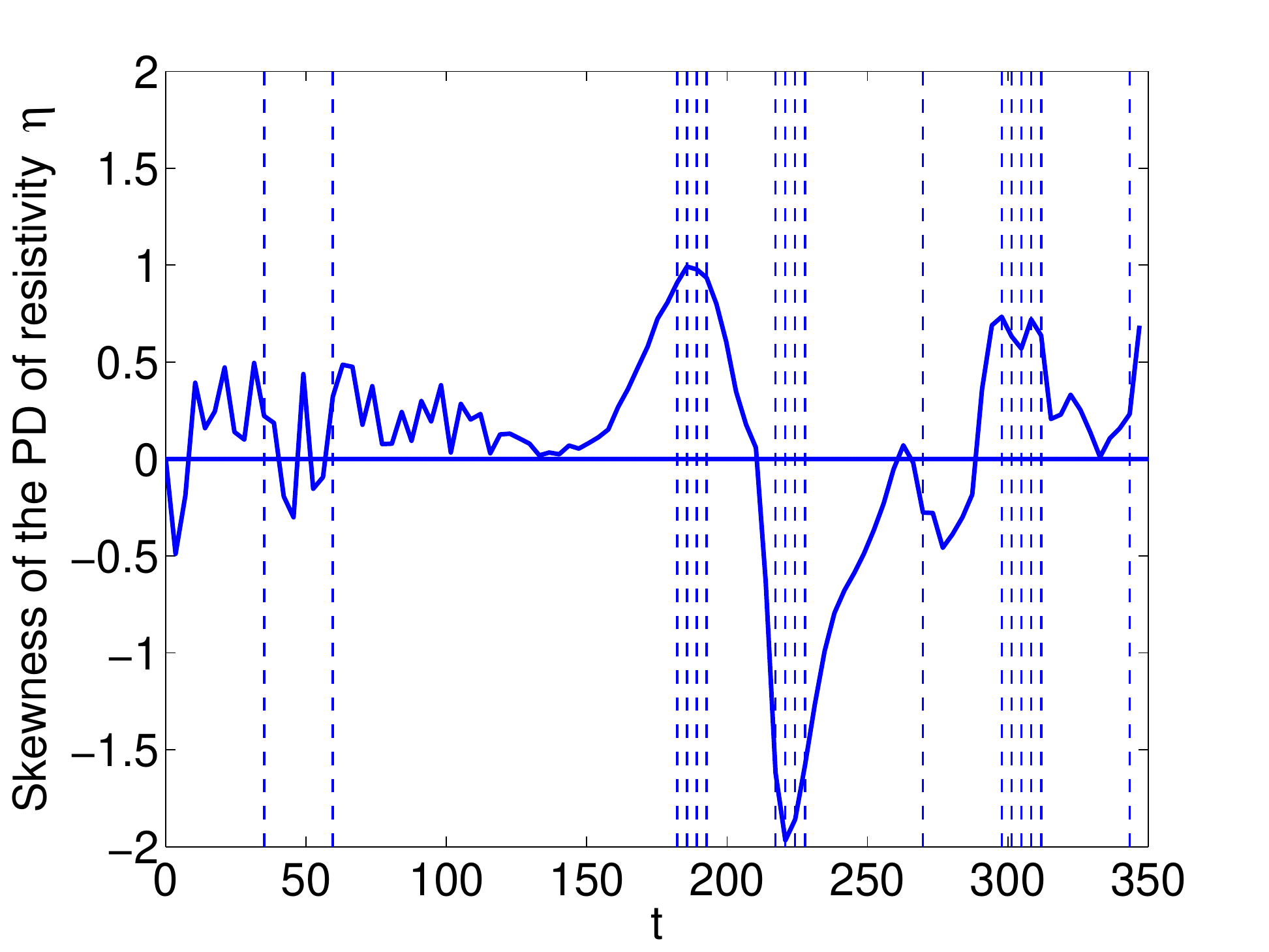}}
\subfigure[Kurtosis] {\includegraphics[angle =90,width = 0.49\textwidth ]{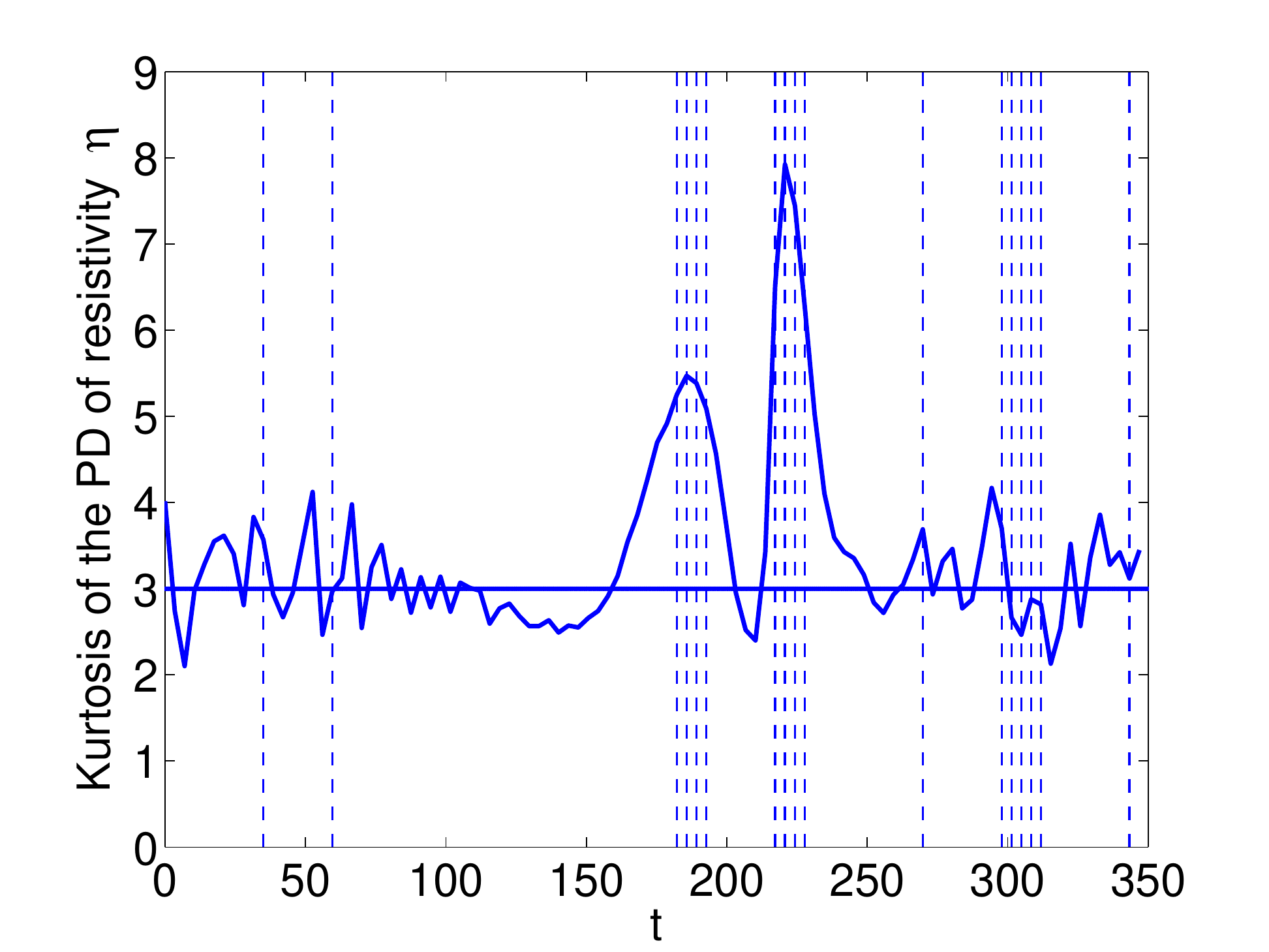}}
%\subfigure[] {\includegraphics[width = 0.45\textwidth]{iaexplicitspatial25.eps}}
\caption{Time evolution of the skewness and the kurtosis of the probability distribution of ensemble resistivity values.}
\label{iaexplicitskewkurt}
\end{figure}

As in \cite{ionacoustic}, we  perform a chi-square test for the standardized resistivity values at  each time step, where we tested the goodness of fit to a Gaussian
distribution of mean zero and unit standard deviation of each one of the probability distributions of standardized  resistivity
values at the 0.05 and 0.01 significance level. Overall, the test fails in 16.8\% of the time steps at the 0.05 level and at 7.1\% at the 0.01 level, (compared to 8.4\% and 1.8\% in \cite{ionacoustic}) which suggests that the distributions of resistivity values may not fit a Gaussian at all times. Figure \ref{iaexplicitskewkurt} shows the times when chi-square test fails at the 0.05 level (dashed lines). Figure \ref{iaexplicitskewkurt} also shows the time evolution of the skewness and kurturtosis of the resistivity distributions. The skewness of the resistivity distributions remains close to 0 for most of the time evolution, implying that the distributions are mostly symmetrical. The kurtosis of the distributions remains close to 3 for most to the time evolution, consistent with a Gaussian. We also observe that  the amplitude of the oscillations in the values of skewness and kurtosis increases as we approach quasi-linear saturation and beyond. This deviation from the expected Gaussian values can also be seen in Figure \ref{iaexplicithist:2}: the skewness is negative, which indicates that the left tail is longer; the kurtosis is well above three indicating a more sharply peaked distribution than a Gaussian. Although the exact values of skewness and Kurtosis are different from \cite{ionacoustic}, the overall shape and the qualitative behavior remains the same. We remark that our plots seem smoother in time because of the Hermite interpolation we used when post-processing the data.

%Figure \ref{iaexplicitelectric25} shows snapshots of the electric field for $80<x<130$ at  three different times during the linear regime of the instability.
%In order to check the phase speeds of these waves, a rough estimate can be made
%by following the peak of a wave as it changes position in time. 
%Selecting the middle {\color{red}  phase speed
%
%
%growth rate
%
%or delete both?}
%peak (solid line, z = 74.5?1
%De), over the time period shown t = 10.1!?1
%pe , it moves
%z = 10.1?1
%De. Hence the estimate of the phase speed is vph = z/t = 3.4?105 ms?1.
%The sound speed for the initial plasma parameters used in this run is 3.36?105 ms?1,
%so this estimate suggests that the ion-acoustic waves generated during the linear stage
%of the instability are traveling at the ion sound speed, as linear theory predicts.
%

%\begin{figure}[!htbp]
%\centering
%%\subfigure[] {\includegraphics[width = 0.45\textwidth]{iaexplicitEfm25.eps}}\\
%\subfigure[] {\includegraphics[width = 0.45\textwidth]{iaexplicitelectric25.eps}}
%%\subfigure[] {\includegraphics[width = 0.45\textwidth]{iaexplicitspatial25.eps}}
%\caption{Three snapshots of the electric field. $\textnormal{\bf Scheme-1}$. $\mu_i=1/25.  $}
%\label{iaexplicitelectric25}
%\end{figure}
%

\section{Concluding Remarks}
\label{sec:conclusion}

In this paper, we develop explicit and implicit energy-conserving Eulerian solvers for the   two-species VA system and apply the methods to simulate  current-driven ion-acoustic instability. The overall results show excellent conservation of the total particle number and total energy regardless of the mesh size as predicted by the theoretical studies. The implicit methods, though do not suffer from CFL restrictions, still require $w_{pe} \Delta t  \lesssim 1$ to fully resolve the electron kinetic effects. For the current-driven ion-acoustic instability, we perform an ensemble of 100 VA simulations with random phase perturbations to investigate the anomalous resistivity with a reduced mass ratio. The results agree well with previous studies. In future work, it would be interesting to generalize such schemes to simulate multi-species systems when the electron kinetic effects are of less importance. A multiscale algorithm would be desired to follow the ion dynamics and to be able to take time step sizes with $w_{pe} \Delta t  \gg 1$.

\section*{Acknowledgments}
 YC is supported by grants NSF DMS-1217563, DMS-1318186, AFOSR FA9550-12-1-0343 and the startup fund from Michigan State University. AJC is supported by AFOSR grants FA9550-11-1-0281, FA9550-12-1-0343 and FA9550-12-1-0455,  NSF grant DMS-1115709
and MSU foundation SPG grant RG100059.  We gratefully acknowledge the support from Michigan Center for Industrial and Applied Mathematics.

%
%\begin{figure}[!htbp]
%\subfigure[total particle number]{\includegraphics[width=0.45\textwidth]{landaustrangmass.eps}}
%\subfigure[Momentum]{\includegraphics[width=0.45\textwidth]{landaustrangmom.eps}}
%\\
%\subfigure[Enstrophy]{\includegraphics[width=0.45\textwidth]{landaustrangens.eps}}
%\subfigure[Total energy]{\includegraphics[width=0.45\textwidth]{landaustrangenergy.eps}}
%\label{landau_strang}
%\caption{Landau damping. $50\times 100$ mesh. $cfl = 10$.$ \epsilon_{tol}=1e-10$.}
%\end{figure}
%
%
%\begin{figure}[!htbp]
%\subfigure[total particle number]{\includegraphics[width=0.45\textwidth]{tsstrangmass.eps}}
%\subfigure[Momentum]{\includegraphics[width=0.45\textwidth]{tsstrangmom.eps}}
%\\
%\subfigure[Enstrophy]{\includegraphics[width=0.45\textwidth]{tsstrangens.eps}}
%\subfigure[Total energy]{\includegraphics[width=0.45\textwidth]{tsstrangenergy.eps}}
%\label{ts_strang}
%\caption{Two-stream instability. $50\times 100$ mesh. $cfl = 5$ for $S_3+P^2$ and $cfl = 10$ for $\textnormal{\bf Scheme-3F}+Q^3$.$ \epsilon_{tol}=1e-10$.}
%\end{figure}
%\begin{figure}[!htbp]
%\subfigure[total particle number]{\includegraphics[width=0.45\textwidth]{bumpstrangmass.eps}}
%\subfigure[Momentum]{\includegraphics[width=0.45\textwidth]{bumpstrangmom.eps}}
%\\
%\subfigure[Enstrophy]{\includegraphics[width=0.45\textwidth]{bumpstrangens.eps}}
%\subfigure[Total energy]{\includegraphics[width=0.45\textwidth]{bumpstrangenergy.eps}}
%\label{bump_strang}
%\caption{Bump-on-tail instability. $50\times 100$ mesh. $cfl = 10$ . $\epsilon_{tol}=1e-10$.}
%\end{figure}

\bibliographystyle{abbrv}
\bibliography{yingda,ref_cheng,refer,ref_cheng_plasma_2,xinghui,ts}

\end{document}